\documentclass[12pt, a4paper]{article}
\usepackage[left=1in,right=1in,top=1.0in,bottom=1in]{geometry}

\usepackage{amsmath,amssymb,amsthm,amsfonts,afterpage,hyperref}
\usepackage{amscd,subeqnarray}

\usepackage{graphicx}
\usepackage{authblk}
\usepackage{tikz}

\usepackage{array}
\usepackage{wrapfig}
\usepackage{color}
\usepackage{ragged2e}
\usepackage{subfig}
\usepackage{stmaryrd}
\usepackage{siunitx}
\usepackage{multirow}
\usepackage{xcolor,cancel}
\usepackage{boldfonts}
\usepackage{symbols}
\usepackage{footmisc}
\usepackage{float}
\usepackage{setspace}
\usepackage{bm}
\usepackage{booktabs}
\usepackage{tikz, xcolor}
\usepackage{soul}

\usepackage[ruled,vlined]{algorithm2e}
\SetKwInput{KwInput}{Input}                
\SetKwInput{KwOutput}{Output} 

\usepackage[sort&compress]{natbib}
\setcitestyle{numbers}

\usetikzlibrary{shapes,arrows, positioning,calc}	
\usetikzlibrary{arrows,snakes,backgrounds}

\newcommand{\bfeps}{\boldsymbol{\epsilon}}		
\newcommand{\bfeta}{\boldsymbol{\eta}}

\newcommand{\Sym}{\text{Sym}}   			%
\newcommand{\Div}{\text{Div}}     				%
\newcommand{\ddiv}{\text{div}}     				%
\newcommand{\tr}{\text{tr}}       				%

\DeclareMathAlphabet{\mathpzc}{OT1}{pzc}{m}{it}

\DeclareMathOperator*{\essinf}{ess\,inf}
\DeclareMathOperator*{\esssup}{ess\,sup}
\newcommand{\bfn}{\boldsymbol{n}}	
\newcommand{\bfh}{\boldsymbol{h}}	
\newcommand{\bfu}{\boldsymbol{u}}	
\newcommand{\bfb}{\boldsymbol{b}}	
\newcommand{\bfv}{\boldsymbol{v}}	
\newcommand{\bfw}{\boldsymbol{w}}
\newcommand{\bfF}{\boldsymbol{F}}	
\newcommand{\bfE}{\boldsymbol{E}}	
\newcommand{\bfC}{\boldsymbol{C}}
\newcommand{\bfB}{\boldsymbol{B}}	
\newcommand{\bfx}{\boldsymbol{x}}	
\newcommand{\bfX}{\boldsymbol{X}}	
\newcommand{\bfS}{\boldsymbol{S}}	
\newcommand{\bfT}{\boldsymbol{T}}		
\newcommand{\bfI}{\boldsymbol{I}}	 
\newcommand{\bfzero}{\boldsymbol{0}}	

\newcommand{\bfD}{\boldsymbol{D}}
\newcommand{\bfM}{\boldsymbol{M}}
\newcommand{\bfN}{\boldsymbol{N}}

\newcommand{\bfA}{\boldsymbol{A}}

\newcommand{\bff}{\boldsymbol{f}}	
\newcommand{\bfg}{\boldsymbol{g}}	

\newtheorem{definition}{Definition}
\newtheorem{theorem}{Theorem}[section]

\newtheorem{lemma}[theorem]{Lemma}
\newtheorem{remark}{Remark}

\newtheorem{Theorem}{Theorem}[section]

\providecommand{\keywords}[1]
{
  \small	
  \textbf{\textit{Keywords---}} #1
}

\title{\textsf{A finite element model for a coupled thermo-mechanical system: 
nonlinear strain-limiting thermoelastic body}}

\author[1,2]{Hyun C. Yoon\thanks{hyun.yoon@kigam.re.kr}}
\author[1]{Karthik K. Vasudeva\thanks{karthikkumarv.kk@gmail.com}}
\author[1]{S. M. Mallikarjunaiah\thanks{m.muddamallappa@tamucc.edu}\thanks{corresponding author}}
\affil[1]{Department of Mathematics \& Statistics,
Texas A\&M University - Corpus Christi
Corpus Christi, TX 78412-5825, USA}
\affil[2]{
Petroleum \& Marine Division,
Korea Institute of Geoscience and Mineral Resources,
124 Gwahak-ro,
Daejeon 34132, Republic of Korea}
\date{}

\begin{document}

\maketitle

\begin{abstract}
{We 
investigate} a specific finite element model to study the thermoelastic behavior of an elastic body within the context of nonlinear strain-limiting constitutive {relation}. 
{As a special subclass of implicit relations,} {the thermoelastic response of our interest is} 
{such that} stresses can be arbitrarily large, but strains remain small, especially in the neighborhood of crack-tips.  {Thus, the proposed model can be inherently consistent with the assumption of the 
small strain theory.}  In the present communication, we consider a two-dimensional coupled system {-- linear and quasilinear partial differential equations for temperature and displacements, respectively.} 
{Two distinct temperature} {distributions} {of the Dirichlet type} are considered for  
boundary condition, {and a standard finite element 
{method} of continuous Galerkin is employed to obtain the numerical solutions for the field variables.} 
For a domain with an edge-crack, we find that the near-tip strain growth {of our model} is much slower than the growth of stress{, which is the salient feature compared to the inconsistent results of the classical linearized description of the elastic body}.  
{Current study can provide} a theoretical and computational framework to develop physically meaningful models {and examine other coupled multiphysics} such as an evolution of complex network of cracks induced by thermal shocks.
\end{abstract}

\noindent \keywords{Thermoelasticity,  Strain-limiting,  Crack-tip singularity,  Nonlinear elasticity,  Finite element method.}

\section{Introduction}

{Modern physics and engineering problems in nature or  man-made system are invariantly multi-physical. For example, the objects experiencing high temperature differences (e.g., high-speed commercial and defense airplanes, space vehicles, nuclear reactors, and jet engines) cannot be formulated by a single equation but can only be described with every physics involved \cite{brown2008}.} 
 {Thus, computing these complex phenomena for various applications requires not only a rigorous algorithm to accurately solve the coupled problem (e.g., the thermal-hydrological-mechanical (THM) processes or the fluild-structure interaction (FSI)) but also} 
an {individual} mathematical and computational  
{model} that can provide a realistic assessment of mechanical or thermo-mechanical state of the body. 
{Regarding a robust 
model, it is structured on some key components.} 
{For instance, in designing materials that can withstand and work at sufficiently high temperatures, a  thermo-mechanical model within the continuum scale} {requires correct} input variables such as 
mechanical properties of the material {(e.g., the Lam$\acute{e}$ parameters, thermal conductivity)} and the external boundary conditions {(e.g., the purely mechanical or a combination of thermo-mechanical loading)}.  
{Ultimately,} an essential constituent of that model would be {the ``constitutive law'' that provides a deterministic {but also physically meaningful or realistic} relation among these input quantities. 
Given all these, the response of the material body can be characterized by {the output fields: displacement, strain, stress, temperature distribution, and heat flux.} 

{However,} some of the celebrated theories in applied mechanics erroneously predict unphysical results {that must be corrected to obtain a reliable solution.} 
{Concerning the heat transfer,  
one} eminent constitutive theory that {can} yield unreasonable results {comes from the theory of heat conduction}, 
{\textit{Fourier's law}} (also referred as \textit{Fourier's ``first'' law}). 
{As} 
given the status of a law {based on an approximation of a linear relationship between the heat flux and the temperature gradient, it yields some gratifying results in many instances.} {Nevertheless,} this widely-used theory cannot help predicting strangely that the temperature propagates with infinite velocity. 
{This fundamental issue is due to the Fourier's law assumption} that the heat flux {soars} 
spontaneously {right} after the temperature gradient applied, which in practice is not a realistic assumption since no physical process {can be} instantaneous in reality. Several corrections have been proposed to this law, {such as the seminal work} by Cattaneo \cite{cattaneo1958} and Vernotte \cite{vernotte1958}, popularly {known} as the CV model \cite{marin2011}. {Still} there are some known issues with the CV model{: for instance, the} necessity for the relaxation time 
as a modeling parameter and a wave-like motion of the temperature variable due to the hyperbolic nature of the constituent partial differential equation. There are several correction models proposed in the relevant literature to ensure the finite speed for the temperature propagation \cite{chester1963,lord1967generalized,chandrasekharaiah1986,chandrasekharaiah1998,agrawal2016,whitaker2013}.  {Numerous works including but not limited to the above} 
consider the elastic bodies described either by higher gradients or within the framework of Cauchy elasticity with additional thermal effects. 

{On the other hand, there is another notorious issue related to the singmarin2011ular value {near any concentrator (e.g., a crack-tip) in an elastic body} -- the crack-edge strain singularity {delineated} \cite{broberg1999,anderson2017fracture} by the linear elastic fracture mechanics (LEFM). {The issue can occur} when
modeling the crack-surface by using the classical zero-traction boundary condition.}  Several new constitutive theories have been proposed with augmentations to the classical LEFM model to eliminate the inconsistent prediction.  In many of these theories, it has long been contended that the inconsistent crack-edge strain singularity is due to the use of zero traction crack-face boundary condition which does not take into account of the work done near the fracture edge. Hence, many of the theories focused on correcting the mechanics of the crack-surface as a modified boundary condition for the mathematical boundary value problem, such as the notion of \textit{cohesive zone} introduced by Barenblatt \cite{barenblatt1962}, the \textit{surface elasticity} archetype proposed by Gurtin and Murdoch \cite{gurtin1975},  a paradigm of modeling crack-face as a \textit{diving surface} as introduced by Sendova and Walton \cite{sendova2010} and a regularization to {the surface} boundary condition for a stable finite element simulation as proposed in \cite{ferguson2015}. However, some of these modifications to crack-edge boundary condition do not always prevent the square-root strain singularity at crack-tip as predicated by the classical LEFM.

{In this contribution, we focus on studying the response of a thermoelastic body and 
{resolving} the singularity, which cannot be underemphasized, near the concentrator such as a crack-tip. Recently, Rajagopal and his co-authors \cite{rajagopal2007elasticity,rajagopal2007response,rajagopal2016,bustamante2020} made a compelling argument that one can arrive at the nonlinear models for {small strain based on infinitesimal theory} to describe {the elastic} behavior {of body, i.e.,} non-dissipative materials. More importantly, the traditional Cauchy and Green elastic models can be regarded as a subclass of models introduced by Rajagopal \cite{rajagopal2003}. An important feature of Rajagopal's models is that strain being limited below a certain value -- which can be set \textit{a priori} -- regardless of the stress value. Such a possibility cannot be obtained within the framework of classical linearized theory of elasticity. 
 Using such algebraic nonlinear models, several studies \cite{rajagopal2011modeling,gou2015modeling,Mallikarjunaiah2015,kulvait2013,kulvait2019,bustamante2010note, bustamante2011solutions,bustamante2015implicit,ortiz2012,bustamante2018nonlinear} have been able to revisit the classical problems of solid mechanics, and poromechanics \cite{fu2019generalized,fu2020constraint}.
In particular, the results {in} \cite{Mallikarjunaiah2015,kulvait2013,kulvait2019,ortiz2012,ortiz2014numerical,hyun-mms-2021} predict crack-tip stress singularity similar to the classical LEFM model for static crack problems.  This implies that the crack-tip is a singular energy sink hence one can utilize a local fracture criterion to study the evolution of cracks within the nonlinear strain-limiting relation \cite{hyun-mms-2021,lee2020nonlinear,yoon2021quasi}. 
{Thus, the objective of this study is to offer a reasonable rationale to use the nonlinear strain-limiting models and 
describe {a reliable} thermoelastic response {particularly around the crack-tip}.} To that end, we 
{take as a starting point} constructing the boundary value problem within nonlinear stress-strain relationship, {of which} one special case 
{reduces to} the classical linearized thermoelasticity constitutive model.} 
{For future work, this study can be a stepping stone to develop physically meaningful models for the crack evolution under thermo-mechanical loading such as by using the phase-field regularization \cite{bourdin2000numerical,lee2020nonlinear,yoon2021quasi}.} 

The organization of the paper is as follows. In Section 2, we provide a brief introduction about basic notations of  classical elasticity, {and} introduce the basic equations of a new class of elastic constitutive relationships, {which are} followed by the definition of constitutive relations for implicit thermoelasticity. We start with an implicit response relation between the Cauchy stress, the Cauchy-Green tensor, and absolute temperature to describe the response of a thermoelastic body {to derive the proposed model in this study}. In Section 3, we formulate a boundary value problem using a special subclass of isotropic relation with a fundamental assumption that the displacement gradients are infinitesimal. In the same section, we propose a numerical method based on the stable finite element discretization of the {coupled} system of {linear and quasilinear} partial differential equations {for temperature and mechanics, respectively}. We also present a weak formulation and a theorem on {the} existence and uniqueness of solution for the quasilinear partial differential equation. Numerical experiments and corresponding results are described in Section 4. Finally, in Section 5, concluding remarks and some ideas for the future work are made.

\section{Formulation of the nonlinear thermoelastic problem}
\subsection{Basic notations}
In spatial dimension $d=2$, we consider $\Omega \subset \mathbb{R}^d$ is an open bounded polygonal domain, which represents the reference configuration of the elastic material with boundary $\partial \Omega$ and an outward normal $\bfeta = \left(\eta_1, \ldots, \eta_d\right)$. The boundary $\partial \Omega$ be a $(d-1)$ dimensional Lipschitz manifold consists of two mutually disjoint measurable sets $\Gamma_N$ and  $\Gamma_D$, such that $\partial \Omega = \overline{\Gamma_N \cup \Gamma_D}$, and at least nonempty Dirichlet boundary $\Gamma_D \neq \emptyset$. We let $\bfx = \left(x_1, \ldots, x_d \right)$ denote a spatial point in the current configuration while $\bfX = \left(X_1, \ldots, X_d \right)$ denote the same point in the reference configuration. Let $\bfu \colon \Omega \to \mathbb{R}^d$ denote the displacement field of the material body and $\theta \colon \Omega \to \mathbb{R}$ is the temperature. In the rest of this paper, we use the usual notations of \textit{Lebesgue} and \textit{Sobolev spaces}. For example, we shall use $L^{p}(\Omega)$ as the space of all \textit{Lebesgue integrable functions} with $p\in[1, \infty)$. In particular when $p=2$, $L^{2}(\Omega)$ denote the space of all square integrable functions with $\left( v, \; w \right):= \int_{\Omega} v \,w \; d\bfx$ as the inner product and $\| v \| := \| v \|_{L^2(\Omega)} = \left(v, \, v \right)^{1/2}$ for the corresponding induced norm. Then, let $C^{m}(\Omega), \; m \in \mathbb{N}_0$ denote the linear space of continuous functions on $\Omega$. Further, let $H^{k}(\Omega)$ denote the classical \textit{Sobolev space} \cite{ciarlet2002finite,evans1998partial} and
\begin{align}
H^{k}(\Omega) &:= W^{k, 2} = \left\{ v \in L^{2}(\Omega) \; \colon \; \sum_{|\zeta| \leq k} D^{\zeta} v \in L^{2}(\Omega) \right\},  \\
H^1_0(\Omega) &:= \mbox{closure of} \; C_{0}^{\infty}(\Omega) \; \mbox{in} \; H^{1}(\Omega), \label{def-H01}
\end{align}
here $\zeta = \left\{\zeta_1, \, \zeta_2, \ldots, \zeta_d \right\}$ is a multi-index,  $|\zeta| = \sum_{i=1}^{d} \zeta_i$, and the derivatives are in the weak sense. Let $H^{-1}(\Omega)$ denote the dual space to $H^{1}(\Omega)$. We also define the following subspaces of $H^1(\Omega)$
\begin{subequations}\label{test-fun-space}
\begin{align}
V^{\bfu} &:= \left\{ \bfu \in \left( H^{1}(\Omega)\right)^d \colon \bfu=\bfzero \; \mbox{on} \; \Gamma_D^{\bfu}\right\}, \\
V^{\theta} &:= \left\{ \theta \in  H^{1}(\Omega) \colon \theta=0 \; \mbox{on} \; \Gamma_D^{\theta}\right\}.
\end{align}
\end{subequations}
{Let} $\Sym(\mathbb{R}^{d \times d})$ is the space of $d \times d$ symmetric tensors equipped with the inner product $\bfA \colon \bfB = \sum_{i, \, j=1}^d \, \bfA_{ij} \, \bfB_{ij}$ and for all $\bfA = (\bfA)_{ij}$ and $\bfB =(\bfB)_{ij}$ in $\Sym(\mathbb{R}^{d \times d})$, {the associated norm}  {$\| \bfA \| = \sqrt{\bfA \colon \bfA}$}. Further, let $\bfF \colon \Omega \to \mathbb{R}^{d \times d}$ denote the deformation gradient,  $\bfC \colon \Omega \to \mathbb{R}^{d \times d}$ denote the right Cauchy-Green stretch tensor, $\bfB \colon \Omega \to \mathbb{R}^{d \times d}$ denote the left Cauchy-Green stretch tensor, $\bfE \colon \Omega \to \mathbb{R}^{d \times d}$ denote the Lagrange strain, $\bfeps \colon \Omega \to  \Sym(\mathbb{R}^{d \times d})$ denote the symmetric linearized strain tensor, respectively defined as:
\begin{align}
\bfF &:= \nabla_r \bfx = \bfI + \nabla \bfu, \\
 \bfB &:=\bfF\bfF^{\mathrm{T}}, \; \bfC :=\bfF^{\mathrm{T}}\bfF,  \; \bfE := \dfrac{1}{2} \left( \bfC - \bfI \right),\\
\bfeps(\bfu) &:= \dfrac{1}{2} \left( \nabla \bfu + \nabla \bfu^{\mathrm{T}}\right),
\end{align}
where $\left( \cdot \right)^{\mathrm{T}}$ denotes the \textit{transpose} operator for the second-order tensors, $\bfI$ is the $d$-dimensional indentity tensor, $\nabla_r$ and $\nabla$ are the gradient operators in the reference and current configuration, respectively.
Under the standard assumption of linearized elasticity that 
\begin{equation}\label{small_grad}
\max_{\bfX \in \Omega} \| \nabla_X \bfu \|  \ll \mathcal{O}(\delta), \quad \delta \leq 1, 
\end{equation}
whence there is no distinction between reference and current configuration. Let $\bfT \colon \Omega \to  \Sym(\mathbb{R}^{d \times d})$ be the \textit{Cauchy stress tensor} in the current configuration and it satisfies the equation of motion 
\begin{equation} 
\rho \, {\ddot{\bfu}} = \ddiv \, \bfT + \rho \, \bfb,
\end{equation}
where $\rho$ is the density of the body and $\bfb \colon \Omega \to \mathbb{R}^d$ is the body force in the current configuration and the notation $\dot{\left(  \;\; \right)}$ denotes the time derivative. The first and second \textit{Piola-Kirchhoff stress tensor tensors},  $\bfS \colon \Omega \to \mathbb{R}^{d \times d}$ and $\overline{\bfS} \colon \Omega \to \mathbb{R}^{d \times d}$, are defined by 
\begin{equation}
\bfS =  \bfT \bfF^{-\mathrm{T}} \, J, \quad \overline{\bfS} :=\bfF^{-1} \bfS, 
\end{equation}
where $J = \det(\bfF)$ and it is assumed to be positive for realistic deformations. 

\subsection{A new class of elastic solids}
In the present communication, we are interested in the studying a subclass of general implicit relations \cite{rajagopal2003,rajagopal2007elasticity,rajagopal2011non,rajagopal2011conspectus,rajagopal2014nonlinear,rajagopal2007response} for thermoelastic body. The classical Cauchy elasticity theory {describes} the Cauchy stress based on an assertion that 
\begin{equation}
\bfT = \mathcal{G}(\bfF).
\end{equation}
In case of isotropic and homogeneous elastic body, the stress tensor is given by 
\begin{equation}
\bfT = \delta_1 \, \bfI + \delta_2 \, \bfB +  \delta_3 \, \bfB^2,
\end{equation}
where $\delta_i, \; i=1, \,2, 3$ depend on $\rho, \, \tr\bfB, \, \tr\bfB^2, \, \tr\bfB^3$. Rajagopal \cite{rajagopal2007elasticity} generalized the constitutive of relationships of Cauchy elasticity by introducing  implicit elastic response relations as 
\begin{equation}\label{implicit-1}
0= \Xi(\bfT,\,\bfB).
\end{equation}
Further, Rajagopal \cite{rajagopal2007elasticity} also considered the special subclass of \eqref{implicit-1}:
\begin{equation}\label{SL1}
\bfB := \mathcal{F}( \bfT).
\end{equation}
For the above class of models \eqref{SL1}, if there exists a constant $M >0$ such that 
\begin{equation} 
\sup_{\bfT \in \Sym} |\mathcal{F}( \bfT)| \leq M,
\end{equation}
and such relations are referred as the \textit{strain-limiting constitutive theories} \cite{MalliPhD2015,Mallikarjunaiah2015}. Under an assumption of classical linear elasticity as in \eqref{small_grad}, the relationship \eqref{SL1} in the infinitesimal regime reduces to 
\begin{equation}\label{SL2}
\bfeps =  \mathcal{F}( \bfT),
\end{equation}
where the response function $\mathcal{F}$ is a mapping defined as:
\begin{equation}
\mathcal{F} \colon \Sym(\mathbb{R}^{d \times d}) \mapsto \Sym(\mathbb{R}^{d \times d}), \quad \mathcal{F}(\bfzero) = \bfzero.
\end{equation}
In this contribution, we have used the following form of the tensor valued function $\mathcal{F}(\cdot)$ which is a nonlinear function of the Cauchy stress tensor as in the constitutive relationship \eqref{SL2} is
\begin{equation}\label{def-F}
\mathcal{F}(\bfT) = \dfrac{ \mathbb{K}[\bfT]}{ \left( 1 + \beta^{a} | \mathbb{K}^{1/2}[\bfT] |^{a}                   \right)^{1/a}}, \quad \beta, \; a >0,
\end{equation}
which is uniformly bounded with $M = 1/\beta$
and the symbol $\mathbb{K}[\cdot]$ is a fourth-order, linearized compliance tensor which is the inverse to the linearized elasticity tensor $\mathbb{E}[\cdot]$.  For {an} isotropic and homogeneous elastic material {with linearization}, these tensors are defined as: 
\begin{subequations}\label{K-E-def}
\begin{align}
\mathbb{E}[\bfeps] &:= 2 \mu\bfeps + \lambda \, \tr(\bfeps) \, \bfI, \\
\mathbb{K}[\bfT] &:= \dfrac{1}{2 \, \mu} \, \bfT - \dfrac{1}{2 \, \mu \left( \lambda + \frac{2}{d} \mu \right)d} \, \tr(\bfT) \, \bfI,
\end{align}
\end{subequations}
in which $\mu$ and $\lambda$ are Lam\'{e} parameters given by:
\begin{equation}
\mu = \dfrac{E}{2 (1 + \nu)}, \quad \lambda = \frac{E \nu}{(1+\nu)(1-2\nu)},
\end{equation}
where $E$ denotes Young's modulus and $\nu$ denotes Poisson's ratio. Our interest in this paper is to study a homogenous material behavior under thermo-mechanical loading, which implies that $\mu$ and $\lambda$ are constants. It is important to note that if the parameter $\beta$ tends to zero in the above equation \eqref{def-F}, then the limiting strain model \eqref{SL2} turns into the classical linearized constitutive relationship. See \cite{MalliPhD2015,Mallikarjunaiah2015}, for the detailed explanation about the qualitative properties of these fourth-order tensors and the construction of the square-root of these higher-order tensors. Moreover, both $\mathbb{K}[\cdot]$ and  $\mathbb{E}[\cdot]$ are  constituted to be positive definite, self-adjoint operators on $\Sym$, and these two tensors satisfy
\begin{equation}
\mathbb{E}[\mathbb{K}[\bfT]] = \bfT, \quad \mathbb{K}[\mathbb{E}[\bfeps]] = \bfeps. 
\end{equation}
Based on \cite{itou2018states}, we list some important properties of the tensor valued function \eqref{def-F} as a lemma.
\begin{lemma}\label{lemma1}
Let the constants $M>0$ and $C_{\beta, \, a}$ exist, such that for all $\bfT_1, \, \bfT_2 \in \text{Sym}(\mathbb{R}^{d \times d})$ the following holds:
\begin{itemize}
\item Uniform boundedness: $\quad \| \mathcal{F} (\bfT) \|  \leq M$,
\item Monotonicity: $\quad  \left( \mathcal{F} (\bfT_1) - \mathcal{F} (\bfT_2) \right) \colon \left( \bfT_1 - \bfT_2 \right)  \geq 0$, 
\item Lipschitz continuity: $\quad \dfrac{\left( \mathcal{F} (\bfT_1) - \mathcal{F} (\bfT_2) \right) \colon \left( \bfT_1 - \bfT_2 \right) }{\| \bfT_1 - \bfT_2 \|^2} \leq C_{\beta, \, a}$,\\
{where} the constant $C_{\beta, \, a}$ depends on the modeling parameters $\beta$ and $a$.
\end{itemize}
\end{lemma}
The conditions in the \textbf{Lemma 2.1.} 
guarantee that $\mathcal{F}(\bfT)$ is a monotone operator for all $\bfT \in \Sym(\mathbb{R}^{d \times d})$. \\

The nonlinear  strain-limiting constitutive relationship \eqref{SL2} has been studied by several researchers to revisit various classical problems of solid mechanics \cite{kulvait2013,kulvait2019,bulivcek2014elastic,ortiz2014numerical,ortiz2012,kannan2014unsteady,kambapalli2014circumferential,fu2019generalized,fu2020constraint} including the evolution of quasi-static cracks \cite{yoon2021quasi,lee2020nonlinear}. The physical meaning of the response relation \eqref{SL2} is explained in \cite{rajagopal2014nonlinear}, and a systematic derivation of thermodynamically consistent constitutive theory within the framework of fully implicit theory of elasticity is established in \cite{rajagopal2003} and also in \cite{bridges2015implicit}.

\subsection{Nonlinear constitutive relations for thermoelastic body}
Our objective in the current contribution is to study the temperature distribution and the state of crack-tip stress-strain in a thermoelastic material body which is described by a special subclass of implicit relations \cite{rajagopal2003,rajagopal2014nonlinear}. For the derivation of the constitutive relation for the thermoelastic body, we closely follow the work in \cite{bustamante2017implicit}. To derive the nonlinear constitutive relation among the field variables, we assume that the material under consideration is an isotropic, homogeneous solid that occupies the entire region defined by $\Omega$. Bustamante and Rajagopal \cite{bustamante2017implicit} showed that a general class of constitutive relations that could be used to describe the response of thermoelastic body, specifically the class of implicit models of the form
\begin{equation}\label{imp-thermo-mod}
\mathfrak{F}(\overline{\bfS}, \, \bfE, \, \theta) =0,
\end{equation}
where $\theta \colon \Omega \to \mathbb{R}$ is the absolute temperature and $\mathfrak{F}$ is the second-order tensor valued relation. A special sub-class of the above fully implicit model \eqref{imp-thermo-mod} is given by 
\begin{equation}\label{exp:model}
\overline{\bfS} = \mathfrak{D}(\bfE, \, \theta). 
\end{equation}
Let $\bfh_r \colon \Omega \to \mathbb{R}^d$ be the heat flux in the reference configuration, $\varepsilon \colon \Omega \to \mathbb{R}$ denotes the internal energy, $\rho_r \in \mathbb{R}_{+}$ is the density of the material in the reference configuration, $w \colon \Omega \to \mathbb{R}^d,$ with $w = \tr (\bfS \, \dot{\bfE} )$ is the rate of work,
and $\Upsilon \colon \Omega \to \mathbb{R}$ is the rate of internal heat generation. Then, the the first law of thermodynamics in the reference configuration is given by \cite{callen1998thermodynamics,bustamante2017implicit}
\begin{equation}
 \dot{\varepsilon} = \frac{1}{\rho_r} (w + \Div\, \bfh_r )+  \Upsilon.
\end{equation}
The dissipation $\Pi \colon \Omega \to \mathbb{R}$ of the material is written as 
\begin{equation}
\Pi = \theta \, \dot{\eta} - \dot{\varepsilon}  + \frac{1}{\rho_r} \, w,
\end{equation}
in which the term $\eta \colon \Omega \to \mathbb{R}$ is the specific entropy density (per unit mass), and moreover the dissipation must satisfy the inequality
\begin{equation}
\Pi \geq 0.
\end{equation}
If $\mathfrak{F}$ is an isotropic function, then it follows that the implicit model defined in \eqref{imp-thermo-mod} takes the following most general form
\begin{align}
\alpha_0 \, \bfI &+ \alpha_1 \, \bfS + \alpha_2 \, \bfS^2 + \alpha_3 \, \bfE + \alpha_4 \, \bfE^2 + \alpha_5 \, ( \bfE\bfS + \bfS\bfE ) + \alpha_6 \, ( \bfE^2\bfS + \bfS\bfE^2 )  \notag \\
&+ \alpha_7 \, ( \bfS^2\bfE + \bfE\bfS^2 ) + \alpha_8 \, ( \bfS^2\bfE^2 + \bfE^2\bfS^2 ) = \bfzero,
\end{align}
where the material moduli $\alpha_i, \; i=0, \cdots, 8$ depend on the following invariants \cite{bustamante2017implicit}
\begin{align}
\Big\{& \tr(\bfS), \, \dfrac{1}{2} \tr(\bfS^2), \, \dfrac{1}{3} \tr(\bfS^3), \, \tr \bfE, \, \dfrac{1}{2} \tr(\bfE^2), \, \dfrac{1}{3} \tr(\bfE^3), \, \tr (\bfS \bfE), \, \tr (\bfS^2 \bfE), \notag \\
&   \tr (\bfS \bfE^2), \, \tr (\bfS^2 \bfE^2), \, \theta \Big\}.
\end{align}
Now, let us consider a special sub-class of the explicit model \eqref{exp:model} which is a foundation for the  same approximation that is used in the classical linearized elasticity, precisely that the displacement gradient is small, i.e., like in \eqref{small_grad}. Therefore, we have following approximations as
\begin{equation}
\nabla_r \bfu \approx  \nabla \bfu, \quad \bfE \approx \bfeps, \quad \bfS \approx \bfT, \quad \nabla_r \theta \approx  \nabla \theta, \quad \bfh_r \approx \bfh,
\end{equation}
in which the symbol $\nabla$ is the gradient operator in the current configuration. After a series of assumptions on the invertibility of the constitutive class of the function and neglecting the higher order $\delta$ terms, we arrive at the most general form of a sub-class of relationship for the field variables as:
\begin{equation}\label{exp-bfeps-1}
\bfeps = \mathcal{F}(\bfT, \, \theta).
\end{equation}
Notice that the model for the isotropic homogeneous elastic materials is a special subclass of the above constitutive class \eqref{exp-bfeps-1} and if $\mathcal{F}(\cdot)$ is an isotropic function, then {a} subclass is given by
\begin{equation}\label{exp-bfeps-2}
\bfeps = \hat{\alpha}_0 \, \bfI + \hat{\alpha}_1 \, \bfT +  \hat{\alpha}_1 \, \bfT^2,
\end{equation}
where $\hat{\alpha}_i, \; i=0, \, 1, \, 2$ depend upon the invariants 
\begin{equation}
\left\{ \tr(\bfT), \, \dfrac{1}{2} \tr(\bfT^2), \, \dfrac{1}{3} \tr(\bfT^3), \, \theta  \right\}.
\end{equation}
The constitutive class as in \eqref{exp-bfeps-2} can be perceived as a counterpart to the classical model derived from the starting point of Cauchy or Green elasticity. An obvious fact in the models of type \eqref{exp-bfeps-2} is that the roles of the stress and strain are reversed. More importantly, the generalized model of the type $\bfeps = \mathcal{F}(\bfT)$ provides the linearized strain as a nonlinear function of the stress, and such a possibility will pave a way that one could have bounded strains when stress blows up. Hence, the relationship \eqref{exp-bfeps-2} is an important prospect to describe the behavior of brittle materials including redefining the nature of stress-strain singularity near the concentrators such as cracks and fractures \cite{rajagopal2011modeling,kulvait2013,ortiz2014numerical,gou2015modeling,Mallikarjunaiah2015,MalliPhD2015,lee2020nonlinear,yoon2021quasi}. Next, a legitimate question is to check about the invertibility of the class {of the} nonlinear relation \eqref{SL2}. It was shown in \cite{mai2015strong,Mallikarjunaiah2015} that if $\mathcal{F}(\cdot)$ is a monotone operator, hence it is invertible with the inverse map as $\mathcal{F}^{-1}(\cdot)$. The constitutive relation equivalent to Green elasticity or hyperelastic formulation, as introduced in \cite{Mallikarjunaiah2015}, is  given by
\begin{equation}\label{T-as-inverse}
\bfT= \mathcal{F}^{-1}(\bfeps).
\end{equation}
It is important to emphasize that the model of type \eqref{SL2} or \eqref{T-as-inverse} converges to the classical linear stress-strain relation as $\beta \to 0$, {which is another motivation to study the model}. 
{As a final form,} the linearized constitutive equations for thermoelastic bodies can be written
\begin{equation}\label{model-T}
\bfT := \Psi(|\mathbb{E}^{1/2} [\bfeps]|) \, \mathbb{E}[\bfeps] - \alpha \, \theta,
\end{equation}
and the function $\Psi(\cdot)$ can be chosen with a special form as:
\begin{equation}\label{eq:Psi}
\Psi(r) = \dfrac{1}{ \left(1 - (\beta \, r)^a          \right)^{1/a}} \quad \mbox{with} \quad \beta \geq 0, \; a >0. 
\end{equation}
In the present numerical study, we utilize the above form of the stress-strain relation to study the response of the 
\textit{{strain-limiting}} thermoelastic body under the combination of thermo-mechanical loading. 

\begin{remark}
One can also derive linearized constitutive equations for thermoelastic body, as described in \cite{bustamante2017implicit}, in which one may not obtain an explicit description for  stress as the one in \eqref{model-T}. In \cite{bustamante2017implicit}, a linearized implicit constitutive relationship is obtained of the form:
\begin{equation}\label{model-eps}
\bfeps = \omega_0 \, \tr(\bfT) \, \bfI + \omega_1 \, \bfT + \omega_2 \, (\theta -1) \, \bfI, 
\end{equation}
where $\omega_0, \, \omega_1, \,\omega_2$ are constants. Models of the type \eqref{model-T} or \eqref{model-eps} are very important in their own rights, as these could be compelling candidates to describe the response of brittle materials such as rock, gum metal, ceramics, glass, and more specifically in the study of evolution of cracks and fractures in brittle bodies. 
\end{remark}

\section{Boundary value problem: Thermoelastic state}
{In many industrial applications}, the expansion and contraction of materials due to thermo-mechanical loading  are of great importance. 
{Generally speaking,} the thermo-mechanical loading occurs simultaneously, therefore the material response will be governed by the coupled displacement and temperature fields. Hence, both field variables have to be modeled as a simultaneous system of equations taking into account of highly nonlinear interaction. In practice, deriving such a governing system is far from obvious. In this section, we formulate a boundary value problem 
of a {loosely} coupled linear-quasilinear system of partial differential equations 
for temperature and displacements pair 
$(\bfu = ({\bfu_x,\, \bfu_y}), \; \theta)$ {in 2D} which {describes} the thermoelastic behavior of the body. {The loosely coupled system here is
implying an one-way coupling system where the heat energy affects the mechanics, not the other way around.}

\subsection{Mathematical model and \textit{Newton's method} for linearization}
 The elastic material body under consideration is homogeneous, thermally conductive, isotropic at zero initial reference temperature and initially unstrained and unstressed. Let $\Omega$ be an open, bounded, Lipschitz, and connected domain with the boundary $\partial \Omega$ consisting of two smooth disjoint parts $\Gamma_N^{\#}$ and $\Gamma_D^{\#}$ such that $\partial \Omega = \Gamma_N^{\#} \cup \Gamma_D^{\#}$, here the symbol $\#$ can be either $\bfu$ or $\theta$. To describe the thermoelastic state of the material under investigation, we consider the following basic equations for the generalized theory of thermoelasticity in the absence of the body force {and variation of internal energy}:
\begin{subequations}
\begin{align}
- \nabla \cdot \left(k \, \nabla \theta \right) &= g \quad \text{in} \quad \Omega, \label{eq1} \\
 - \nabla \cdot \bfT + \alpha \nabla \theta &= \boldsymbol{0} \quad \text{in} \quad \Omega,  \label{eq2}
\end{align}\label{eq:gov}
\end{subequations}
where $k \in \mathbb{R}$ is the heat conduction coefficient,  $\alpha = \alpha_T \, (3\lambda + 2 \mu)$, $\alpha_T$ is linear thermal expansion coefficient, $\lambda$ and $\mu$ are the Lam$\acute{e}$ parameters, and $g \colon \Omega \to \mathbb{R}$ is the heat source term. To model the stress response of the material, we make use of \textit{strain-limiting} theory as in \cite{rajagopal2007elasticity,Mallikarjunaiah2015}:
\begin{equation}
\bfeps = \frac{\mathbb{K}[\bfT]}{\left( 1 + \beta^a \, | \mathbb{K}^{1/2}[\bfT] |^a  \right)^{1/a}}. \label{eq3}
\end{equation}
We note that the above constitutive relationship is invertible \cite{mai2015monotonicity,mai2015strong,Mallikarjunaiah2015}, and the inverted relation written for the Cauchy stress as:
\begin{equation}\label{def-T}
 \bfT=\frac{\mathbb{E}[\bfeps]}{(1-\beta^a |\mathbb{E}^{1/2}[\bfeps]|^{a})^{1/a}},
\end{equation}
where the fourth-order tensors $\mathbb{K}[\cdot]$ and $\mathbb{E}[\cdot]$ are defined in \eqref{K-E-def}. Using the above inverted constitutive relationship \eqref{def-T} in \eqref{eq2}, we obtain the following {updated equations for the model}: 
\begin{subequations}\label{system1}
\begin{align}
- \nabla \cdot \left(k \, \nabla \theta \right) &= g \quad \text{in} \quad \Omega, \label{eq1-updated} \\
    -\nabla \cdot \Bigg[\frac{\mathbb{E}[\bfeps]}{(1-\beta^a |\mathbb{E}^{1/2}[\bfeps]|^{a})^{1/a}}\Bigg] + \alpha \nabla \theta &= \boldsymbol{0}\quad \text{in} \quad \Omega. \label{eq2-updated}
\end{align}
\end{subequations}
The above system of partial differential equations (PDEs) need to be supplemented by appropriate boundary conditions for both the unknown variables. The following are the boundary conditions: 
\begin{subequations}\label{bcs-v1}
\begin{align}
\bfT  \bfn &= \textbf{0}, \; \bfx \in \Gamma_N^{\bfu}, \quad  \bfu = \bfu_0, \; \bfx \in \Gamma_D^{\bfu}, \\
-k \, \bfn \cdot \nabla \theta &=0, \;  \bfx \in \Gamma_N^{\theta}, \quad  \theta = \theta_0, \; \bfx \in \Gamma_D^{\theta}, 
\end{align}
\end{subequations}
where $\bfn$ is the {outward} unit normal vector to the surface, and the boundary values $\bfu_0 \in (L^2(\Omega))^d$, $\theta_0 \in L^2(\Omega)$.  In the system of equations \eqref{system1}, the first equation in \eqref{eq1-updated} is a classical Laplacian for the temperature variable if the heat conduction coefficient is constant on the entire domain, and the second equation  in \eqref{eq2-updated} is a quasi-elliptic PDE. The solution algorithm we use is  \textit{Newton's method} at the PDE level, which results in a sequence of linear PDEs that can be discretized within the framework of standard finite element method. A critical issue with the Newton's method is about achieving the second-order convergence rate, and such a rate can be achieved only if the initial guess, for the solution variable, is \textit{sufficiently} close to the exact solution. The results presented in this paper are purely computational, although we have examined a standard routine of convergence test with a manufactured solution.  

Consider the quasilinear part of the second equation in \eqref{system1}
\begin{equation}\label{eq-L}
\mathcal{L}(\bfeps(\bfu)) := \Bigg[ \frac{\mathbb{E}[\bfeps(\bfu)]}{(1-\beta^a |\mathbb{E}^{1/2}[\bfeps(\bfu)]|^{a})^{1/a}}\Bigg],
\end{equation}
with 
\begin{subequations}\label{eq-E}
\begin{align}
\mathbb{E}[\bfeps(\bfu)] &:= 2 \mu\bfeps(\bfu) + \lambda \, \tr(\bfeps(\bfu)) \, \bfI \\
&=\mu \left( \nabla \bfu + \nabla \bfu^{\mathrm{T}}\right) + \lambda (\nabla \cdot \bfu) \, \bfI,
\end{align}
\end{subequations}
and 
\begin{subequations}\label{eq-SquareRootE}
\begin{align}
| \mathbb{E}^{1/2} [\bfeps(\bfu)]^2 &:=  \mathbb{E}^{1/2} [\bfeps(\bfu)] \colon \mathbb{E}^{1/2} [\bfeps(\bfu)]  \\
&= \bfeps(\bfu) \colon \mathbb{E} [\bfeps(\bfu)] \\
& = \frac{1}{2} \left( \nabla \bfu + \nabla \bfu^{\mathrm{T}}\right) \colon \left[ \mu  \left( \nabla \bfu + \nabla \bfu^{\mathrm{T}}\right)  \right] + \lambda \left( \nabla \cdot \bfu \right)^2.
\end{align}
\end{subequations}
The following lemma {is} for the tensor-valued mapping $\mathcal{L}(\cdot)$ as {assumed} and given in \eqref{eq-L}. 

\begin{lemma}
Let $m >0$, and $\mathfrak{{G}}:=\left\{ \bfD \in \Sym(\mathbb{R}^{d \times d}) \colon \| \bfD \| \leq m \right\}$, then consider the mapping 
\begin{equation}
\bfD \in \mathfrak{{G}} \mapsto \mathcal{L}(\bfD) := \dfrac{\bfD}{ \left( 1 - \beta^a \| \bfD \|^a       \right)^{1/a}} \in \Sym(\mathbb{R}^{d \times d}).
\end{equation}
Then, for each $\bfD, \, \overline{\bfD} \in \mathfrak{G}$, we have 
\begin{align}
|\mathcal{L}(\bfD) -  \mathcal{L}(\overline{\bfD}) | &\leq  \dfrac{ \widehat{C}_{\beta, \, \alpha}  \left( | \bfD - \overline{\bfD}  |  \right)}{\left(1 - \beta^{\alpha} (|\bfD| + | \overline {\bfD}|)^{\alpha}               \right)^{1/\alpha}}, \\
\left( \mathcal{L}(\bfD) -  \mathcal{L}(\overline{\bfD}) \right) \colon \left( \bfD - \overline{\bfD} \right) &\geq  \| \bfD - \overline{\bfD}  \|^2.
\end{align} 
The above {conditions 
imply} that $\mathcal{L}(\bfeps(\bfu))$ as defined in \eqref{eq-L} is a monotone operator. 
\end{lemma}
\begin{definition} (Uniformly Elliptic Operator)\\
The operator $\mathcal{L}(\cdot)$ is said to be uniformly elliptic if and only if for some constants $c_1, \; c_2 >0$, known as ellipticity constants, we have
\begin{equation}
c_1 \, \sup_{\| \bfv \| =1} \, \|\bfM \, \bfv\| \leq \mathcal{L}(\bfM + \bfN) -  \mathcal{L}(\bfN) \leq c_2 \, \sup_{\| \bfv \| =1}  \, \|\bfN \, \bfv\| \quad \forall \, \bfM, \; \bfN \in \Sym(\mathbb{R}^{d \times d}).
\end{equation}
\end{definition}
If $\mathcal{L}(\cdot)$ is \textit{Fre$\acute{e}$chet differentiable} and for some $\bfM \in \Sym(\mathbb{R}^{d \times d})$, then the derivative of $\mathcal{L}(\cdot)$ in the direction of $\bfN$ is given by
\begin{equation}
D\mathcal{L}(\bfM) \bfN = \mathcal{L}(\bfM) \colon \bfN.
\end{equation}
Suppose $\mathcal{L}(\cdot)$ is differentiable, then it satisfies the \textit{ellipticity condition} \cite{evans1998partial}, i.e., there exists $c_0 >0$ such that for each $\bfM \in \mathfrak{{G}}$ (where $\mathfrak{{G}}$ as defined in the above \textbf{Lemma 3.1.})
\begin{equation} 
\bfv \cdot \mathcal{L}(\bfM) \bfv \geq c_0 \, \| \bfv \|^2 \quad \forall \, \bfv \in \mathbb{R}^d,
\end{equation}
where the constant $c_0$ is independent of $\bfM$.\\

\noindent\textbf{Linearization using Newton's method: }The smoothness assumption that we made for the Sobolev space \eqref{def-H01} and the quasilinearity of the operator in \eqref{eq-L} bestow us  to apply \textit{Newton's method} to obtain the linearized version of the PDE for \eqref{eq2-updated}.  \\

Given the initial guess $\bfu^0 \in \left(C^2(\Omega)\right)^d$, with $\mathbb{E}[\bfeps(\bfu^0)] \in \mathcal{C} \subseteq \Sym(\mathbb{R}^{d \times d})$, for each $n \in \mathbb{N}_0$, find $\bfu^{n+1} \in  \left(C^2(\Omega)\right)^d$ with $\mathbb{E}[\bfeps(\bfu^{n+1})] \in \mathcal{C}$ such that 
\begin{equation}
D \mathcal{L}(\bfeps(\bfu^n)) \left( \bfu^{n+1} - \bfu^n \right) = - \mathcal{L}(\bfeps(\bfu^n)),
\end{equation}
where $D\mathcal{L}(\bfeps(\bfu))$ denotes \textit{Fre$\acute{e}$chet derivative}, and is defined as
\begin{equation}\label{FreDerivative}
D \mathcal{L}(\bfeps(\bfu))  \bfv = \lim_{\xi \to 0} \dfrac{\mathcal{L}(\bfeps(\bfu) + \xi \, \bfeps( \bfv))- \mathcal{L}(\bfeps(\bfu))}{\xi},
\end{equation}
for each $\bfv \in \left(C^2(\Omega)\right)^d$. Using \eqref{eq-L} in the definition of the \textit{Fre$\acute{e}$chet derivative} as in \eqref{FreDerivative}, we obtain
\begin{equation}\label{frechet}
 D \mathcal{L}(\bfeps(\bfu)) \delta  \bfu  =  \Bigg[
\frac{\mu \left({ \nabla \delta\bfu^n + \left(\nabla  \delta\bfu^n \right)^T}\right)  + \lambda \, (\nabla \cdot  \delta \bfu^n)\,  \bfI }{\left( 1 - \left( \beta \, | \mathbb{E}^{1/2} \left[ \bfeps \right]|\right)^{a}\right)^{1/a}}
\\ + 
\frac{\beta^{a} \,   \Theta_{1}\{  \bfu\} \, \Theta_{2}\{  \bfu^n,  \delta\bfu^n \}\,  \mathbb{E}[ \bfeps ]}{\left( 1 - \beta^{a} | \mathbb{E}^{1/2} \left[  \bfeps \right]|^{a}\right)^{1 + 1/a}} \Bigg]\, ,
\end{equation}
{where the upper script $n$ is the iteration number for the \textit{Newton's method}, and} the functions $\left| \mathbb{E}^{1/2} \left[\bfeps   \right]\right|$, $\Theta_{1}\{ \bfu^n \}$,  $\Theta_{2}\{ \bfu^n,  \delta  \bfu^n  \}$ are defined as:
\begin{align}
\left| \mathbb{E}^{1/2} \left[  \bfeps   \right] \right|^{2} &= \mathbb{E}^{1/2}[ \bfeps ] \colon \mathbb{E}^{1/2} \left[  \bfeps  \right]  \notag  \\
&= \bfeps \colon \mathbb{E}^{1/2}[\mathbb{E}^{1/2} \left[  \bfeps   \right] ] \notag  \\
&= \bfeps  \colon \mathbb{E} \left[  \bfeps   \right] \notag  \\
& = \left(\frac{ \nabla  \bfu^n   + (\nabla  \bfu^n){^{T}} }{2} \right)\,  \colon \mu\left({\nabla  \bfu^n   + (\nabla  \bfu^n){^{T}}}\right)   +  \lambda\, \left( \nabla \cdot   \bfu^n   \right)^{2},\label{eqn_temp}  \\
\Theta_{1}\{ \bfu^n  \} &:= \left| \mathbb{E}^{1/2} \left[  \bfeps   \right] \right|^{a-2},  \\
\Theta_{2}\{ \bfu^n,  \delta  \bfu^n  \} &:=\frac{1}{2} \left(   \left| \mathbb{E}^{1/2} \left[  \bfeps   \right] \right| \right)^{\prime}  \notag \\
& = \left(\frac{ \nabla  \bfu^n   + (\nabla  \bfu^n){^{T}} }{2} \right) \colon \mu \left({\nabla \delta \bfu^n   + (\nabla \delta \bfu^n)^{T}}\right)   +  \lambda \, \left( \nabla \cdot   \bfu^n   \right) \left( \nabla \cdot  \delta  \bfu^n    \right).
\end{align}
Then combining the equations \eqref{frechet} and \eqref{eq2-updated}, we obtain the updated model equations:
\begin{subequations}
\begin{align}
&- \nabla \cdot \left(k \, \nabla \theta \right) = g,  \label{strong-form-1} \\
&-\nabla \cdot \Bigg[
\frac{\mu \left({ \nabla \delta\bfu^n + \left(\nabla  \delta\bfu^n \right)^T}\right)  + \lambda \, (\nabla \cdot  \delta \bfu^n)\,  \bfI }{\left( 1 - \left( \beta \, | \mathbb{E}^{1/2} \left[ \bfeps \right]|\right)^{a}\right)^{1/a}}  \notag \\
&+ 
\frac{\beta^{a} \,   \Theta_{1}\{  \bfu\} \, \Theta_{2}\{  \bfu^n,  \delta\bfu^n \}\,  \mathbb{E}[ \bfeps ]}{\left( 1 - \beta^{a} | \mathbb{E}^{1/2} \left[  \bfeps \right]|^{a}\right)^{1 + 1/a}} \Bigg] = - \alpha \nabla \theta - \nabla \cdot \Bigg[\frac{\mathbb{E}[\bfeps(\bfu^n)]}{(1-\beta^a |\mathbb{E}^{1/2}[\bfeps(\bfu^n)]|^{a})^{1/a}}\Bigg]. \label{strong-form-2}
\end{align}
\end{subequations}

\begin{remark}
{Note in \eqref{eq1} particularly for the heat flux that we assume the classical heat transfer constitutive relation in this study. Since the heat conduction coefficient $k$ is assumed to be constant over a domain, the well-known unphysical computation of the infinite speed of heat propagation aforementioned can be drawn from the model. To ovecome the issue, various models have been proposed starting from the relaxation time term \cite{cattaneo1958,vernotte1958,marin2011} that makes the parabolic equation to the hyperbolic equation for the constitutive relation, along with the dependent, thus non-constant, heat conduction coefficient. In this paper, we focus on the nonlinear relation of stress and strain only from the mechanical side, and the simplified linear heat transfer without any variation of internal energy considered. Further, since we consider the thermal effect being loosely coupled to the mechanics, an explicit form of \eqref{model-T} can be a simplified and reasonable choice out of an implicit form such as \eqref{model-eps}.}   
\end{remark}

\subsection{Continuous weak formulation}
{For the rest of the paper}, we will use the following subspaces of $H^{1}(\Omega)$ for the function spaces to approximate the field variables, i.e., the displacements $\bfu$ and the temperature increment $\theta$:
\begin{align}
\widehat{V}^{\bfu} &= \left\{  \bfv \in \left(H^1(\Omega)\right)^{d} \colon \bfv(\bfx)=\bfu_0, \quad \forall \, \bfx \in \Gamma_D^{\bfu}\right\}, \label{eq:space1} \\
\widehat{V}^{\theta} &= \left\{  q \in H^1(\Omega) \colon q(\bfx)=\theta_0, \quad \forall \, \bfx \in \Gamma_D^{\theta}\right\}, \label{eq:space2}
\end{align}
and the test function spaces are defined as in \eqref{test-fun-space}. To pose a weak formulation, we multiply the equations in the strong formulation \eqref{strong-form-1} and  \eqref{strong-form-2} with the test functions from $V^{\bfu}$ and $V^{\theta}$ and then integrating by parts using Green's formula together with the boundary conditions given in \eqref{bcs-v1}, we arrive at the following weak formulation. \\

\noindent\textbf{Continuous weak formulation: }Given $\bfu^{0} \in \widehat{V}^{\bfu}$, for $n=0, 1, 2, \cdots$, find $\bfu^{n+1}:=\bfu^n + \delta \bfu^n \in \widehat{V}^{\bfu}$ and $\theta \in \widehat{V}^{\theta}$, such that 
\begin{subequations}\label{eq:weak_formulation}
\begin{align}
    A_{\theta}(\theta,q) &= L_{\theta}(q)\quad \forall\, q \in {V}_{\theta}, \\
    A_{\bfu}(\bfu^n; \, \delta\bfu^n,\bfv) &= L_{\bfu}(\bfu^n, \, \bfv)\quad \forall\, \bfv \in {V}_{\bfu},
\end{align}
\end{subequations}
where the bilinear terms $A_{\theta}(\theta, \, q), \,A_{\bfu}(\bfu^n; \delta\bfu^n, \, \bfv)$ and the linear term $L_{\theta}(q), \,L_{\bfu}(\bfv)$ are given by
\begin{align}
    A_{\theta}(\theta,q) &= \int_{\Omega} k \; \nabla\theta \,\cdot \nabla q  \;   d\Omega\, , \\
 L_{\theta}(q) &= \int_{\Omega} g \, q \; d\Omega\, ,
\end{align}

\begin{multline}\label{A-L-Def}
A_{\bfu}(\bfu^n; \, \delta\bfu^n,\, \bfv) = \int_{\Omega}  \Bigg[ \Bigg[
\frac{ \mu\left({ \nabla \delta\bfu^n + {\left(\nabla  \delta\bfu^n \right)^T}}\right)  + \lambda\, (\nabla \cdot  \delta \bfu^n)\,  \bfI }{\left( 1 - \left( \beta  | \mathbb{E}^{1/2} \left[ \bfeps \right]|\right)^{a}\right)^{1/a}} \\
+ \frac{\beta^{a}    \Theta_{1}\{  \bfu^n\}  \Theta_{2}\{  \bfu^n,  \delta\bfu^n \}  \mathbb{E}[ \bfeps]}{\left( 1 - \beta^{a} | \mathbb{E}^{1/2} \left[  \bfeps \right]|^{a}\right)^{1 + 1/a}} \Bigg] \colon \bfeps(\bfv) \bigg] d\Omega \, ,
\end{multline}
\begin{align}\label{A-L-Def-2}
    L_{\bfu}(\bfu^n, \,  \bfv)&= - \int_{\Omega} \alpha \;  \nabla\theta \; \bfv \; d\Omega  - \int_{\Omega} \left[\Bigg[\frac{ \mu \left( { \nabla  \bfu^n + \left(\nabla  \bfu^n \right)^T }  \right)  + \lambda \, (\nabla \cdot  \bfu^n)\,  \bfI  }{\left( 1 - \left( \beta  | \mathbb{E}^{1/2} \left[ \bfeps \right]|\right)^{a}\right)^{1/a}}\Bigg] \colon  \bfeps(\bfv)   \right]  \, d\Omega\, .
\end{align}
On the coarse grid, the initial guess $\bfu^1 \in \widehat{V}^{\bfu}$ for the \textit{Newton's method} is obtained by solving the linear problem with choosing $\beta =0$ in \eqref{A-L-Def}. Then, this value is projected onto the finite element mesh. 
{When $(\cdot)^n$ is denoted as the iteration number for the \textit{Newton's method},} a residual based stopping criterion can be  used to terminate the method.  

\subsection{Existence and uniqueness theorem}
We make the following assumptions on the data:
\begin{itemize}
\item[A1: ] For homogeneous material model, the heat conduction coefficient, $k \colon \Omega \to \mathbb{R}$ is assumed to be a global constant, but bounded with $0 < k_1 < k < k_2 < \infty$. For the heterogeneous material, $k$ is assumed to be $k(\bfx) \colon \Omega \to L^{\infty}(\Sym(\mathbb{R}^{d \times d}))$ and 
\[
0 < k_1 := \essinf_{\bfx \in \Omega}  \inf_{\bfv \in \mathbb{R}^d, \, \bfv \neq \bfzero} \dfrac{(k(\bfx)\bfv, \, \bfv )}{(\bfv, \, \bfv)}, \quad \infty > k_2 := \esssup_{\bfx \in \Omega}  \sup_{\bfv \in \mathbb{R}^d, \, \bfv \neq \bfzero} \dfrac{(k(\bfx)\bfv, \, \bfv )}{(\bfv, \, \bfv)}.
\]
\item[A2: ] For homogeneous solid model, the material constants, $\mu, \, \lambda, \, \alpha$ are all constants.  For heterogeneous material the parameters are:  $\mu, \, \lambda, \, \alpha \in L^{\infty}(\Omega, \, \mathbb{R)}$ and 
\[
0 < \mu_1 := \essinf_{\bfx \in \Omega} \mu(\bfx) \leq \esssup_{\bfx \in \Omega}  \mu(\bfx) =: \mu_2 < \infty.
 \]
Similarly, the other material parameters $\lambda, \, \alpha$ have the respective upper and lower bounds. 
 \item[A3: ]$g \in L^{2}(\Omega)$, $\bfu_0 \in \left( L^2(\Omega) \right)^d$, $\theta_0 \in L^2(\Omega)$.
\end{itemize}
In \cite{bulivcek2014elastic,bulivcek2015existence,beck2017existence,bulivcek2015analysis,bonito2020finite}, a detailed mathematical proof of the existence and uniqueness of weak solution to the strain-limiting constitutive description of the solid have been presented.  In this paper, we reframe our problem, especially the equation \eqref{eq2} so that we can compare our problem with the previous studies, in particular with \cite{beck2017existence}, concerning the existence and uniqueness of solution to the strain-limiting description: 
\begin{Theorem} [Main theorem in \cite{beck2017existence}]
 Let $\Omega \subset \mathbb{R}^d$ be bounded, connected, Lipschitz domain with open sets Dirichlet boundary $\Gamma_D$ and Neumann boundary $\Gamma_N$ such that $\Gamma_D \, \cap \, \Gamma_N = \emptyset$ and $\overline{\Gamma_D \cup \Gamma_N} = \partial \Omega$, given a vector field $\bff \colon \Omega \to \mathbb{R}^d$, a given traction $\bfg: \Gamma_N \to \mathbb{R}^d$, a given boundary load $\bfu_0 \colon \Gamma_D \to \mathbb{R}^d$, and a given bounded mapping $\mathcal{F} \colon \Sym(\mathbb{R}^{d \times d}) \to \Sym(\mathbb{R}^{d \times d})$, then find the pair $(\bfu, \, \bfT)$ such that $\bfu \colon \overline{\Omega} \to \mathbb{R}^{d}$, and $\bfT \colon \overline{\Omega} \to \Sym(\mathbb{R}^{d \times d})$, and 
\begin{align}
- \div \; \bfT &= \bff \quad \mbox{in} \quad \Omega, \notag \\
 \bfeps(\bfu) &= \mathcal{F}(\bfT) := \dfrac{\bfT}{(1 + \left( \beta \, \| \bfT \| \right)^a)^{1/a}}\quad \; a>0, \,\beta \geq 0 \quad \mbox{in} \quad \Omega, \notag \\
\bfu &= \bfu_0, \\
{\bfT  \bfeta} &= \bfg, \notag 
\end{align}
with following assumptions on the data: 

\noindent [B1] $\bff \in L^{2}(\Omega)^d$, $\bfg \in L^{2}(\Gamma_N)^d$, and 
\[
\bfzero = \int_{\Omega} \bff \, d\bfx + \int_{\partial \Omega} \bfg \, dS \quad in \quad \Gamma_N = \partial \Omega.
\]
[B2] $\bfu_0 \in W^{1, \, \infty}(\Omega)^d$, with $\nabla \bfu_0(\bfx)$ for almost every $\bfx \in \overline{\Omega}$ contained in a compact set in $\mathbb{R}^{d \times d}$. \\

\noindent Assume that the data $(\bff, \; \bfg, \; \bfu_0)$ satisfy the above assumptions [B1]-[B2], and consider $a >0$, 
then there exists a pair $(\bfu, \, \bfT) \in W^{1, \; \infty}(\Omega)^{d} \times \Sym(\L^1(\Omega))^{d \times d}$ satisfying 
\begin{equation}
\int_{\Omega} \bfT \cdot \bfeps(\bfw) \, d\bfx = \int_{\Omega} \bff \cdot \bfw \, d\bfx \quad \forall \; \bfw \in C^{1}_{0}(\Omega)^d.
\end{equation}
\end{Theorem}
\noindent The above problem studied in \cite{beck2017existence} is equivalent with the problem that we study in this paper, with $\bff = -\alpha \nabla \theta, \; \bfg = \bfzero$, and a following major change in the definition of the Cauchy stress tensor:
\begin{equation}
\bfT =\mathcal{L}(\bfeps(\bfu)) = \dfrac{ \mathbb{E}[\bfeps]         }{(1 - \left( \beta \, | \mathbb{E}^{1/2}[\bfeps] | \right)^a)^{1/a}}, \quad a>0, \,\beta \geq 0,
\end{equation}
where $\mathcal{L}(\cdot)$ is a uniform monotone operator with at most linear growth at infinity, and further, assuming $k, \,  \mu, \, \lambda, \, \theta_0, \, \bfu_0$ all satisfying the assumptions (A1)-(A3).  Hence, there exists a unique pair of $(\bfu, \; \bfT) \in (H_0^1(\Omega))^d \times \Sym(L^1(\Omega)^{d \times d})$.

\subsection{Finite element discretization of the model }
{In this section,} we first recall some basic notions and structure of the classical finite element method to discretize  the weak formulation \eqref{eq:weak_formulation}. The meshes 
used for computation in all the numerical examples presented in this paper are quadrilaterals.  Let $\left\{ \mathcal{T}_h \right\}_{h >0}$ be a conforming, shape-regular (in the sense of Ciarlet \cite{ciarlet2002finite}) family of triangulation of the domain $\Omega$; $\mathcal{T}_h$  is a finite family of sets such that:
\begin{itemize}
\item[(1)]  $K \in \mathcal{T}_h $ 
implies $K$ is an open simplex with the mesh size $h_{K} := \text{diam}(K)$ for each $K$. Furthermore, we denote the largest diameter of the triangulation by 
\[
h:=\max_{K \in \mathcal{T}_h } \;  h_K,
\]
\item[(2)] for any  $K_1, K_2 \in \mathcal{T}_h$, we have that $\overline{K}_1 \cap \overline{K}_2$ is either a null set or a vertex or an edge or the whole of $\overline{K}_1$ and $\overline{K}_2$,  
\item[(3)] $\bigcup\limits_{K \in \mathcal{T}_h} \overline{K} = \overline{\Omega}$,
\item[(4)] we now define the classical piece-wise affine \textit{finite element spaces}
\begin{align}\label{eq:FEM_Q}
S_h^{\bfu} &= \left\{  \bfu_h \in  \left( C(\overline{\Omega})\right)^d \colon \left. \bfu_h\right|_K \in \mathbb{Q}_k^d, \; \forall K \in \mathcal{T}_h \right\}, \\
S_h^{\theta}&= \left\{  q_h \in  C(\overline{\Omega}) \colon \left. q_h\right|_K \in \mathbb{Q}_k, \; \forall K \in \mathcal{T}_h \right\}, 
\end{align} 
where $\mathbb{Q}_k$ is a set containing the tensor-product of polynomials in $d$ variables up to order $k$ on the reference cell $\widehat{K}$. Then, the discrete approximation spaces are:
\begin{equation}\label{app-spaces}
\widehat{V}_h^{\bfu} = S_h^{\bfu} \, \cap  \, V^{\bfu}, \quad \widehat{V}_h^{\theta} = S_h^{\theta} \, \cap  \, V^{\theta}.
\end{equation}
\end{itemize}

\noindent\textbf{Discrete weak formulation: }{Given} $\bfu^{0}_h \in \widehat{V}_h^{\bfu}$ and the $n$-th Newton's iteration solution, i.e., $(\bfu^n_h, \, \theta_h) \in \widehat{V}_h^{\bfu} \times \widehat{V}_h^{\theta}$, for $n=0, 1, 2, \cdots,$ we 
find $\bfu^{n+1}_h =\bfu^n_h + \delta \bfu^n_h \in \widehat{V}_h^{\bfu}$ and $\theta_h \in \widehat{V}_h^{\theta}$ such that 
\begin{subequations}\label{discrete-wf}
\begin{align}
   A_{\theta}(\theta_h,q_h) &= L_{\theta}(q_h), \quad\forall\, q_h \in \widehat{V}_h^{\theta}, \label{eq:pr4}  \\
   A_{\bfu}( \bfu^n_h; \;        \delta\bfu^n_h,\bfv_h) &= L_{\bfu}( \theta_h, \, \bfu^n_h;  \;       \bfv_h), \quad\forall\, \bfv_h \in \widehat{V}_h^{\bfu}.\label{eq:pr5}  
\end{align}
\end{subequations}
On the finite element grid, the initial value for the Newton's method is obtained by solving the linear problem (i.e., with $\beta =0$ in \eqref{eq:pr5}) and then the solution is set to be a $L^2$-projection for the subsequent iterations. {The bilinear and linear terms} 
in the above discrete weak formulation \eqref{discrete-wf}  can be obtained by their corresponding definitions in the continuous variational formulation as defined in \eqref{eq:weak_formulation}. Finally, the matrices and the vectors in the both sides of \eqref{discrete-wf} are built by using the finite element basis functions defined in \eqref{app-spaces}. Some important remarks concerning the \textit{coercivity}, linear solvers, and stopping criterion for the \textit{Newton's method} are in order. 

\begin{remark}
{For the coercivity in the discrete weak form for \eqref{eq:pr5} (See \eqref{A-L-Def} and \eqref{A-L-Def-2}) and we note that 
\[
\left( 1-\beta^a |\mathbb{E}^{1/2}[\bfeps(\bfu_{h}^n)]|^{a})^{1/a} \right)    > 0 
\]
must be satisfied, thus a careful choice of the nonlinear parameters, $\beta$ and $a$, is required along with $|\mathbb{E}^{1/2}[\bfeps(\bfu_{h}^n)]|$ related to the Lam$\acute{e}$ parameters \cite{lee2020nonlinear}}
\end{remark}

\begin{remark}
In the above discrete weak formulation, the {the temperature ($\theta$) equation}~\eqref{eq:pr4} is solved first due its decoupling with the displacement equation \eqref{eq:pr5}. Moreover, the corresponding matrix assembled for the bilinear term in the equation~\eqref{eq:pr4} is a sparse, symmetric and positive definite. Hence, a Conjugate Gradient iterative solver \cite{saad2003iterative} can be utilized to find the temperature on the entire finite element mesh, {and ultimately, one needs to devise a special preconditioned CG solver. For the current application, we use ``symmetric successive over-relaxation (SSOR)'' as a preconditioner for the CG linear solver \cite{eisenstat1981efficient}.} {On the other hand,} it is worth noting that the matrix corresponding to the displacement {($\bfu$)} formulation in \eqref{eq:pr5} depends on the temperature values and a careful observation about the assembled matrix indicates that the matrix {may} become ill-conditioned when the boundary values are changed slightly and also the matrix clearly looses positive definite property when the gradients become higher. Hence, the CG solver may not be preferred for solving the linear system in \eqref{eq:pr5}. Therefore, a direct solver can be utilized to solve the displacement formulation.
\end{remark}

\begin{remark}\label{rmk5}
For a possible stopping criterion and also to monitor the convergence of the linear solver for \eqref{eq:pr5}, we need to compute the norm of the discrete residual, i.e., 
\begin{equation}\label{eqn:Residual}
\left( \mathcal{L}(\bfeps(\bfu_{h}^n)), \; \boldsymbol{\psi}_h \right) = \int_{\Omega} \frac{\mathbb{E}[\bfeps(\bfu_h^n)] \colon \nabla \boldsymbol{\psi}_h }{(1-\beta^a |\mathbb{E}^{1/2}[\bfeps(\bfu_h^n)]|^{a})^{1/a}} \; d\bfx.
\end{equation}
After each Newton iteration, we need to check this value before proceeding to next iteration or one can fix the number of iterations \textit{a priori}.
\end{remark}

The detailed numerical algorithm used in the computations for the nonlinear strain-limiting thermo-mechanical model is described in Algorithm~\ref{algo001}. \\

\begin{algorithm}[H]
\SetAlgoLined
\KwInput{Choose the parameters like $\beta, \, {a}, \, {\alpha_T}, \, \lambda, \, \mu, \, k$}
Assemble the {temperature equation} (Equation~\eqref{eq:pr4}) \;
  Solve for the temperature variable $\theta_h$ \; 
  Project $\theta_h$  on to the mesh \;
 \While{$[\text{Iteration Number} < \text{Max. Number of Iterations}]$.AND.$[\text{Residual} > \text{Tol.}]$}{
  Assemble the {displacement equation} (Equation~\eqref{eq:pr5}) \;
  Use a \textit{direct solver} to solve for $\delta \bfu^{n}$ \;
 Construct the solution at $n$th iteration using $\bfu_h^{n+1} = \bfu_h^{n} + \delta \bfu_h^{n}$ \;
  Calculate \textit{Residual} using  Equation~\eqref{eqn:Residual}\;
  \If{$\text{Residual}\leq\text{Tol.}$}{
   Break\;
   } }
 Write the solution pair $(\bfu_h^{n+1}, \; \theta_h)$ to output files for post-processing\;
 \caption{Algorithm for the nonlinear strain-limiting {thermo-mechanical} model}
 \label{algo001}
\end{algorithm}

\section{Numerical experiments}
{The primary purpose} of the present contribution is to test the validity of the proposed mathematical description against the classical thermoelasticity model.
{By achieving the goal of this research,
we can describe {the correct} material response {particularly around the crack-tip}, and {ultimately} to develop physically meaningful models for the crack evolution under thermo-mechanical loading.} 
{Within the framework of strain-limiting elasticity theory, {the newly proposed equations of thermoelasticity consists of} a coupled  system for linear-quasilinear partial differential equations \eqref{system1}.}
To date there is no closed form solution available to this system, therefore numerical algorithms can effectively be used for solving the problem formulated in this paper. 
{To this end, for the iteration parameters in the \textit{Newton's method}, we set the tolerance value (\textit{Tol})  as $10^{-8}$ for the convergence criterion and the maximum number of iterations (\textit{Max. Number of Iterations}) as $50$ (See Algorithm~\ref{algo001} and Remark~\ref{rmk5}).}

{In the following sections, {we use} the algebraic nonlinear strain-limiting theory to demonstrate our modeling and computational approach. To study the thermoelasticity problem,} we present  {numerical experiments for displacement, stress, strain, and temperature distributions with two examples of computational domain: a square domain in 2D without and with a slit, corresponding to Example 1 (Section~\ref{ex1}) and Example 2 (Section~\ref{ex2}), respectively}. In each example, a tensile top load is applied and two different temperature boundary conditions {(\texttt{CASE} 1 and \texttt{CASE} 2)} for the bottom part are considered. 
The numerical code is developed using an open source finite element library, {deal.II} \cite{dealII90} and all the computations are done utilizing the High Performance Computer cluster at Texas A\&M University, Corpus Christi.

\subsection{$h$-convergence study}
{First,} we perform the $h$-convergence study for the nonlinear mechanics solver described in the previous section {to verify the developed code with the algorithm}. 
We set a unit square as the domain, and a manufactured solution is chosen as
\begin{equation}\label{m-solution}
\bfu(x,y)=(\sin x\sin y, \cos x\cos y),
\end{equation}
and the boundary conditions are set based on the boundary trace of the above function \eqref{m-solution}. In this test, the temperature variable $\theta$ is taken to be zero on the entire computational mesh. Thus, we only solve the nonlinear mechanics, as in \eqref{eq:pr5}. With the pair of modeling parameter as $(a, \beta)=(1.0, 0.5)$ and the {Lam\'{e}} coefficients being set to be unity, {$6$ global mesh refinements are performed for the unit square. For each global refinement of the computational domain decreasing $h$ in half, we obtain $L^2$ error and the convergence rate}. Table~\ref{tab:ex1} depicts the optimal convergence rate of $2$, verifying the code and algorithm. 

\begin{table}[!h]
\centering
\small
\begin{tabular}{|c|l||l|l|}
\hline
{Cycle of refinement} & \multicolumn{1}{c||}{{$h$}} &   \multicolumn{1}{c|}{$L^2$ Error} & \multicolumn{1}{c||}{Rate}   \\    \hline\hline                        
1                      &  0.5            & 0.030660254881              & -      
 \\ \hline
2                      & 0.25             &  0.007495773956           & 2.7576
 \\ \hline
3                      & 0.125             & 0.001858903095             & 2.3722
 \\ \hline
4                      & 0.0625           & 0.000463682925               & 2.1833                                     \\ \hline
5                      & 0.03125            & 0.000115857772               & 2.0908
 \\ \hline
6                      &0.015625           & 0.000028961325             & 2.0452                                    \\ \hline
\end{tabular}
\caption{The results of $L^2$ error and rate of optimal convergence.}
\label{tab:ex1}
\end{table}%

\subsection{Computational domains and boundary conditions}
For our numerical simulations, we mainly consider the following two computational domains: $\Omega_1$ (Example 1) and $\Omega_2$ (Example 2) as portrayed in the following two figures (See Figure~\ref{Fig:ex1_setup} and Figure~\ref{Fig:ex2_setup}, respectively). {The only} difference in these domains is that Figure~\ref{Fig:ex2_setup} has a  ``slit'' at the right-hand side to the center (the blue bold line denoted by $\Gamma_C$), i.e., a line joining the points $(0.5,0.5)$ and $(1.0,0.5)$. The ``red'' dotted line is kept as a reference line along which we compute the stress-strain leading up to the crack tip, $(0.5,0.5)$. In all our simulations, the computational domain is discretized 7 times globally {(i.e., uniform meshes)} to get a mesh size $h = 0.0078125$. 
\begin{figure}[!h]
\centering
\begin{minipage}{0.45\textwidth}
\centering
\begin{tikzpicture}
\coordinate (A) at (0,0);
\coordinate (B) at (3,0);
\coordinate (C) at (3,3);
\coordinate (D) at (0,3);
  \draw (A) -- (B);
\draw (B) -- (C);  
\draw (C) -- (D);
 \draw (D) -- (A);
\node at (-0.3,-0.25)   {$(0,0)$};
\node at (3.3,3.2)   {$(1,1)$};
 \node [below] at (1.5,0.5) {$\Gamma_1$};
 \node [below] at (1.5,-0.1) {$\theta(\neq0)$};
 \node [right] at (3,1.5) {$\Gamma_2$};
 \node [left] at (0,1.5) {$\Gamma_4$};
 \node [above] at (1.5,2.4) {$\Gamma_3$};
 \node [above] at (1.5,1.2) {$\Omega_1$};
 \draw[->] (0.2,3.1) -- (0.2,3.5);
\draw[->] (1.,3.1) -- (1.,3.5);
\draw[->] (2.,3.1) -- (2.,3.5);
\draw[->] (2.8,3.1) -- (2.8,3.5);
\draw(0.25,-0.1) circle(0.1cm);
\draw(1.25,-0.1) circle(0.1cm);
\draw(1.75,-0.1) circle(0.1cm);
\draw(2.75,-0.1) circle(0.1cm);
 \node at (1.5,3.7) {$\bar{\bfu}=(0,1)$};
 \end{tikzpicture}
 \caption{Example 1. }
\label{Fig:ex1_setup}
\end{minipage}
\hspace{-2em}
\begin{minipage}{0.45\textwidth}
\centering
\begin{tikzpicture}
\draw (0,0) -- (3,0) -- (3,3) -- (0,3) -- (0,0);
\draw [line width=0.5mm, blue]  (1.5,1.5) -- (3,1.5);
\draw [loosely dotted, line width = 0.5mm, red] (0,1.5) -- (1.5,1.5);
\node at (-0.3,-0.25)   {$(0,0)$};
\node at (3.3,3.2)   {$(1,1)$};
\node at (2.2, 1.8)   {$\Gamma_{C}$};
\node [below] at (1.5,0.5) {$\Gamma_1$};
 \node [right] at (3,1.5) {$\Gamma_2$};
 \node [left] at (0,1.5) {$\Gamma_4$};
 \node [above] at (1.5,2.4) {$\Gamma_3$};
 \node [above] at (1.5,0.8) {$\Omega_2$};
 \draw[->] (0.2,3.1) -- (0.2,3.5);
\draw[->] (1.,3.1) -- (1.,3.5);
\draw[->] (2.,3.1) -- (2.,3.5);
\draw[->] (2.8,3.1) -- (2.8,3.5);
\node [below] at (1.5,-0.1) {$\theta(\neq0)$};
\draw(0.25,-0.1) circle(0.1cm);
\draw(1.25,-0.1) circle(0.1cm);
\draw(1.75,-0.1) circle(0.1cm);
\draw(2.75,-0.1) circle(0.1cm);
 \node at (1.5,3.7) {$\bar{\bfu}=(0,1)$};
 \end{tikzpicture}
\caption{Example 2. }
\label{Fig:ex2_setup}
\end{minipage}
\end{figure} 
 
\noindent The boundary conditions for the field variables are as follows:
\begin{subequations}\label{bcs-u}
 \begin{align}
 \text{on} \quad \;  \Gamma_{2}, \; \Gamma_4 \quad {\bfT  \bfn} &= \bfzero , \\
 \text{on} \quad \Gamma_{3} \quad \bar{\bfu} &= \left( 0, \, 1 \right), \\
 \text{on} \quad \Gamma_{1} \quad u_y &= 0, \label{eq:commom_boundary1}\\
 \text{on} \quad \Gamma_{2}, \;  \Gamma_{3}, \; \Gamma_{4} \quad  \nabla \theta \cdot \bfn  &= 0, \\
 \text{on} \quad \Gamma_{1} \quad \theta &=\theta(x),\label{eq:common_boundary2}
 \end{align}
 \end{subequations}
where values of the temperature variable ($\theta$) on $\Gamma_1$ are depending on the cases:
\begin{equation}\label{bcs-theta}
\begin{cases}
\texttt{CASE 1}:~\theta_1(x)=100,\\\texttt{CASE 2}:~\theta_2(x)=500 \, x \, (1-x).
\end{cases}
\end{equation}
\newline
On the slit/crack, $\Gamma_C $ in Example 2, the homogeneous Neumann boundary condition is imposed,
\begin{align}
    {\bfT   \bfn} &= \bfzero, \\
    \nabla\theta\cdot\bfn &= 0, \label{Ex2-temperature-slit}
\end{align}
 for the displacement and the temperature variables, respectively. The {Lam\'{e}} coefficients are set as $(\lambda, \, \mu)=(1.0,\, 1.0)$. Next, the modeling parameter pair $(a, \, \beta)$ for the nonlinear case is set to be $(0.5, \, 0.02)$ which is chosen purely for a quick computational purpose. Note that for the linear case, $\beta$ is set to be $0$. 
All parameters involved in the model along with their corresponding {SI units} are presented in the Table~\ref{Tab:Exs-12}. {A detailed study of the effect of these parameters on the behavior of the solution near the crack-tip is an important and interesting topic that will be investigated in future contributions.}
 Our emphasis in this paper is to show that the model studied in this paper is a natural extension of the {established} nonlinear strain-limiting approach: {from elasticity}
to thermoelasticity. {With the verification of the model,} another focus of this contribution is 
comparison with the classical model {(or the LEFM)}.
\newline
\begin{table}[H]
\small
\begin{center}
\begin{tabular}{ccl}
\toprule
\textbf{Parameter} & \textbf{Value} & {\textbf{Unit}}\\
\midrule
Lam\'{e} coefficient $\lambda$ & $1.0$ & Pa\\
Lam\'{e} coefficient $\mu$ & $1.0$ & Pa \\
Nonlinear parameter $a$ & $0.5$ &  - \\
Nonlinear parameter $\beta$ & $0.02$ & - \\
Heat source $g$ & $-10$ & $J/(m^2s)$ \\
Heat conduction coefficient $k$ & $20$ & $J/(m^2sK)$ \\
Linear thermal expansion coefficient $\alpha_{T}$ & $0.1$  &  $\ \ /K$ 
\\
\bottomrule
\end{tabular}
\caption{Common parameters for numerical simulation of the proposed model.}
\label{Tab:Exs-12}
\end{center}
\end{table}

\subsection{Numerical experiments and discussion}\label{sec:Numerical_results}
{In this section}, we present the numerical solutions to the boundary value problems proposed in the previous section: {Example 1 (Section~\ref{ex1}) and  Example 2 (Section~\ref{ex2}). }  
In each example, both linear and nonlinear problems are solved on the same mesh. {Note that we have the same order of $k=1$, i.e., $\mathbb{Q}_1$, for polynomial functions ($\mathbb{Q}_k$ in \eqref{eq:FEM_Q}) in the finite element space approximating the temperature and the displacement.}

The variation in the modeling parameter $\beta$ gives rise to plethora of nonlinear models. {Again} in our computations, we choose $\beta = 0.02$ to depict the effect of $\beta$ 
{to compare with the model when}  $\beta =0$, i.e., the classical linearized model,  {
which implies that strains are still large uniformly including near the concentrators such as crack-tips.}  A detailed numerical study is required for the optimal choice of the parameter $\beta$ that not only ``limits’’ the strains but also ensures the convergence of the numerical method.

\subsubsection{Example 1: Domain without a slit}\label{ex1}
{Example 1 is addressing a homogeneous unit square domain without any slit or crack.} As in Figure~\ref{Fig:ex1_setup}}, all the boundary parts of the domain (i.e., a unit square of $[0, 1] \times [0, 1]$) are marked for applying boundary conditions for both the field variables. {First,} we study the temperature distribution in a thermoelastic solid described by a special constitutive relationship \eqref{model-T}. 
\begin{figure}[H]
\centering
\subfloat[\texttt{CASE} 1]{\includegraphics[width=0.45\textwidth]{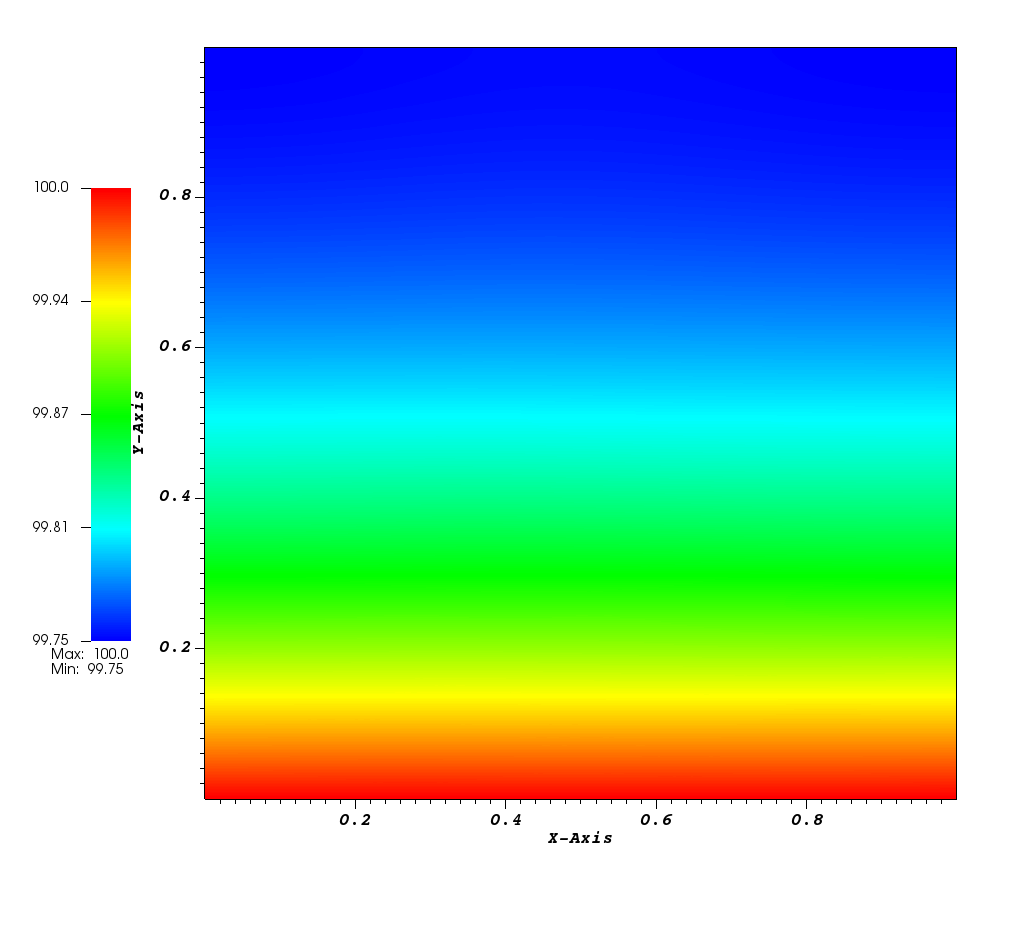}} 
\hspace{0.1in}
\subfloat[\texttt{CASE} 2]{\includegraphics[width=0.45\textwidth]{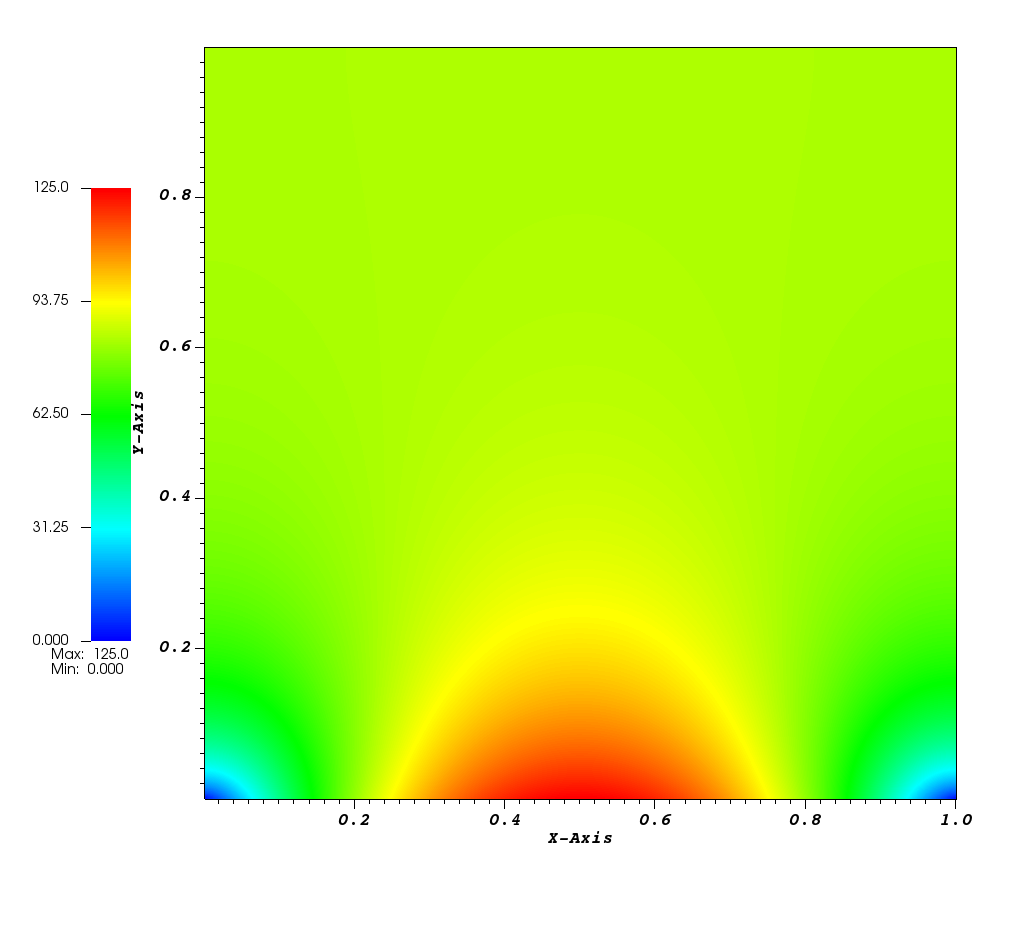}} 
\caption{(Example 1) Temperature distribution for \texttt{CASE} 1 and \texttt{CASE} 2.}
\label{figs:T_Ex1}
\end{figure}
\noindent \textbf{Temperature\quad} The temperature distribution in this example {with both cases (\texttt{CASE} 1 and 2 in \eqref{bcs-theta})} are shown in Figure \ref{figs:T_Ex1}. Since \texttt{CASE} $1$ (Figure~\ref{figs:T_Ex1} (a)) has a constant temperature value on the boundary $\Gamma_1$, i.e., $\theta_1(x)=100, \; 0 \leq x \leq 1$, {and a uniform heat source $g=-10$ working as a sink over the domain}, therefore we can see a uniform distribution of temperature {in each layer from the bottom}.   
{For the same value of heat source,} \texttt{CASE} $2$ (Figure~\ref{figs:T_Ex1} (b)) has a quadratic function $\theta_2(x)=500x\,(1-x), \; 0 \leq x \leq 1$, on the bottom boundary  $\Gamma_1${, resulting in the non-uniform temperature distribution at each height near the bottom}. 
{With these known solutions of temperature effecting on the momentum, next} we present the numerical solutions for the displacement variable in the material body.  
\\
\newline
(\texttt{CASE} $1$) \textbf{Linear\quad} By taking $\beta =0$ in \eqref{eq:pr5}, we obtain the solution to the linear model. Figure \ref{figs:U_Ex1_Case1_Linear} depicts the displacement components, (a) $\bfu_x$ and (b) $\bfu_y$ 
in the whole domain.   
Figure \ref{figs:S_E_Ex1_Case1_Linear} presents the axial stress ($\bfT_{yy}$) and the axial strain ($\bfeps_{yy}$) plots, respectively. 

\begin{figure}[H]
\centering
\subfloat[{$\bfu_x$}]{\includegraphics[width=0.45\textwidth,trim=4 4 4 4,clip]{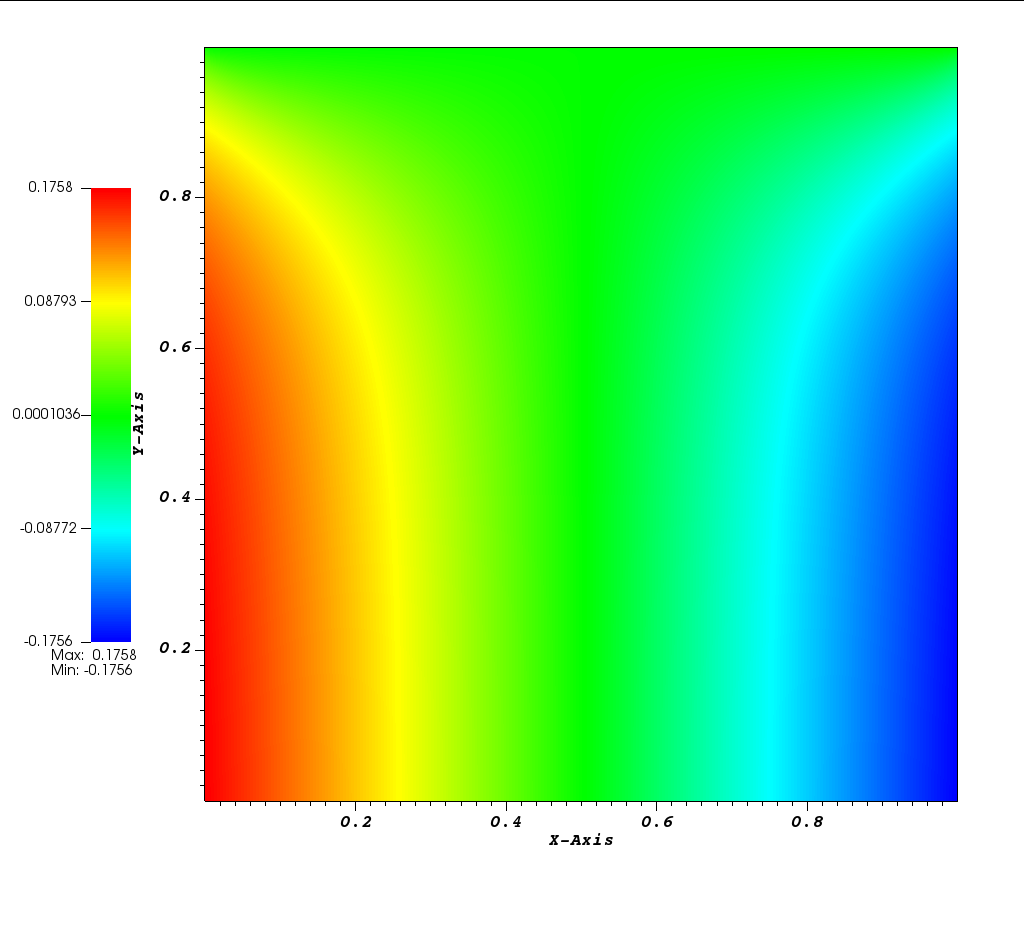}}
\hspace{0.1in}
\subfloat[{$\bfu_y$}]{\includegraphics[width=0.45\textwidth,trim=4 4 4 4,clip]{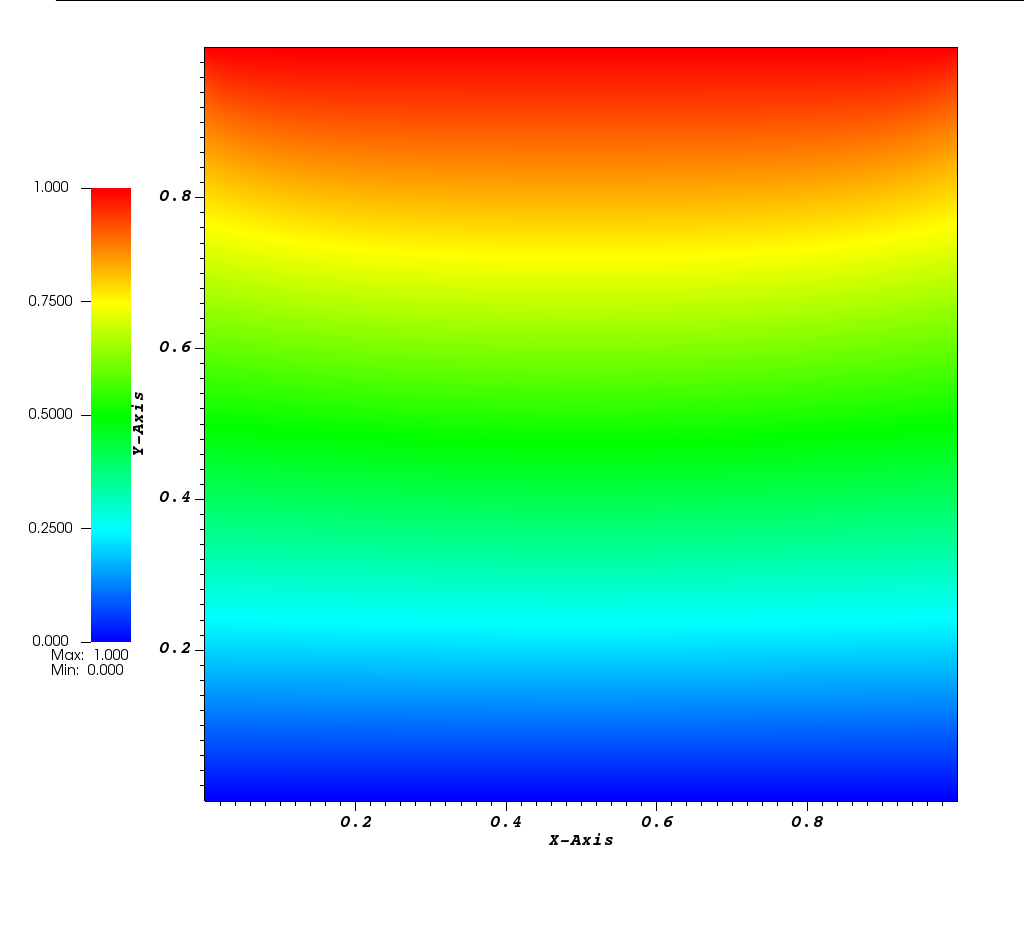}}
\caption{(Example 1 with \texttt{CASE} 1) Displacement for the Linear: (a) x-displacement and (b) y-displacement. }
\label{figs:U_Ex1_Case1_Linear}
\end{figure}

\begin{figure}[H]
\centering
\subfloat[{$\bfT_{yy}$}]{\includegraphics[width=0.45\textwidth,trim=4 4 4 4,clip]{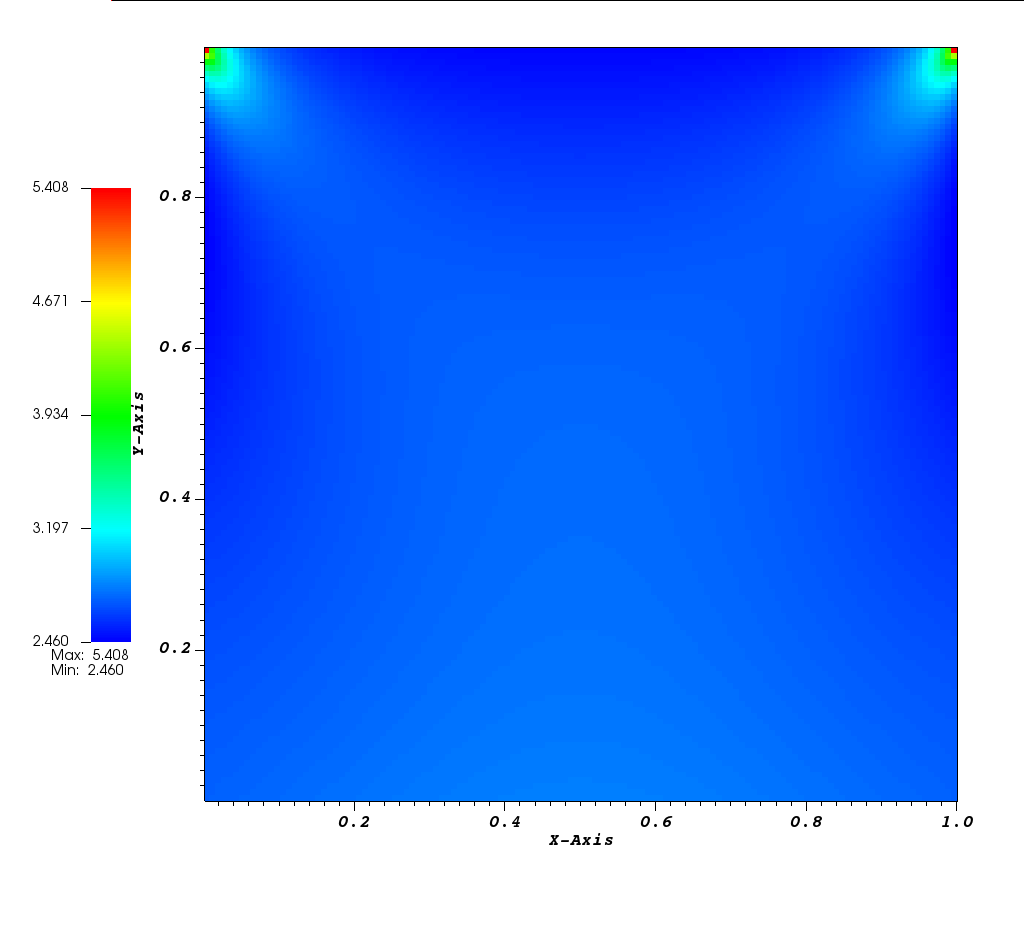}}
\hspace{0.1in}
\subfloat[{$\bfeps_{yy}$}]{\includegraphics[width=0.45\textwidth,trim=4 4 4 4,clip]{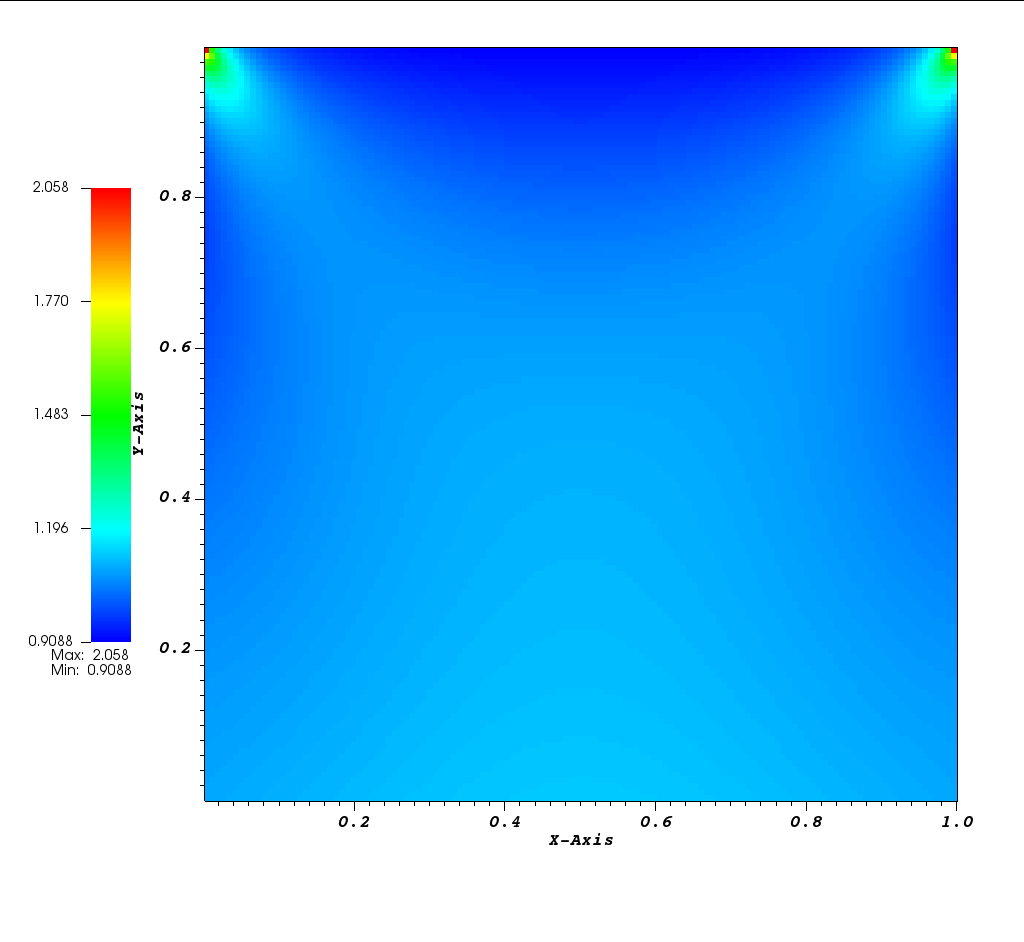}}
\caption{(Example 1 with \texttt{CASE} 1) Stress and Strain for the Linear: (a) axial stress and (b) axial strain.  }
\label{figs:S_E_Ex1_Case1_Linear}
\end{figure}
\begin{figure}[H]
\centering
\subfloat[{$\bfu_x$}]{\includegraphics[width=0.45\textwidth,trim=4 4 4 4,clip]{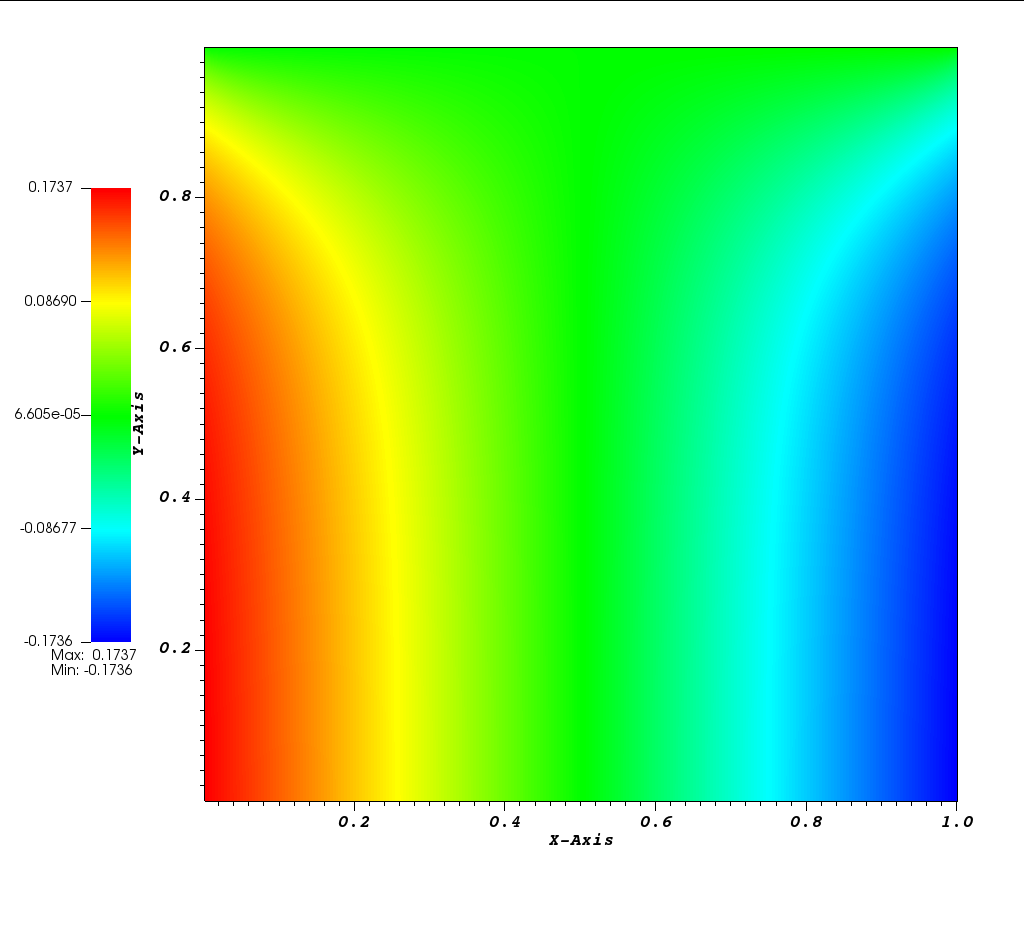}}
\hspace{0.1in}
\subfloat[{$\bfu_y$}]{\includegraphics[width=0.45\textwidth,trim=4 4 4 4,clip]{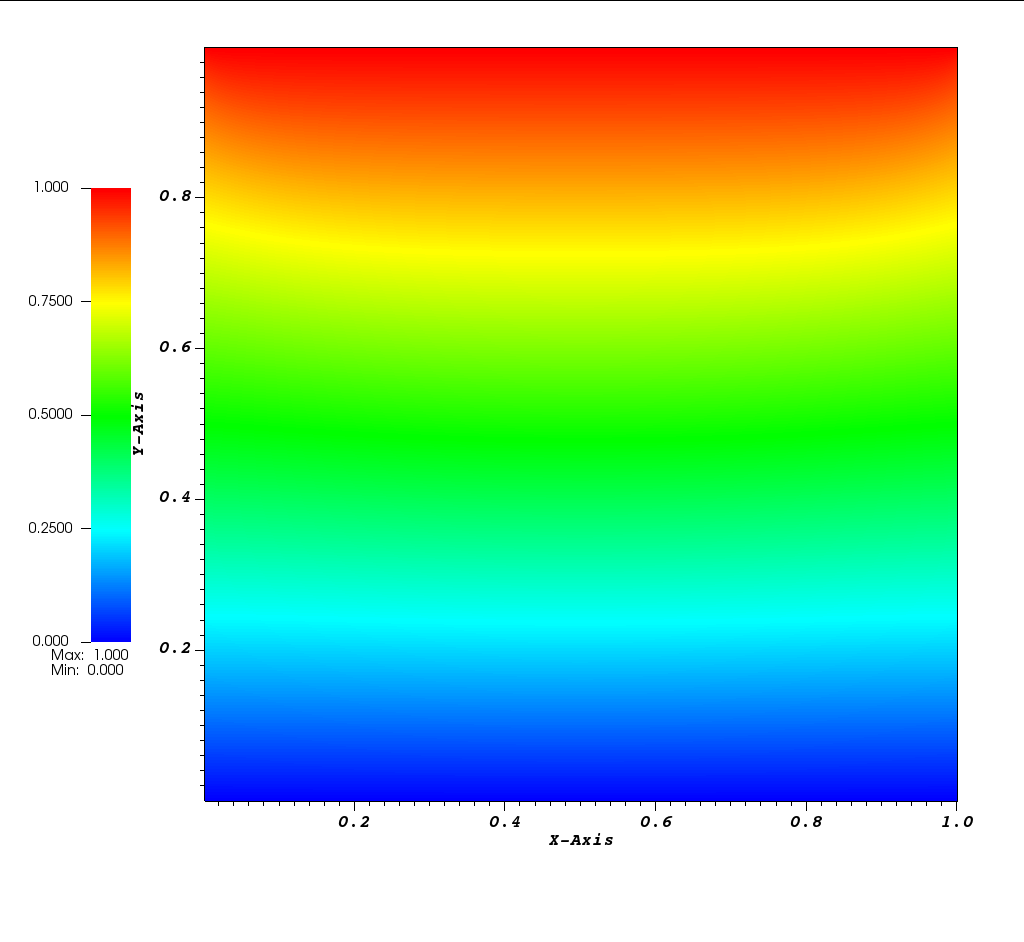}}
\caption{(Example 1 with \texttt{CASE} 1) Displacement for the Nonlinear: (a) x-displacement and (b) y-displacement.}
\label{figs:U_Ex1_Case1_NL_A0p5_B0p01}
\end{figure}
\begin{figure}[H]
\centering
\subfloat[{$\bfT_{yy}$}]{\includegraphics[width=0.45\textwidth,trim=4 4 4 4,clip]{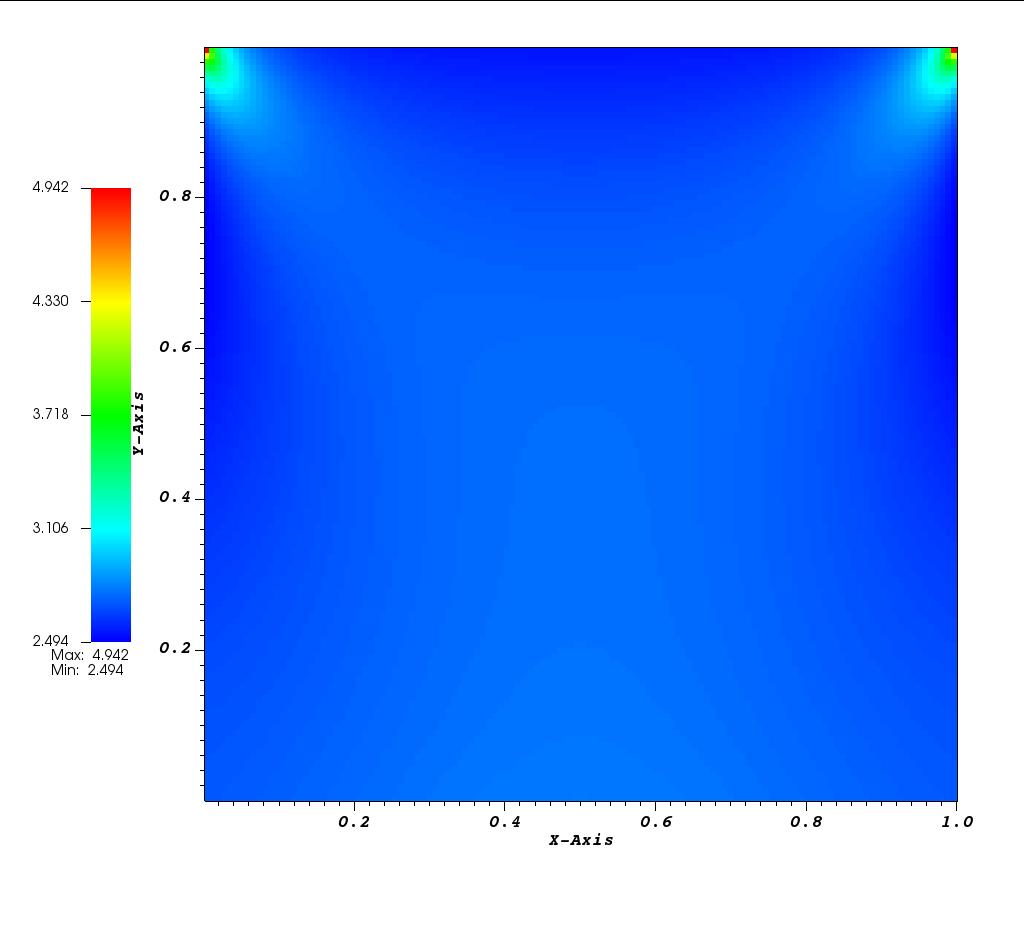}}
\hspace{0.1in}
\subfloat[{$\bfeps_{yy}$}]{\includegraphics[width=0.45\textwidth,trim=4 4 4 4,clip]{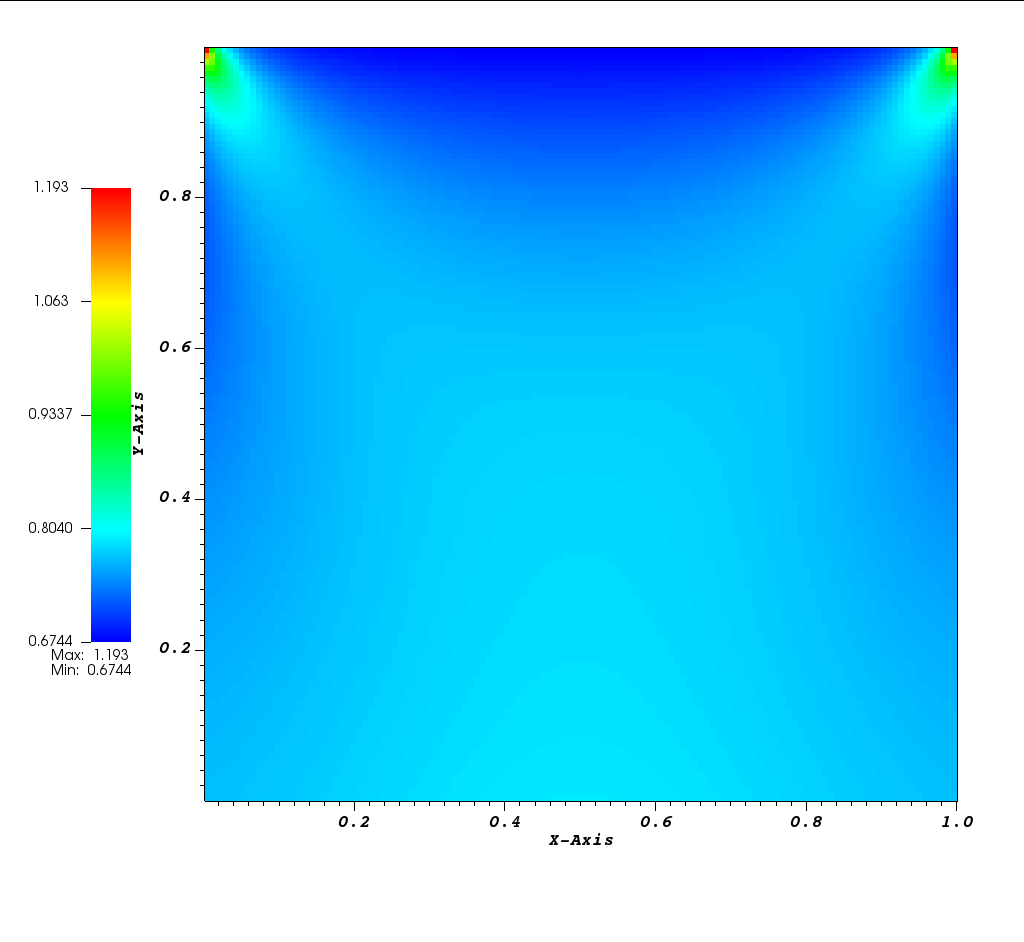}}
\caption{(Example 1 with \texttt{CASE} 1) Stress and Strain for the Nonlinear: (a) axial stress and (b) axial strain.} 
\label{figs:S_E_Ex1_Case1_NL_A0p5_B0p01}
\end{figure}

\noindent(\texttt{CASE} $1$) \textbf{Nonlinear\quad}
Here, we solve the nonlinear model with  \texttt{CASE} $1$ type temperature boundary condition. 
Figure \ref{figs:U_Ex1_Case1_NL_A0p5_B0p01} depicts the displacement components $\bfu_x$ and $\bfu_y$ and the axial stress ($\bfT_{yy}$) and the axial strain ({$\bfeps_{yy}$}) components are presented in Figure \ref{figs:S_E_Ex1_Case1_NL_A0p5_B0p01}. 
{Note that for displacement, the overall pattern of the distributions and values (e.g., the maximum and the minimum values in $\bfu_x$) is almost the same as the one for (\texttt{CASE} $1$) \textbf{Linear} (Figure \ref{figs:U_Ex1_Case1_Linear}). In case for the axial stress ($\bfT_{yy}$) and the axial strain ({$\bfeps_{yy}$}) for the nonlinear model, the patterns are similar to the linear model, but the the maximum and the minimum in the axial strain is much lower for the nonlinear model as expected. Compared to the axial strain, we see that the difference in axial stress between the two models is not much intense.}

\begin{figure}[H]
\centering
\subfloat[{$\bfu_x$}]{\includegraphics[width=0.45\textwidth,trim=4 4 4 4,clip]{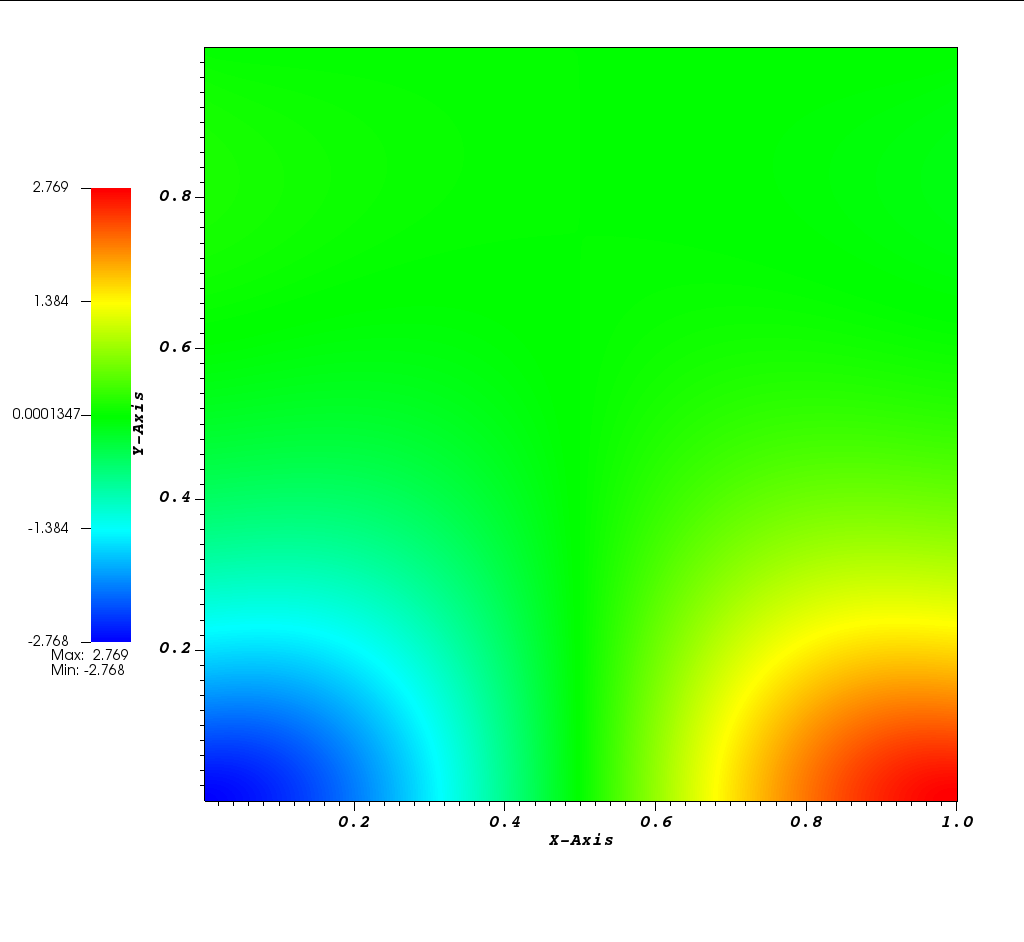}}
\hspace{0.1in}
\subfloat[{$\bfu_y$}]{\includegraphics[width=0.45\textwidth,trim=4 4 4 4,clip]{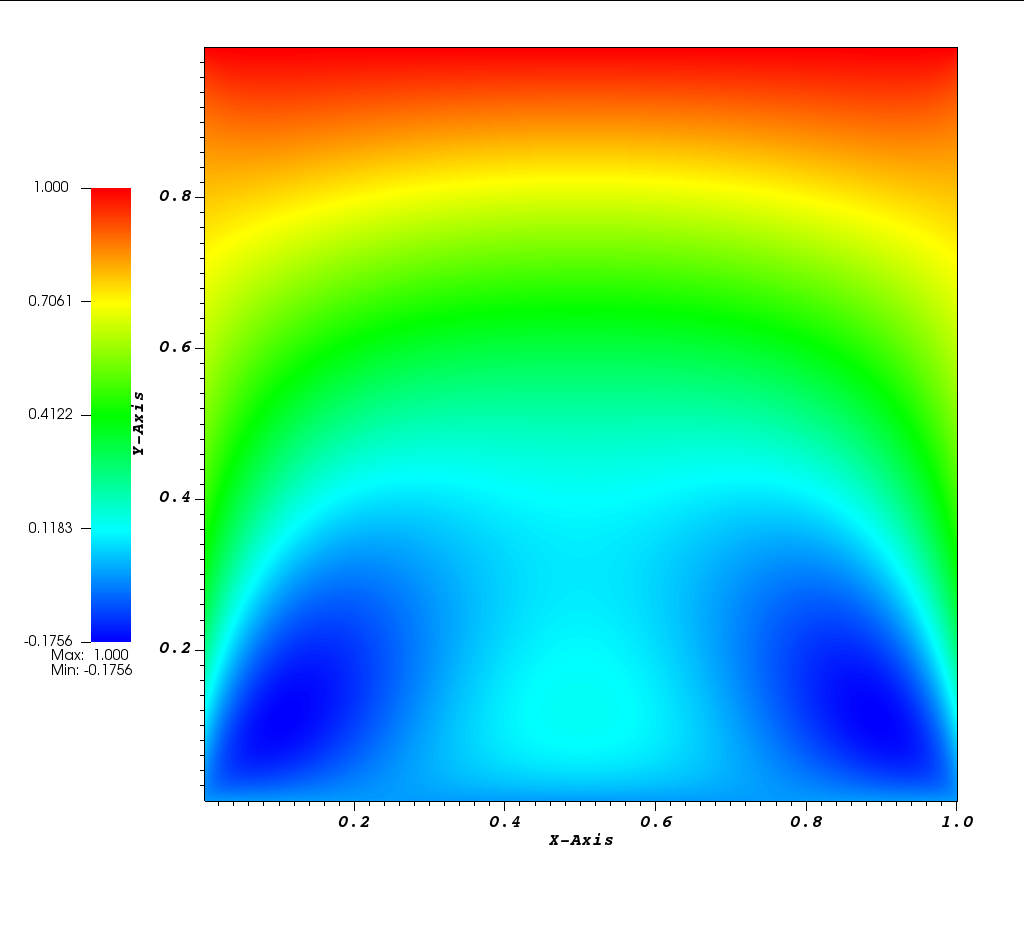}}
\caption{(Example 1 with \texttt{CASE} 2) Displacement for the Linear: (a) x-displacement and (b) y-displacement }
\label{figs:U_Ex1_Case2_Linear}
\end{figure}

\begin{figure}[H]
\centering
\subfloat[{$\bfT_{yy}$}]{\includegraphics[width=0.45\textwidth,trim=4 4 4 4,clip]{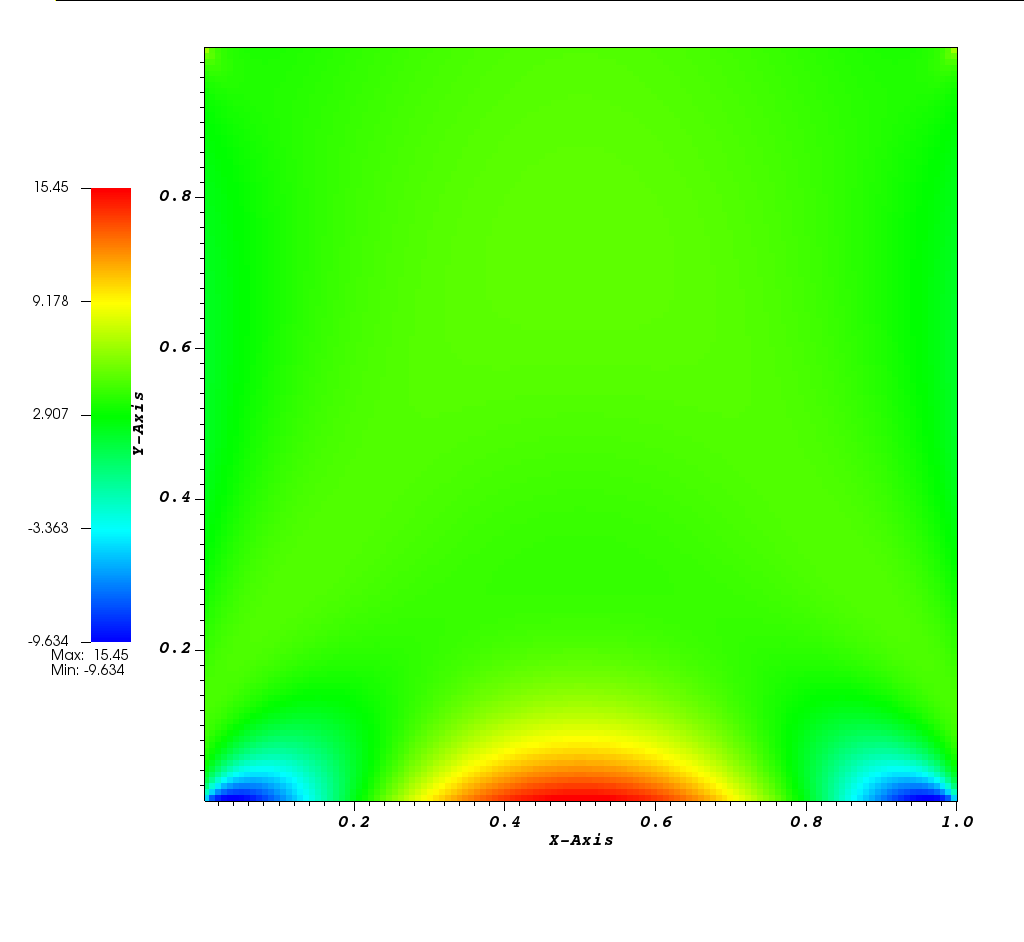}}
\hspace{0.1in}
\subfloat[{$\bfeps_{yy}$}]{\includegraphics[width=0.45\textwidth,trim=4 4 4 4,clip]{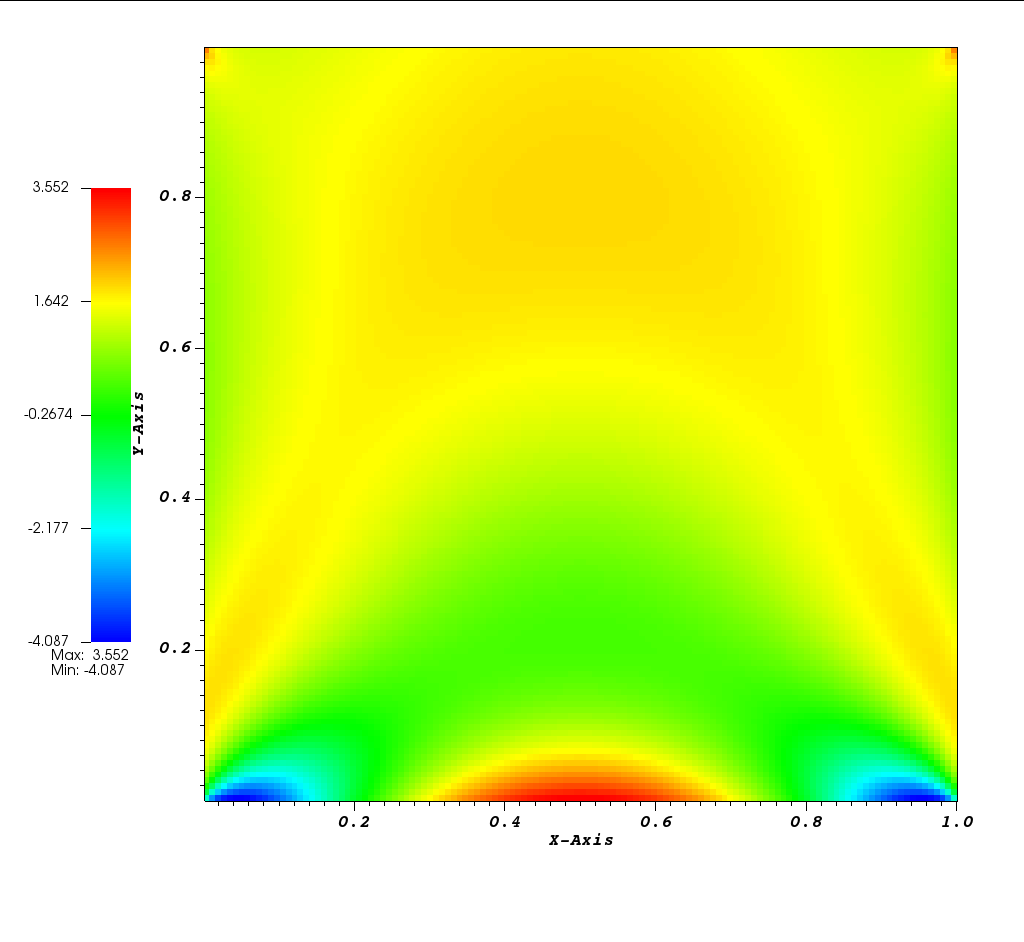}}
\caption{(Example 1 with \texttt{CASE} 2) Stress and Strain for the Linear: (a) axial stress and (b) axial strain.} 
\label{figs:S_E_Ex1_Case2_Linear}
\end{figure}

\noindent(\texttt{CASE} $2$) \textbf{Linear\quad}{Next,} we present the computational results for the second type of temperature boundary condition (\texttt{CASE} 2) 
with the linear model. Figure \ref{figs:U_Ex1_Case2_Linear} shows the displacement field in the entire body. We see that the displacement component $\bfu_y$ has a parabolic type distribution on the bottom face due to the applied temperature boundary condition on $\Gamma_1$. 
Corresponding distribution of axial stress-strain in $y$-direction is shown in Figure \ref{figs:S_E_Ex1_Case2_Linear}. 

\begin{figure}[H]
\centering
\subfloat[{$\bfu_x$}]{\includegraphics[width=0.45\textwidth,trim=4 4 4 4,clip]{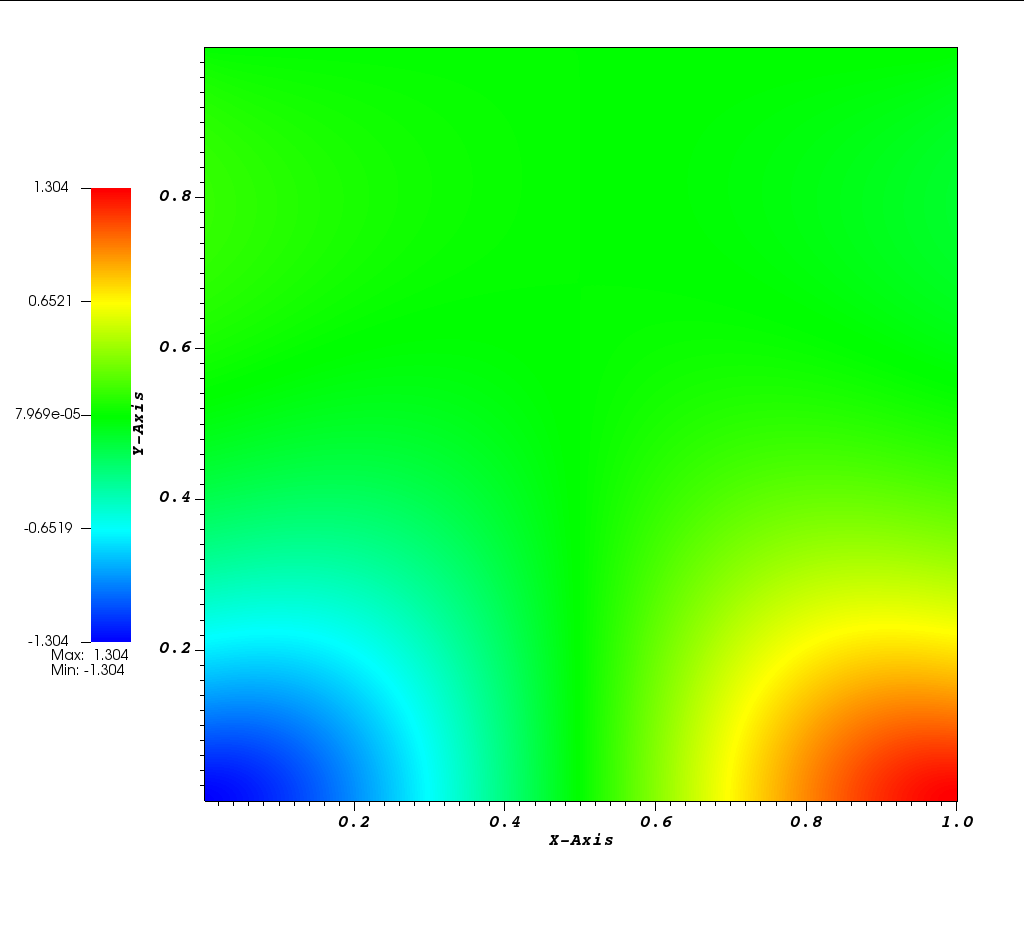}}
\hspace{0.1in}
\subfloat[{$\bfu_y$}]{\includegraphics[width=0.45\textwidth,trim=4 4 4 4,clip]{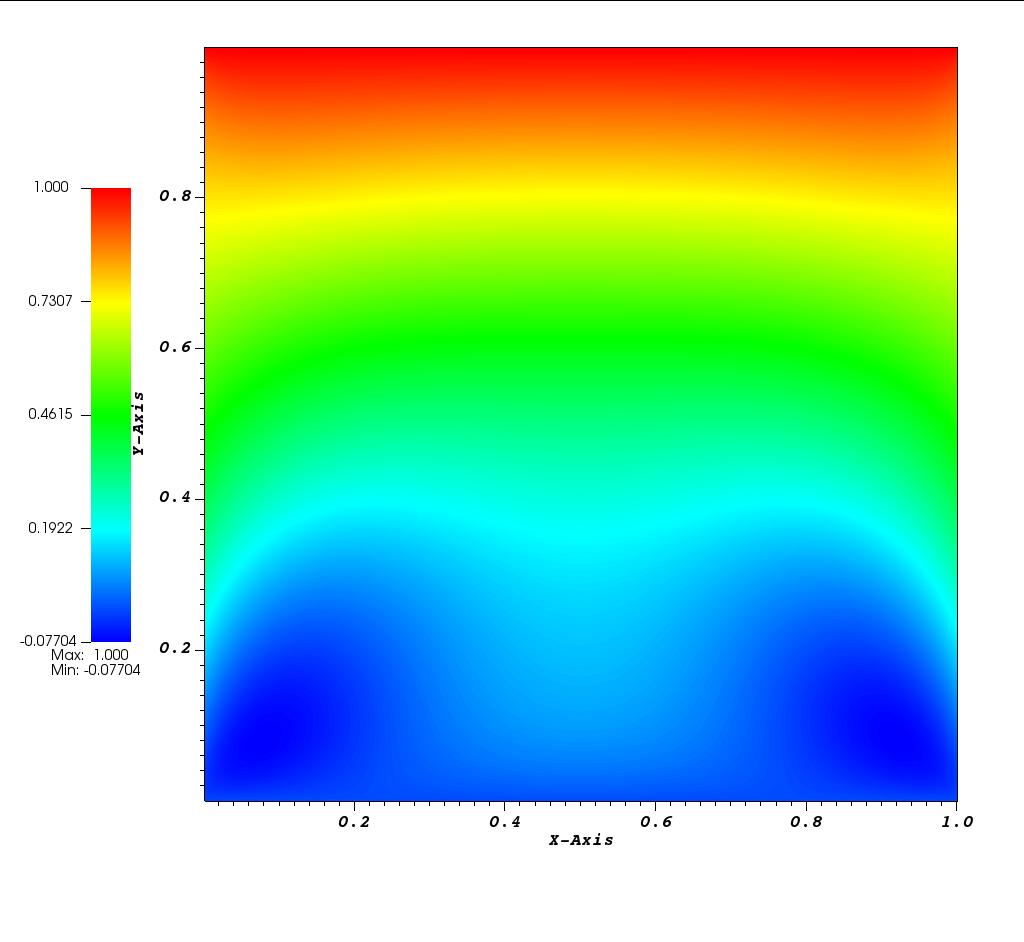}}
\caption{(Example 1 with \texttt{CASE} 2) Displacement for the Nonlinear: (a) x-displacement and (b) y-displacement.}
\label{figs:U_Ex1_Case2_NL_A0p5_B0p01}
\end{figure}

\begin{figure}[H]
\centering
\subfloat[{$\bfT_{yy}$}]{\includegraphics[width=0.45\textwidth,trim=4 4 4 4,clip]{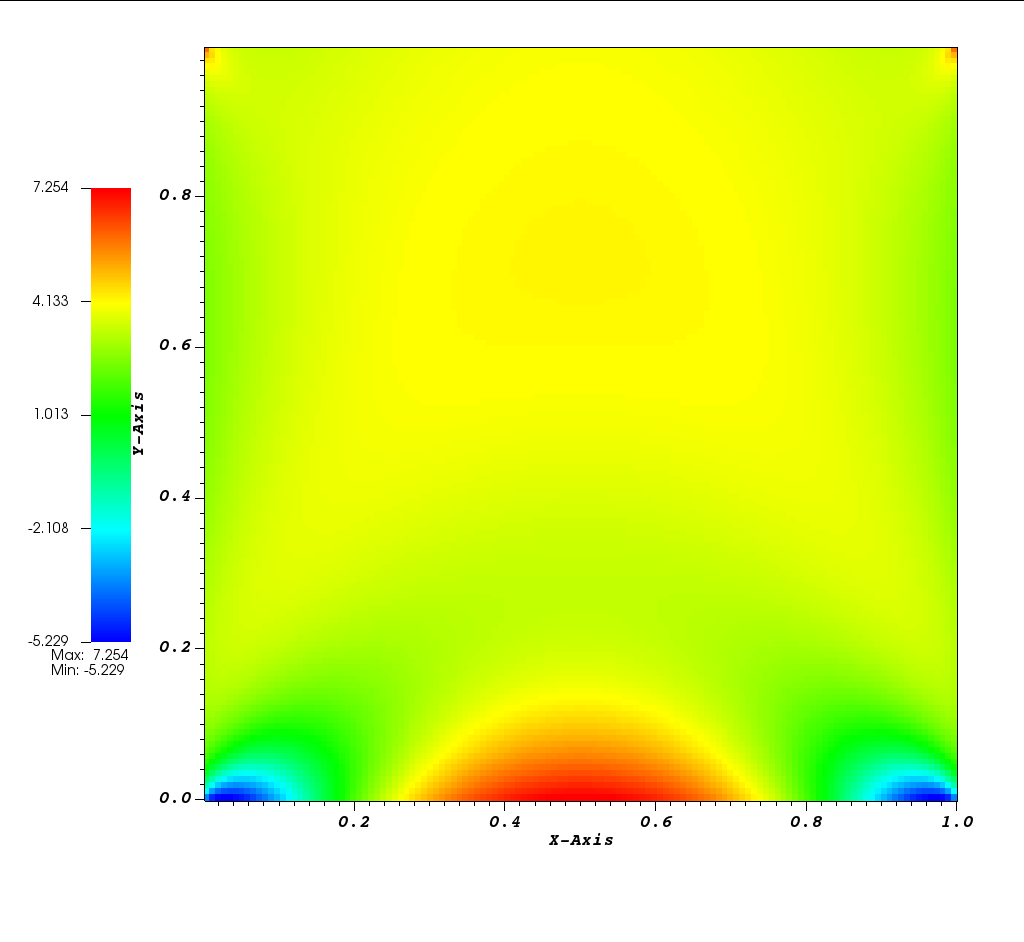}}
\hspace{0.1in}
\subfloat[{$\bfeps_{yy}$}]{\includegraphics[width=0.45\textwidth,trim=4 4 4 4,clip]{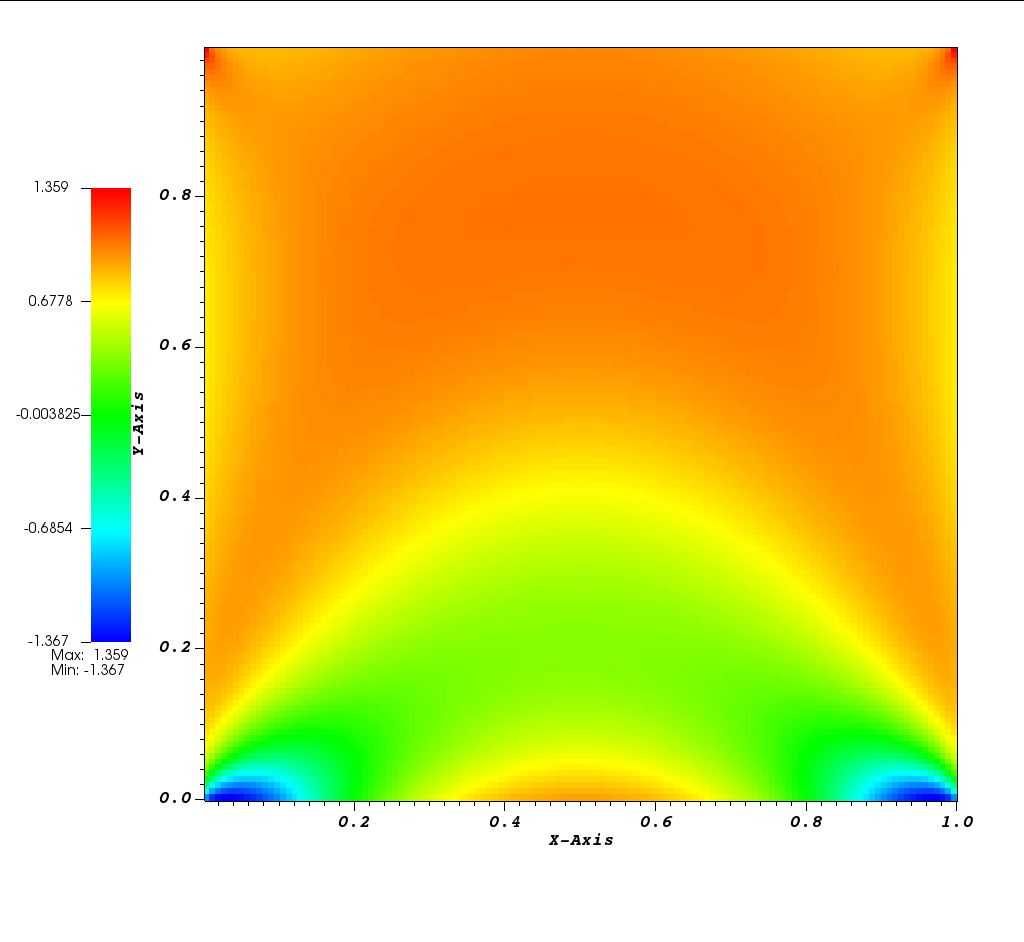}}
\caption{(Example 1 with \texttt{CASE} 2) Stress and Strain for the Nonlinear: (a) axial stress and (b) axial strain.}  
\label{figs:S_E_Ex1_Case2_NL_A0p5_B0p01}
\end{figure}

\noindent(\texttt{CASE} $2$) \textbf{Nonlinear\quad}We explore the effect of the quadratic temperature boundary condition on the response of the material body characterized by the nonlinear relation between stress and strain. 
Figure \ref{figs:U_Ex1_Case2_NL_A0p5_B0p01} represents the displacement field in this case. {Unlike \texttt{CASE} 1, the parabolic distribution of bottom temperature in \texttt{CASE} 2 within the nonlinear model has noticeble effects on the mechanics and brings different maximum and minimum values between the models for each displacement component (Figure~\ref{figs:U_Ex1_Case2_NL_A0p5_B0p01}). 
Similarly, distinct patterns are found in axial stress and strain distributions (Figure~\ref{figs:S_E_Ex1_Case2_NL_A0p5_B0p01}).} 
{For the linear model, the maximum and minimum values for both stress and strain aggregate densely closer to the bottom boundary $\Gamma_1$ where the non-uniform (or closer to a point/region) heat source is located. Meanwhile, the nonlinear model has more distributed pattern over the domain with smaller values, especially for the strain. Similar type of effects of distribution patterns for the linear and the nonlinear models can be found in Example 2 in the following where the tip area of slit works as another boundary related to the stress. Compared to \texttt{CASE} 1, overall we see that \texttt{CASE} 2 has more different distributions of the field variables between the linear and the nonlinear models.}


\subsubsection{Example $2$: Domain with a slit}\label{ex2} 
\begin{figure}[H]
\centering
\subfloat[\texttt{CASE} 1]{\includegraphics[width=0.45\textwidth,trim=4 4 4 4,clip]{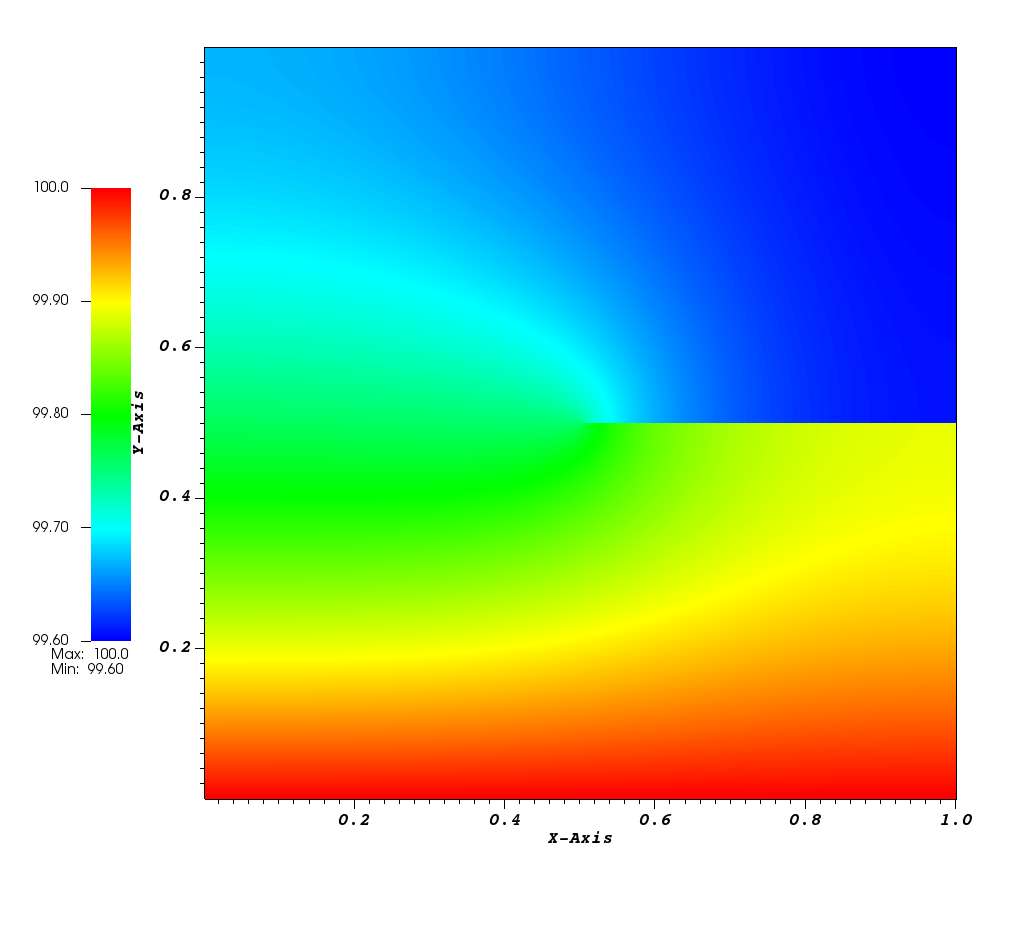}}
\hspace{0.1in}
\subfloat[\texttt{CASE} 2]{\includegraphics[width=0.45\textwidth,trim=4 4 4 4,clip]{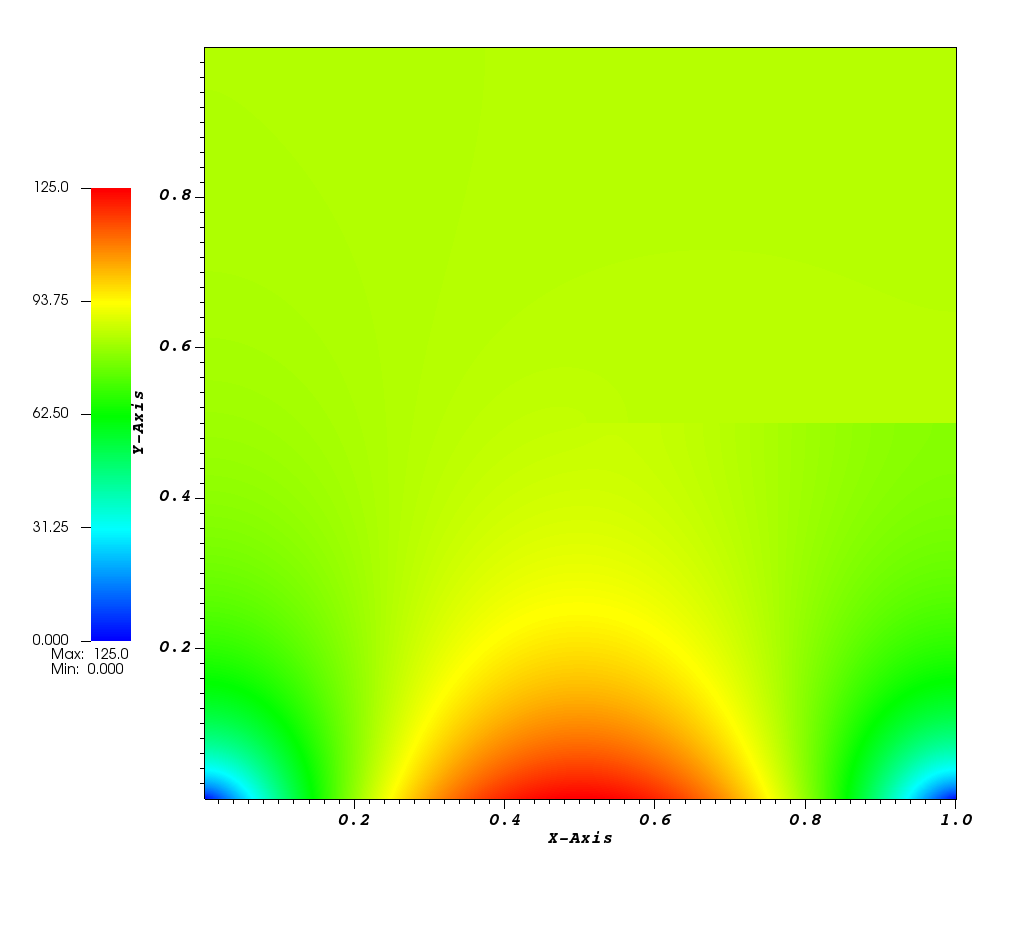}}
\caption{(Example 2) Temperature distribution for \texttt{CASE} 1 and \texttt{CASE} 2.}
\label{figs:T_Ex2}
\end{figure}

In this example, we have considered has a computational domain with a ``slit'' as shown in Figure \ref{Fig:ex2_setup}. {The slit is expressed with two individual vertices at the same coordinate of $(1, 0.5)$ on the right boundary $\Gamma_2$, resulting in a ``crack-opening'' deformation on  $\Gamma_C$.} 
A zero traction boundary condition is imposed on the upper and lower faces of the slit. All other boundary conditions for both variables on all other remaining parts of the boundary are the same as the ones used in Example 1 (as given in \eqref{bcs-u}-\eqref{bcs-theta}). We also consider the two different temperature boundary conditions given in \eqref{bcs-theta}. 
\\
\newline
\textbf{Temperature\quad} {Figure} \ref{figs:T_Ex2} shows two cases of temperature distribution in this example with two different temperature boundary conditions. {Be reminded that \texttt{CASE} 1 has constant temperature boundary condition on the boundary $\Gamma_1$. Compared to Example 1 with \texttt{CASE} 1 (Figure~\ref{figs:T_Ex1} (a))}, Example 2 
with \texttt{CASE} 1 has the temperature distribution {less} uniform 
due to the slit at the center of the half right region. {Note that we have the homogeneous Neumann boundary condition for temperature on the slit as \eqref{Ex2-temperature-slit}.} For \texttt{CASE} 2, {quite similar type of} the parabolic shape for heat distribution from the bottom face can be seen {as Example 1 with \texttt{CASE} 2 (Figure~\ref{figs:T_Ex1} (b))}.

\begin{figure}[H]
\centering
\subfloat[{$\bfu_x$}]{\includegraphics[width=0.45\textwidth,trim=4 4 4 4,clip]{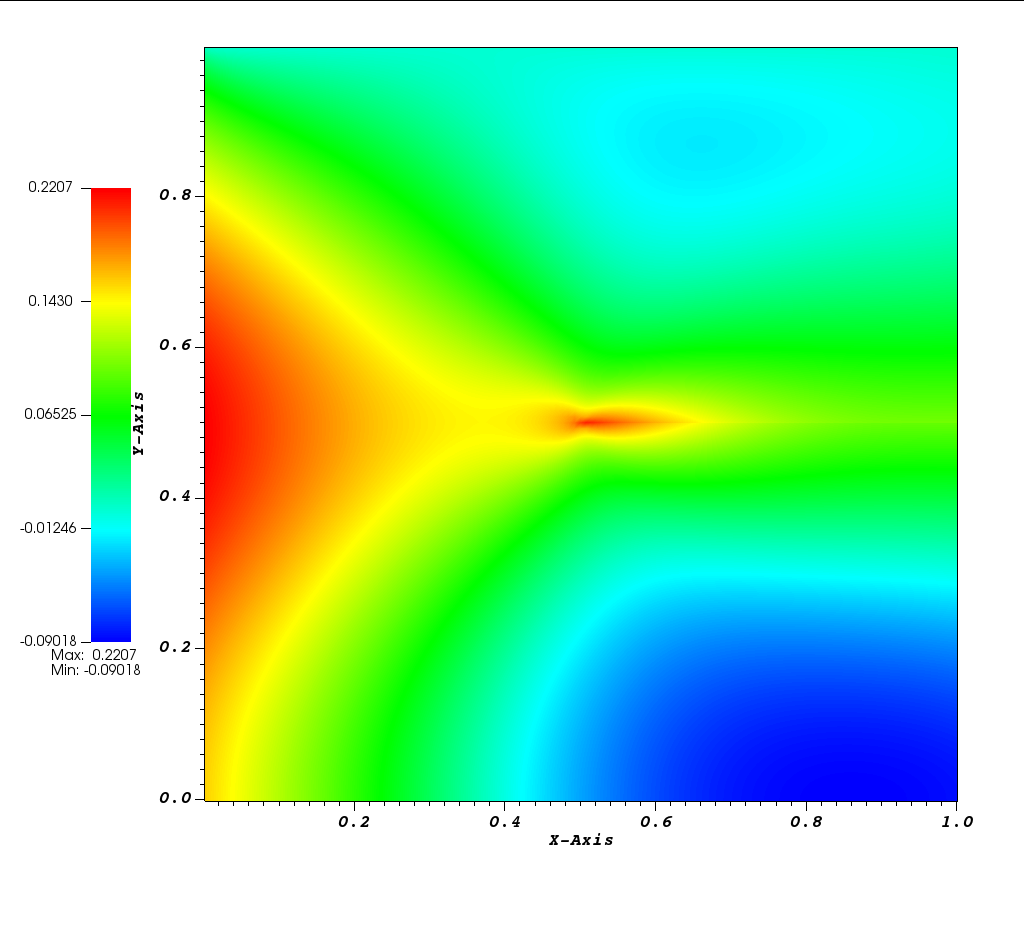}}
\hspace{0.1in}
\subfloat[{$\bfu_y$}]{\includegraphics[width=0.45\textwidth,trim=4 4 4 4,clip]{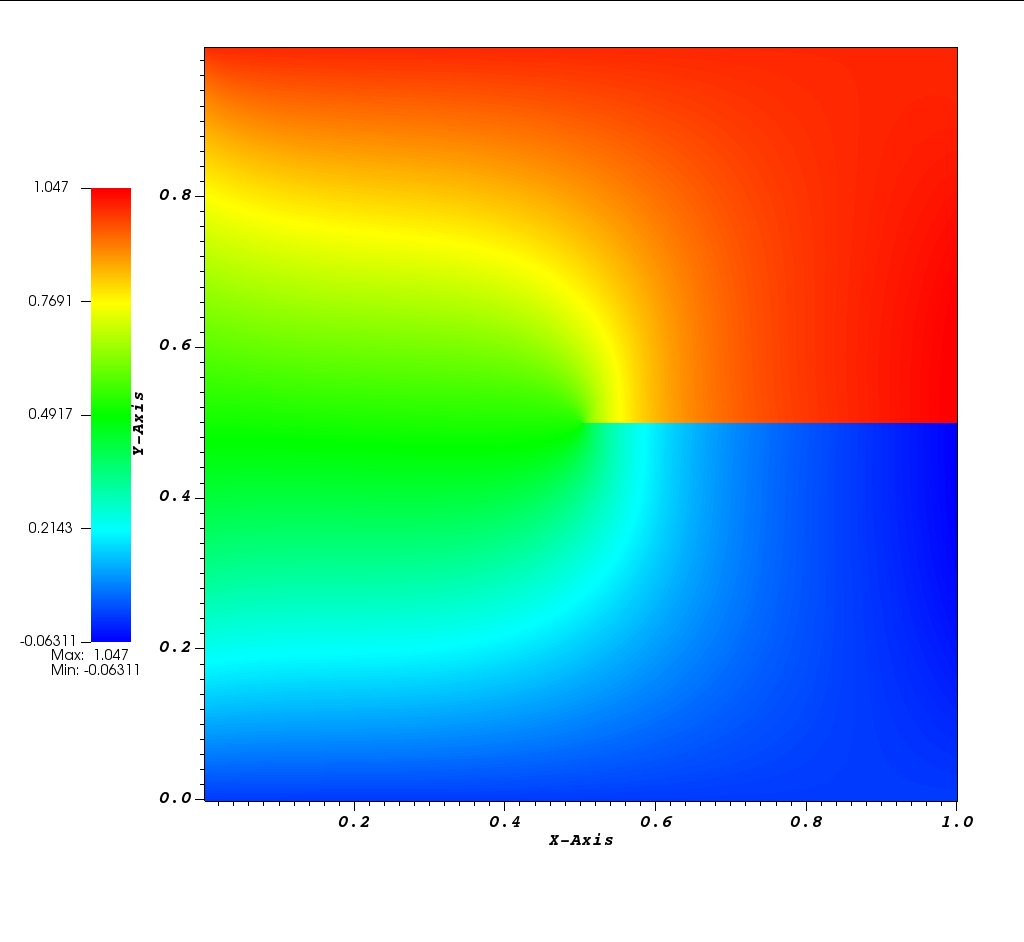}}
\caption{(Example 2 with \texttt{CASE} 1) Displacement for the Linear: (a) x-displacement and (b) y-displacement. }
\label{figs:U_Ex2_Case1_Linear}
\end{figure}

\begin{figure}[H]
\centering
\subfloat[{$\bfT_{yy}$}]{\includegraphics[width=0.45\textwidth,trim=4 4 4 4,clip]{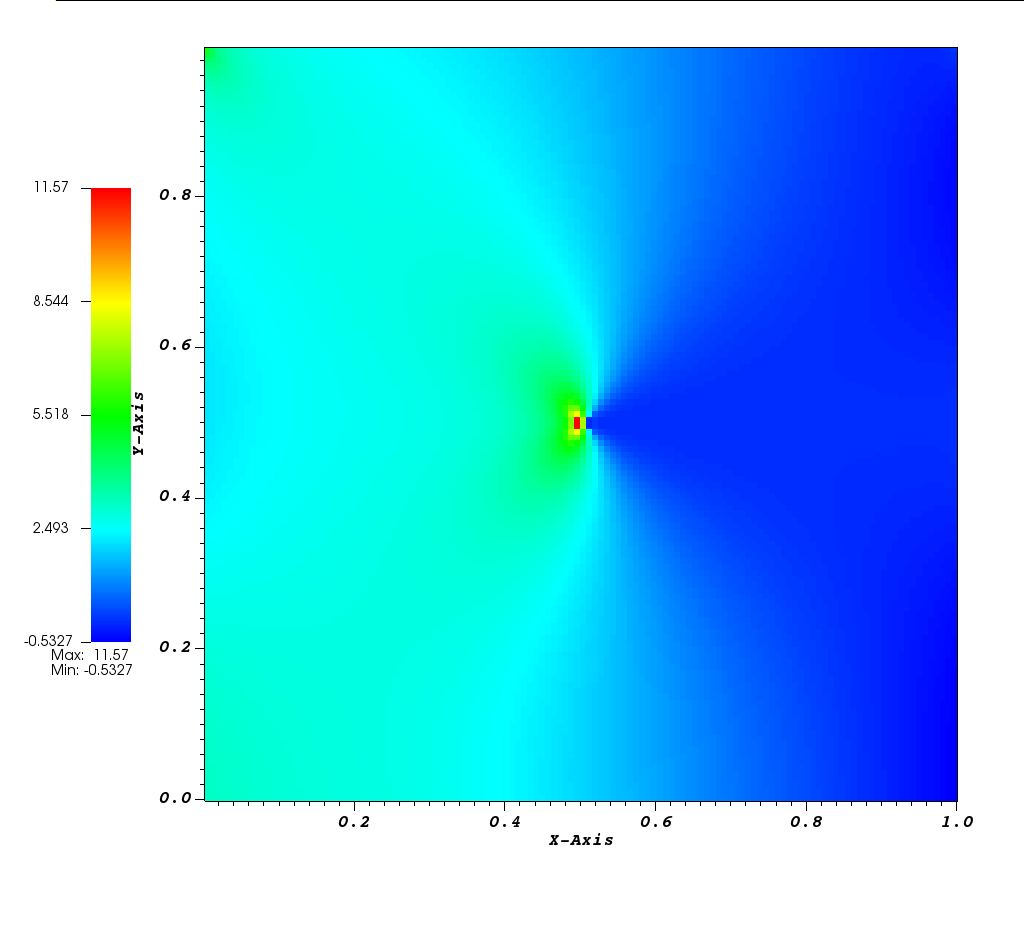}}
\hspace{0.1in}
\subfloat[{$\bfeps_{yy}$}]{\includegraphics[width=0.45\textwidth,trim=4 4 4 4,clip]{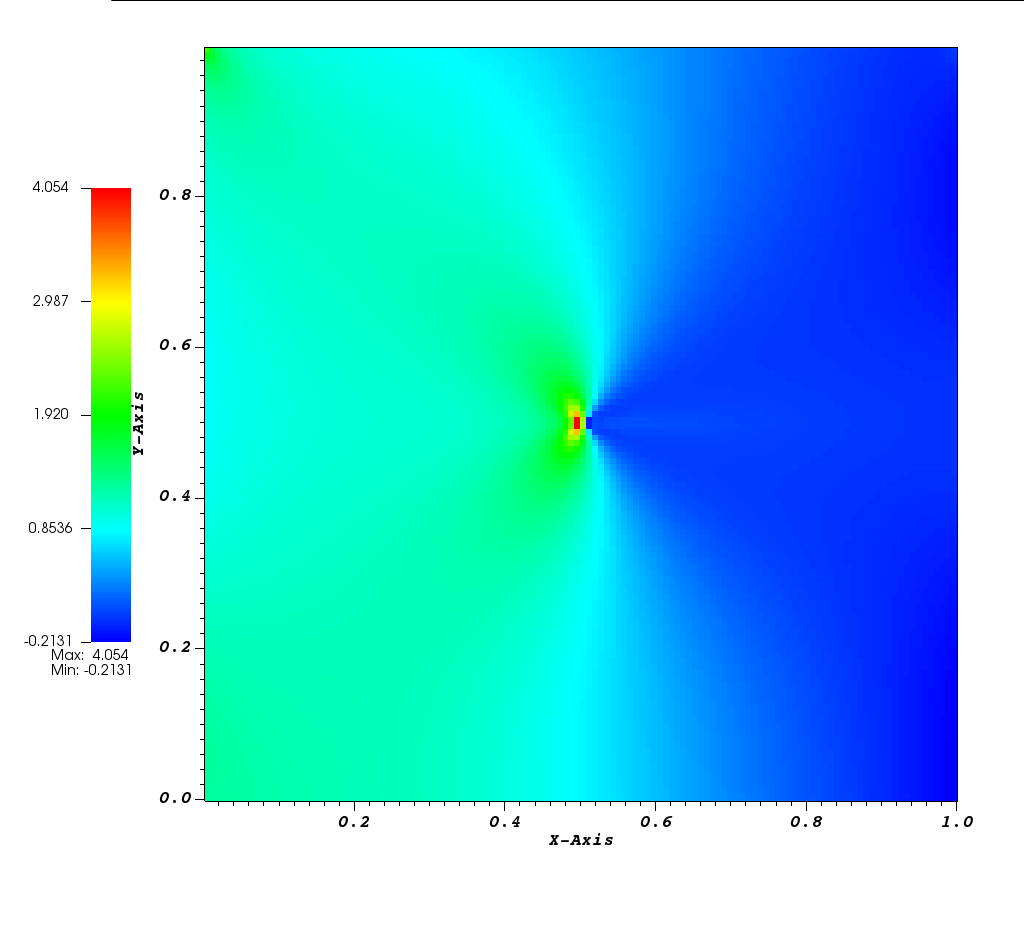}}
\caption{(Example 2 with \texttt{CASE} 1) Stress and Strain for the Linear: (a) axial stress and (b) axial strain.}
\label{figs:S_E_Ex2_Case1_Linear}
\end{figure}

{{After unveiling the  temperature solution, next we present the numerical results for the displacements for both cases of temperature boundary condition on the bottom boundary as Example 1. The ultimate goal of this research is to report some information about the nature of the solution near the crack-tip when the body is under a combination of thermo-mechanical loading.  
The numerical results will present a remarkable distinction in the crack-tip stress-strain solutions for the nonlinear strain-limiting and the linearized elasticity models.}}
\newline
\noindent(\texttt{CASE} $1$) \textbf{Linear\quad}In this case, 
a uniform and constant temperature boundary condition is applied on the bottom boundary and the classical constitutive relation {for linearized elasticity} is used to model the behavior of the bulk material. {Likewise,} we first present the displacements plot ($\bfu_x$ and $\bfu_y$ depicted in Figure \ref{figs:U_Ex2_Case1_Linear}) followed by the axial stress-strain distribution plots ($\bfT_{yy}$ and $\bfeps_{yy}$ depicted in Figure \ref{figs:S_E_Ex2_Case1_Linear}).
{Known as the classical square-root order singularity, we can see that} stress-strain values in the vicinity of the tip are large. 

\begin{figure}[H]
\centering
\subfloat[{$\bfu_x$}]{\includegraphics[width=0.45\textwidth,trim=4 4 4 4,clip]{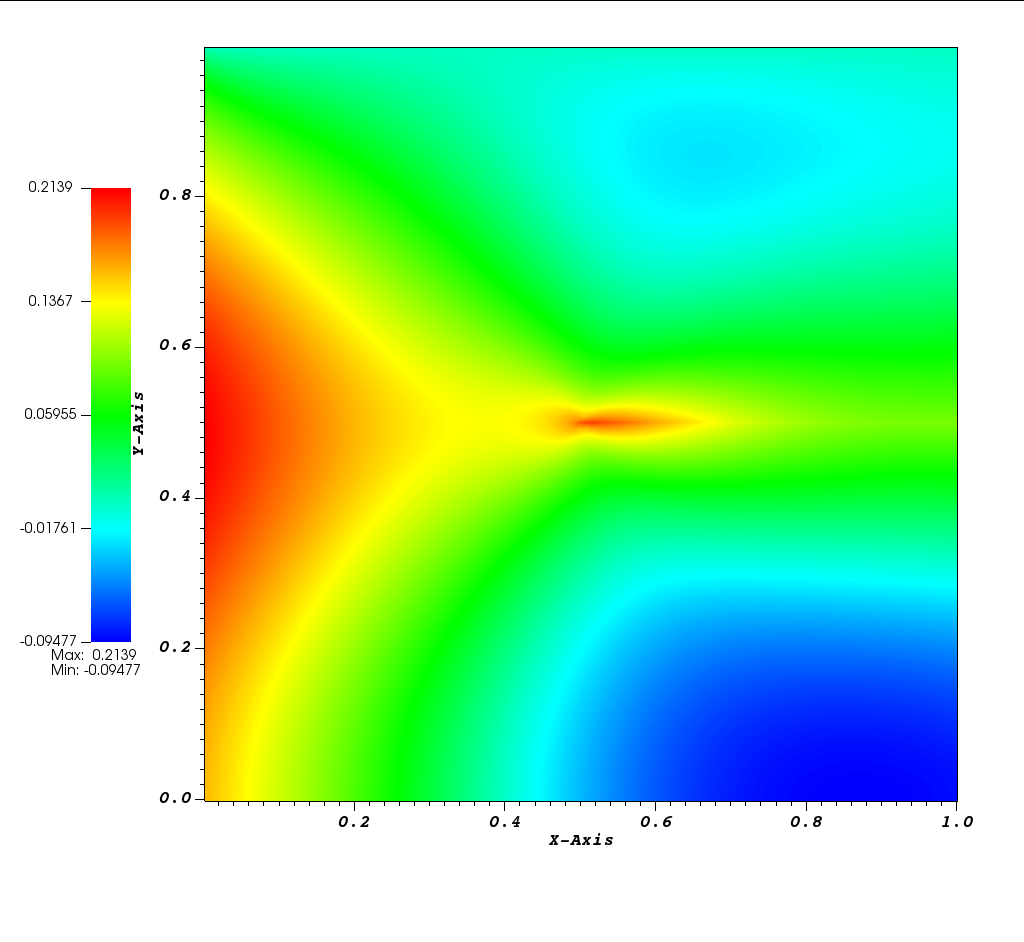}}
\hspace{0.1in}
\subfloat[{$\bfu_y$}]{\includegraphics[width=0.45\textwidth,trim=4 4 4 4,clip]{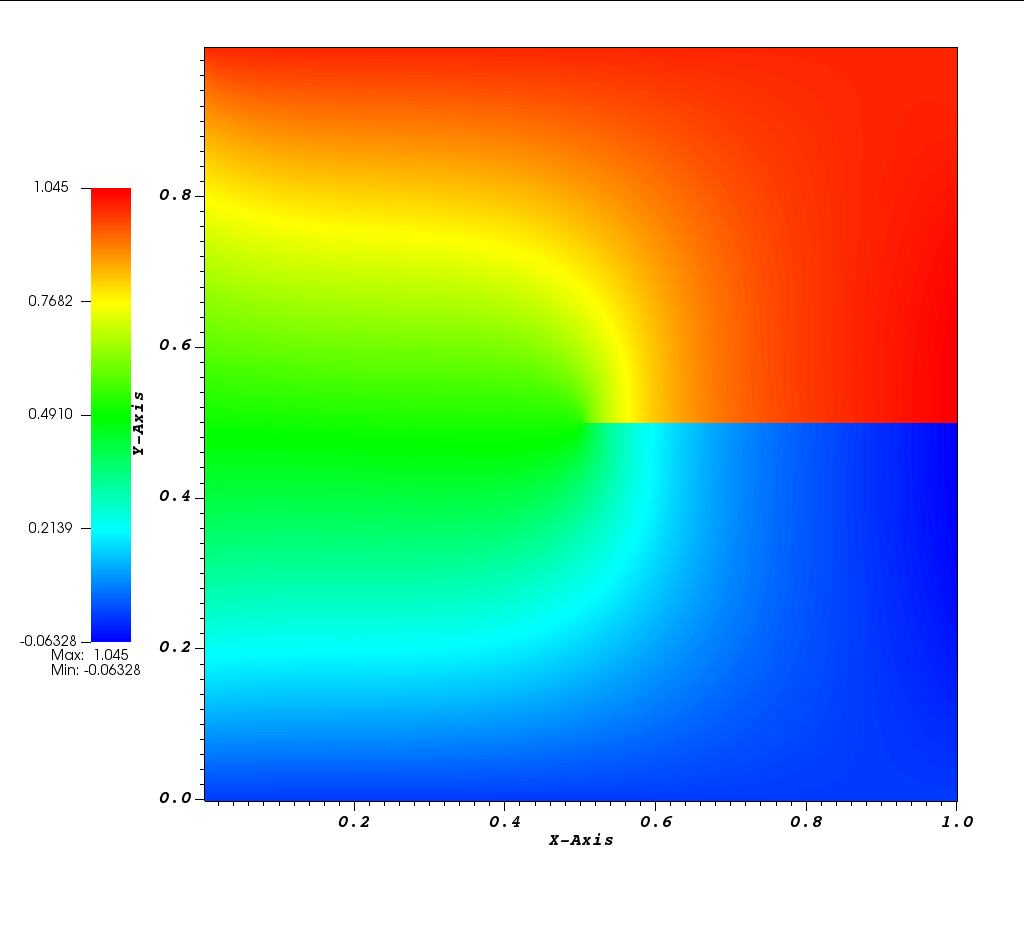}}
\caption{(Example 2 with \texttt{CASE} 1) Displacement for the Nonlinear: (a) x-displacement field and (b) y-displacement. }
\label{figs:U_Ex2_Case1_NL_A0p5_B0p01}
\end{figure}
\begin{figure}[H]
\centering
\subfloat[{$\bfT_{yy}$}]{\includegraphics[width=0.45\textwidth,trim=4 4 4 4,clip]{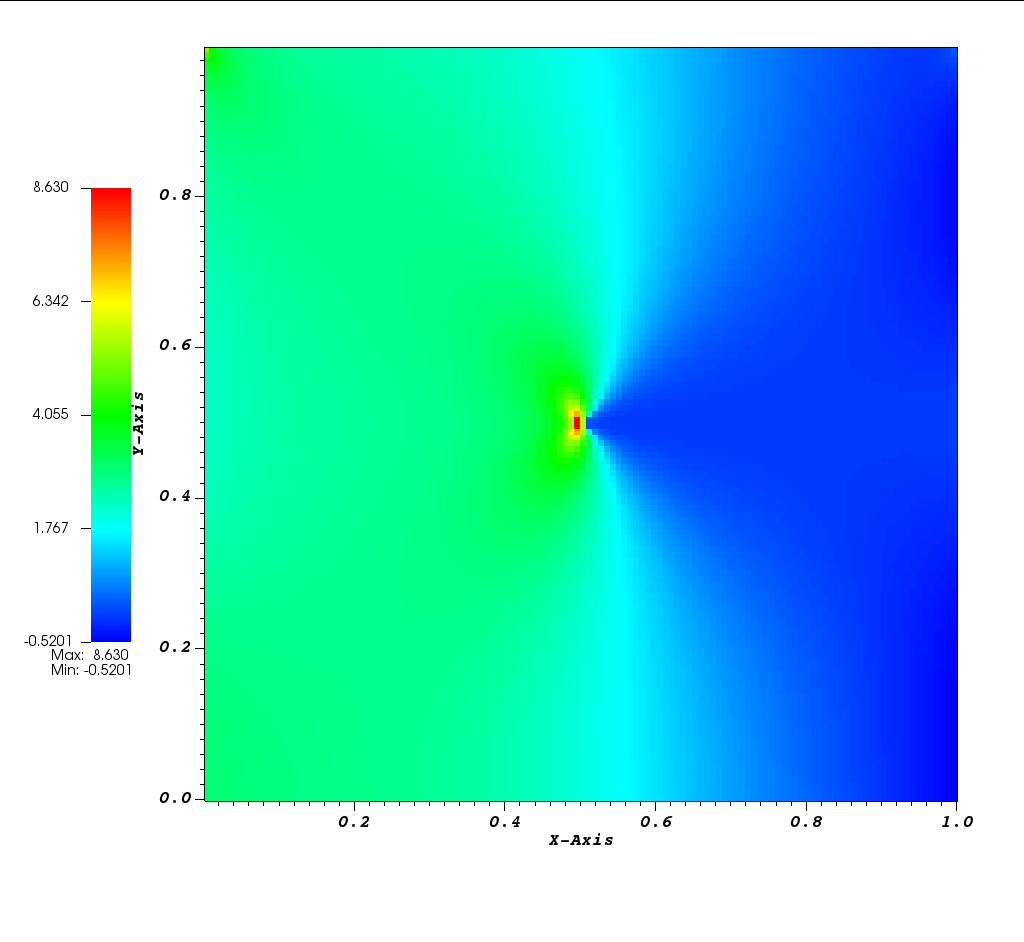}}
\hspace{0.1in}
\subfloat[{$\bfeps_{yy}$}]{\includegraphics[width=0.45\textwidth,trim=4 4 4 4,clip]{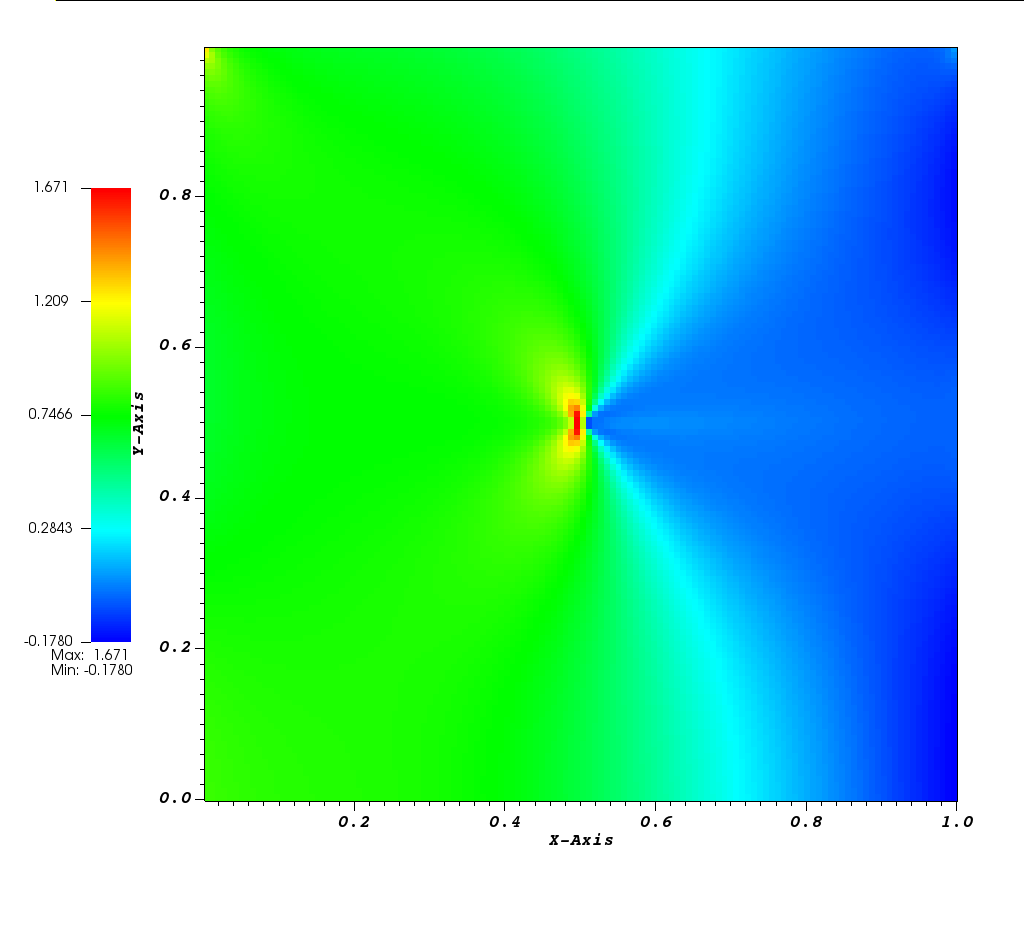}}
\caption{(Example 2 with \texttt{CASE} 1) Stress and Strain for the Nonlinear: (a) axial stress and (b) axial strain.}
\label{figs:S_E_Ex2_Case1_NL_A0p5_B0p01}
\end{figure}

\begin{figure}[H]
\centering
\includegraphics[width=1.0\textwidth]{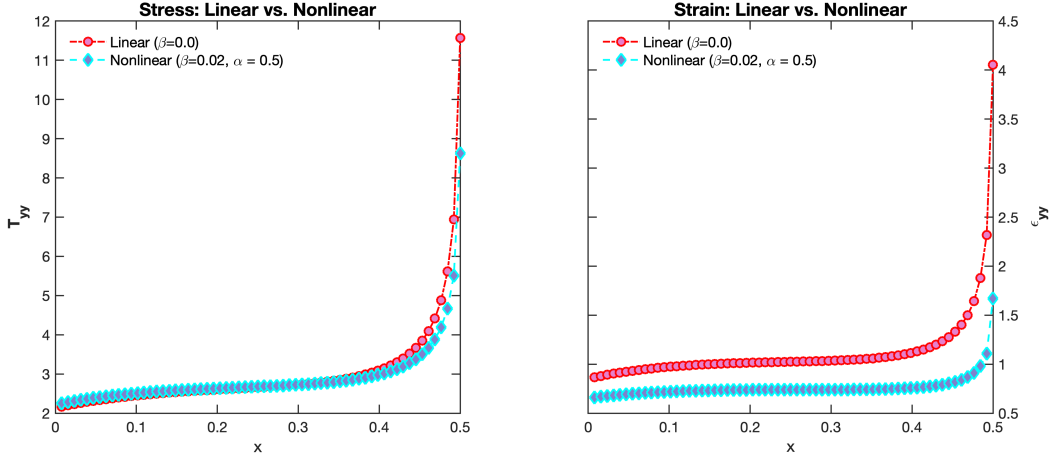}
\caption{(Example 2 with \texttt{CASE} 1) Linear vs. Nonlinear: stress {$\bfT_{yy}$} (left) and strain {$\bfeps_{yy}$} (right) on the reference line.}
\label{figs:S_E_Center_Ex2_Case1_L_NL}
\end{figure}
\begin{figure}[H]
\centering
\subfloat[{$\bfu_x$}]{\includegraphics[width=0.45\textwidth,trim=4 4 4 4,clip]{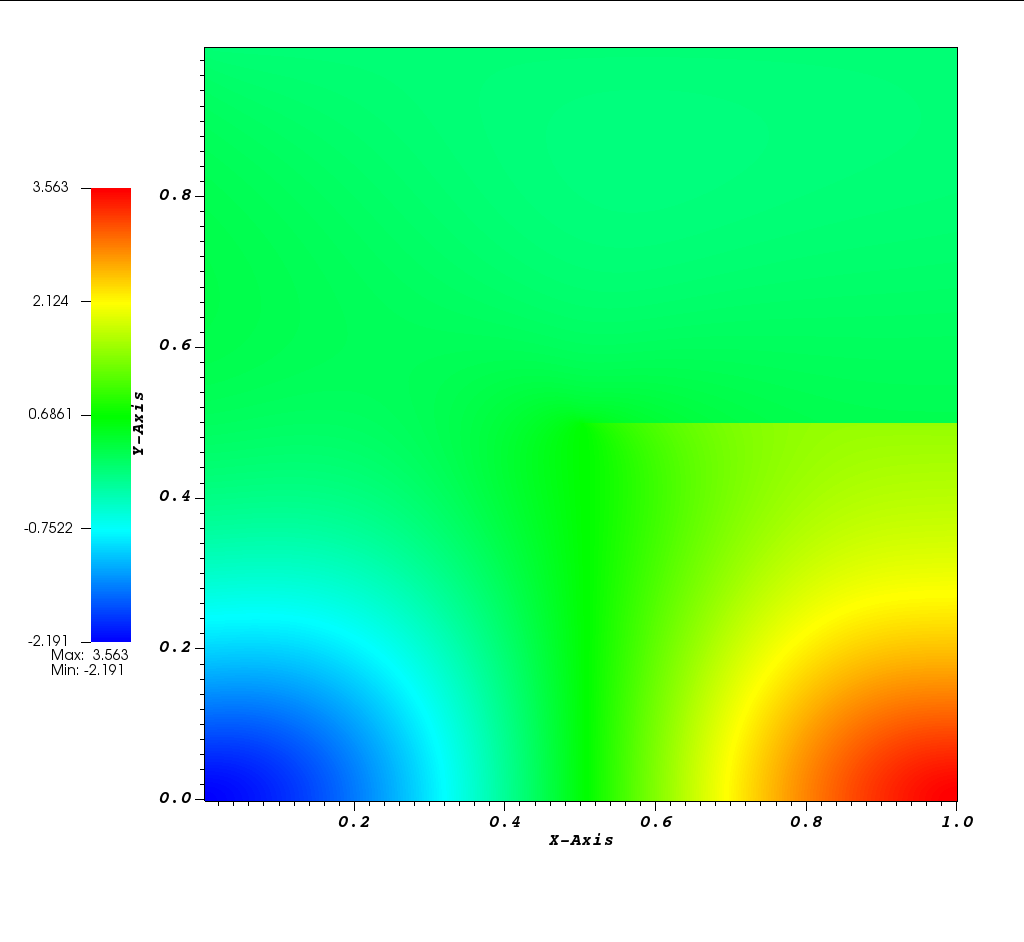}}
\hspace{0.1in}
\subfloat[{$\bfu_y$}]{\includegraphics[width=0.45\textwidth,trim=4 4 4 4,clip]{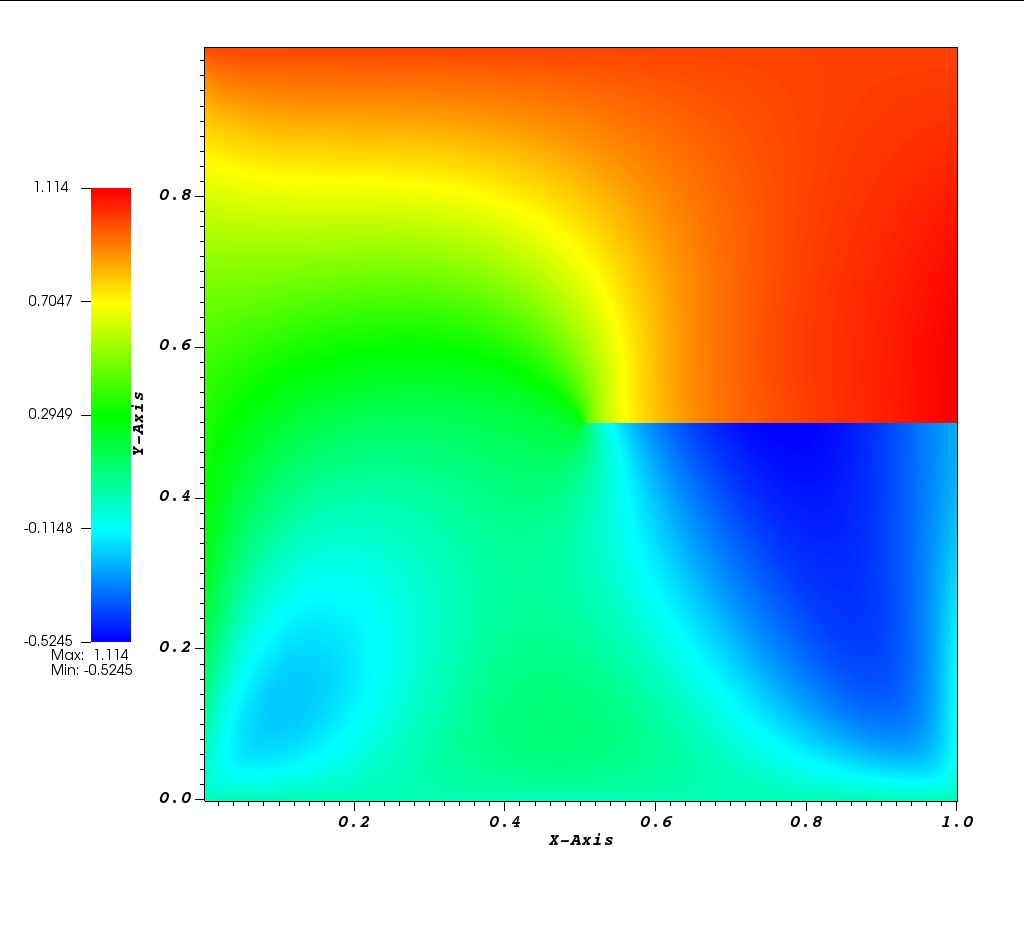}}
\caption{(Example 2 with \texttt{CASE} 2) Displacement for the Linear: (a) x-displacement and (b) y-displacement.}
\label{figs:U_Ex2_Case2_Linear}
\end{figure}
\begin{figure}[H]
\centering
\subfloat[{$\bfT_{yy}$}]{\includegraphics[width=0.45\textwidth,trim=4 4 4 4,clip]{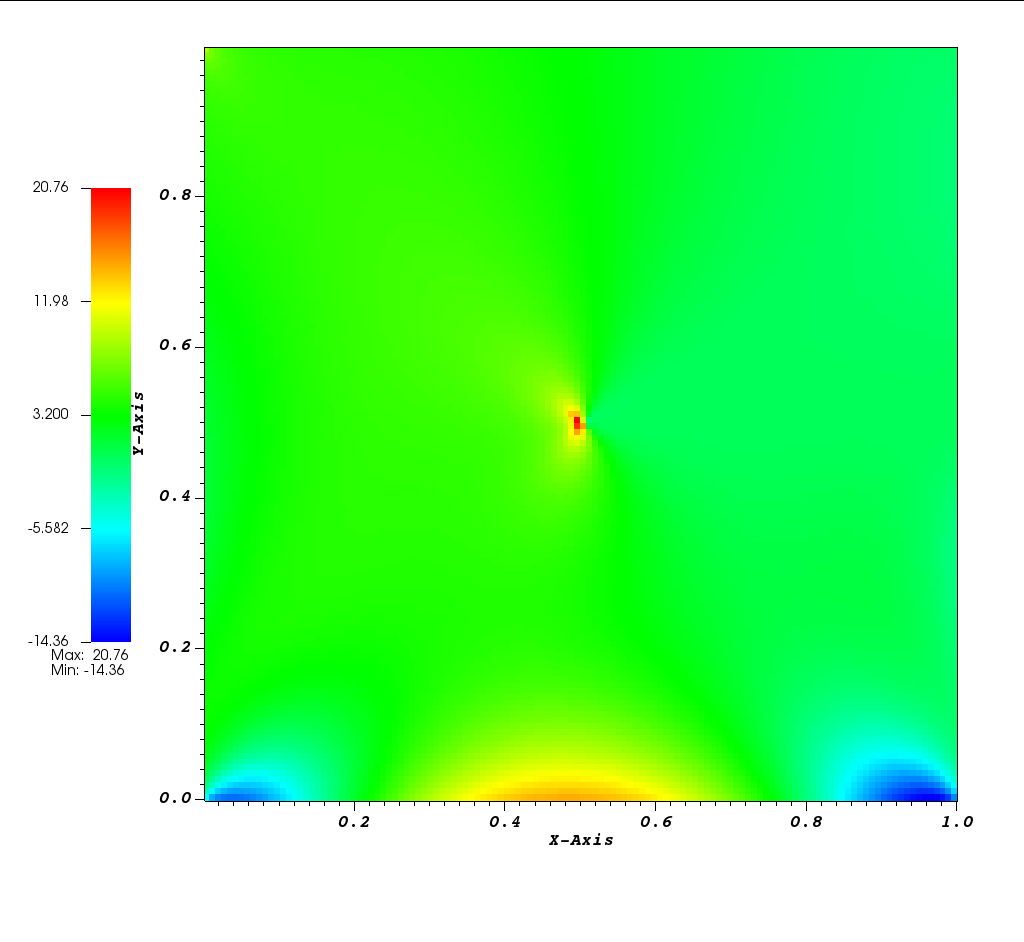}}
\hspace{0.1in}
\subfloat[{$\bfeps_{yy}$}]{\includegraphics[width=0.45\textwidth,trim=4 4 4 4,clip]{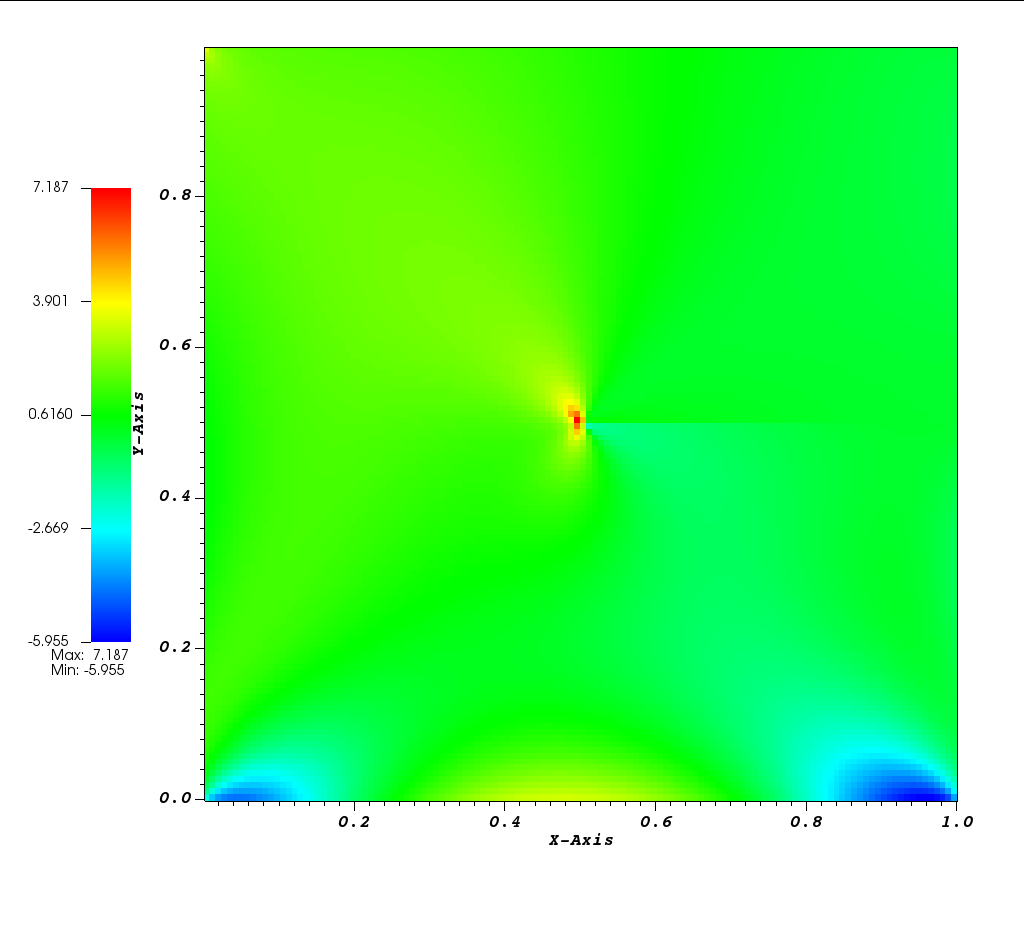}}
\caption{(Example 2 with \texttt{CASE} 2) Stress and Strain for the Linear: (a) axial stress and (b) axial strain.}
\label{figs:S_E_Ex2_Case2_Linear}
\end{figure}

\noindent(\texttt{CASE} $1$) \textbf{Nonlinear\quad}{Here,} we solve the nonlinear boundary value problem under the constant temperature condition at the bottom.  
Figure \ref{figs:U_Ex2_Case1_NL_A0p5_B0p01} depicts the displacement field $(\bfu_x,\bfu_y)$ and Figure \ref{figs:S_E_Ex2_Case1_NL_A0p5_B0p01} represents the distribution of axial stress and strain components for \texttt{CASE} 1 in Example 2. {As similar to Example 1 with \texttt{CASE} 1, the nonlinear model for Example 2 with \texttt{CASE} 1 shows not much difference in its displacement pattern and values from those of the linear model (See Figure~\ref{figs:U_Ex2_Case1_Linear} and Figure~\ref{figs:U_Ex2_Case1_NL_A0p5_B0p01}.)} {Likewise for the axial stress}, we can clearly see that the distribution seems similar to the linear model especially when  approaching the crack tip. {Figure~\ref{figs:S_E_Center_Ex2_Case1_L_NL} shows for the axial stress (left) and the axial strain (right) of the two models on the reference line approaching the tip}. For strains, we can see {much smaller} values along the radial line leading up to the crack-tip, {as expected from the nonlinear strain-limiting model}. {Also note that compared to Figure~\ref{figs:S_E_Ex2_Case1_Linear} (b),  Figure~\ref{figs:S_E_Ex2_Case1_NL_A0p5_B0p01} (b) has a relatively sharper change in color for the legend (e.g., see the cyan-colored line in between the blue and the green) over the domain and in the vicinity of the tip, implying that the the contour lines for the strain near the tip along with the shape of crack for the strain-limiting is much different from the linear model \cite{hyun-mms-2021}.} 
\newline
(\texttt{CASE} $2$) \textbf{Linear\quad}In this last example with the parabolic type of Dirichlet boundary condition for the temperature, we also consider the linear problem first. 
Figure \ref{figs:U_Ex2_Case2_Linear} depicts the displacement fields of $\bfu_x$ and $\bfu_y$, respectively. Figure \ref{figs:S_E_Ex2_Case2_Linear} depicts the axial stress and strain fields in the entire material body. It is clear from the stress-strain plots that the values of 
{stress} have the corresponding maximum ``at the tip'', and strain values are ``comparable'' with the stress concentration in the neighborhood of the crack-tip.  {Also note that the next group of large values appear near the bottom boundary $\Gamma_1$ due to the parabolic distribution of temperature, which is similar as shown in Example 1 with \texttt{CASE} $2$ (Figure~\ref{figs:S_E_Ex1_Case2_Linear}). }

\begin{figure}[H]
\centering
\subfloat[$\bfu_x$]{\includegraphics[width=0.45\textwidth,trim=4 4 4 4,clip]{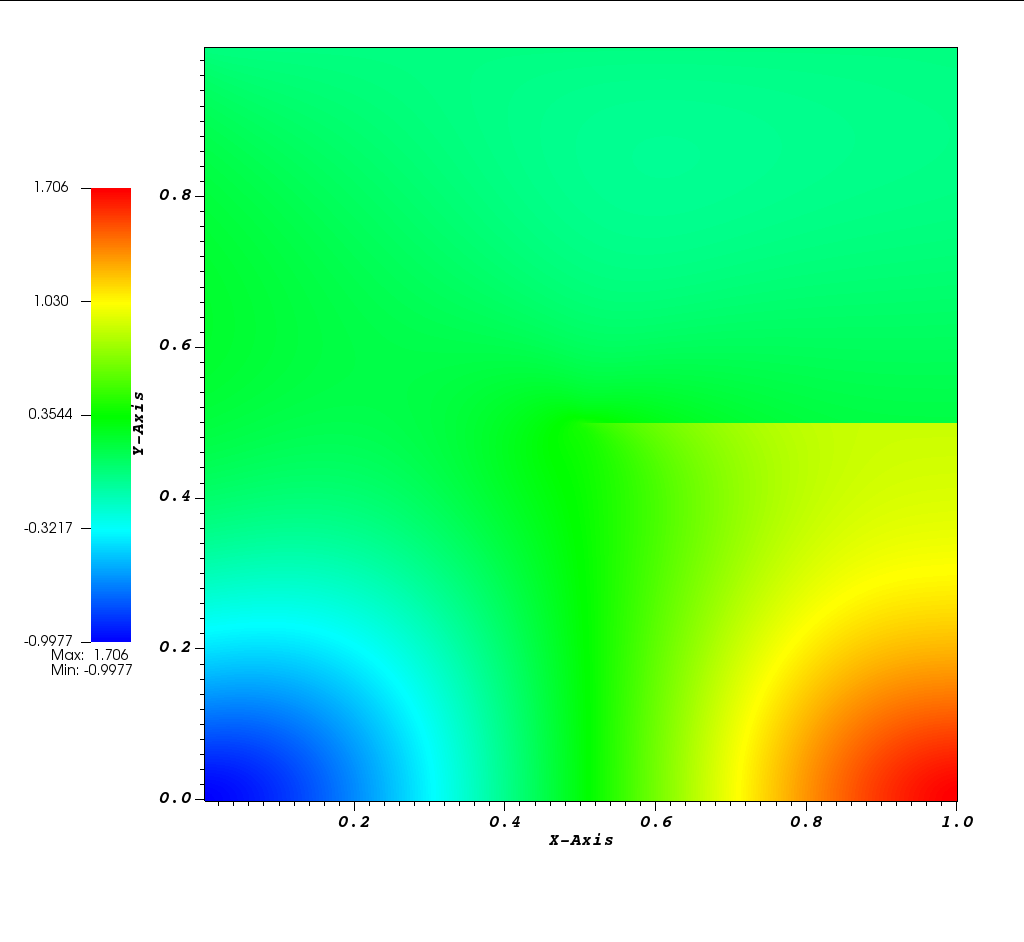}}
\hspace{0.1in}
\subfloat[$\bfu_y$]{\includegraphics[width=0.45\textwidth,trim=4 4 4 4,clip]{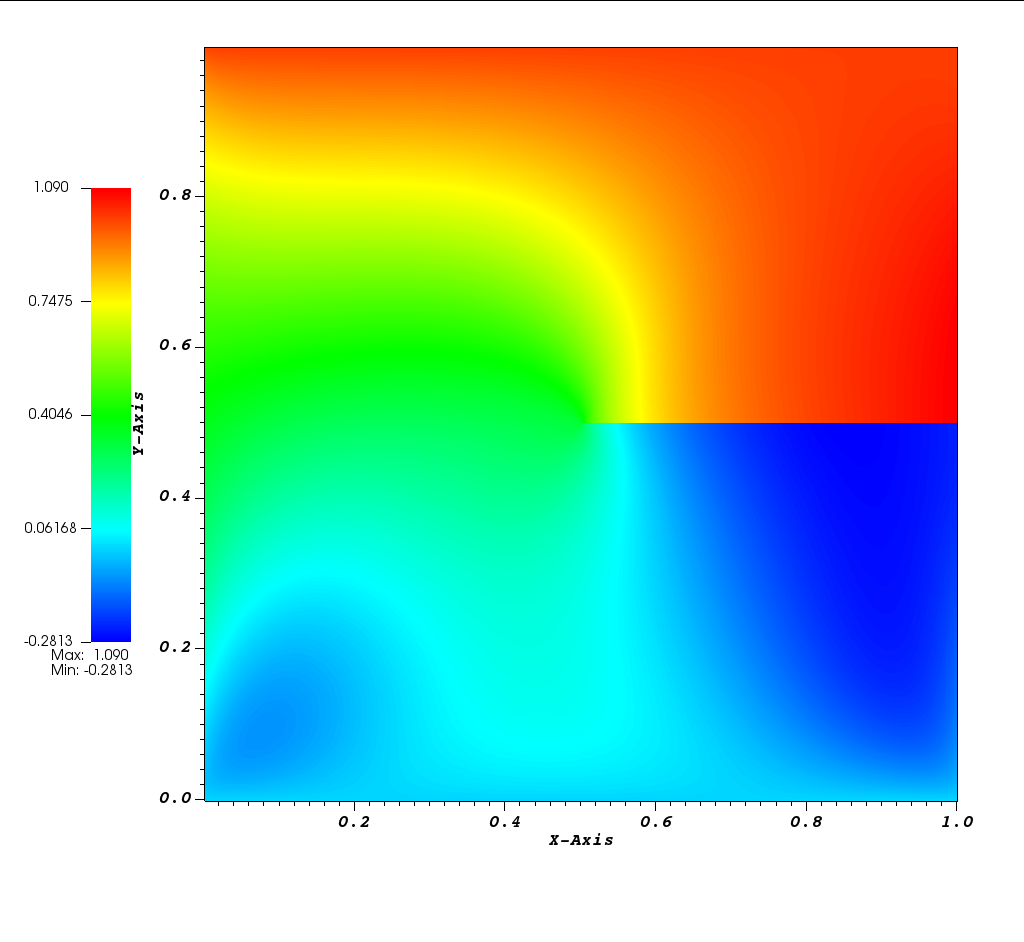}}
\caption{(Example 2 with \texttt{CASE} 2) Displacement for the Nonlinear: (a) x-displacement and (b) y-displacement.}
\label{figs:U_Ex2_Case2_NL_A0p5_B0p01}
\end{figure}

\begin{figure}[H]
\centering
\subfloat[{$\bfT_{yy}$}]{\includegraphics[width=0.45\textwidth,trim=4 4 4 4,clip]{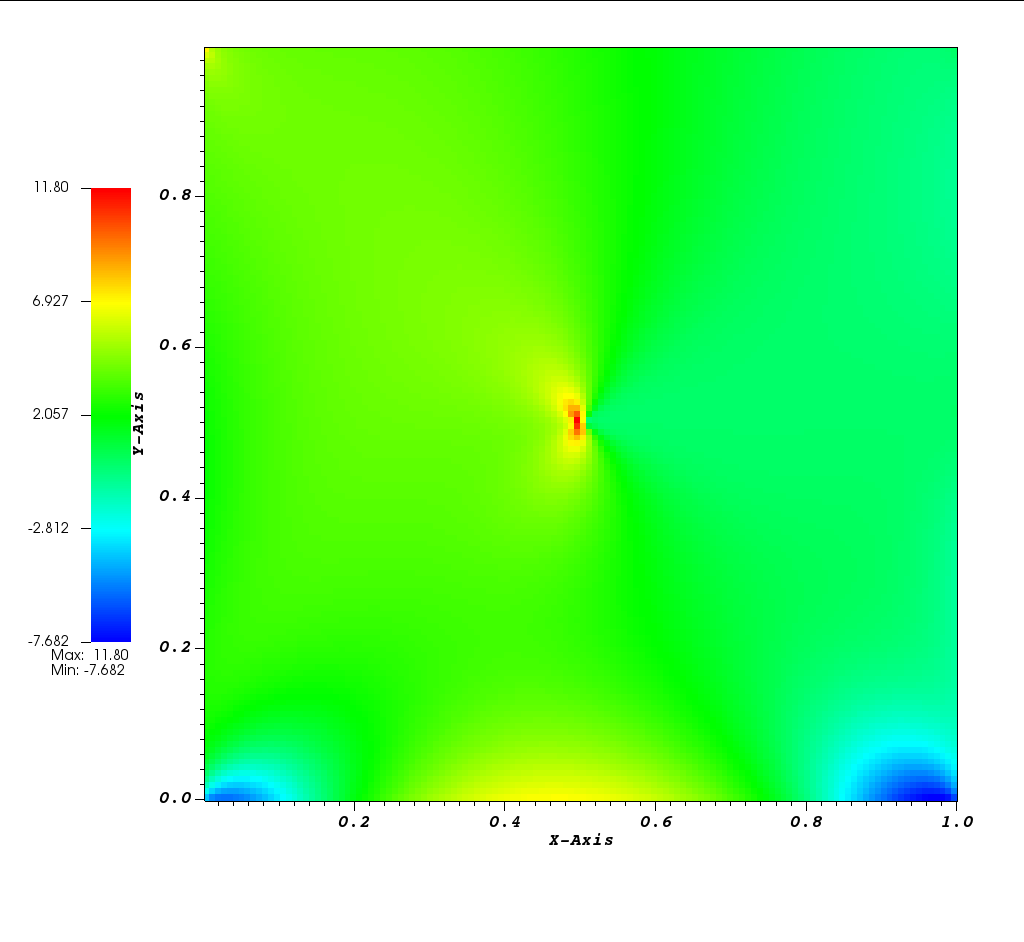}}
\hspace{0.1in}
\subfloat[{$\bfeps_{yy}$}]{\includegraphics[width=0.45\textwidth,trim=4 4 4 4,clip]{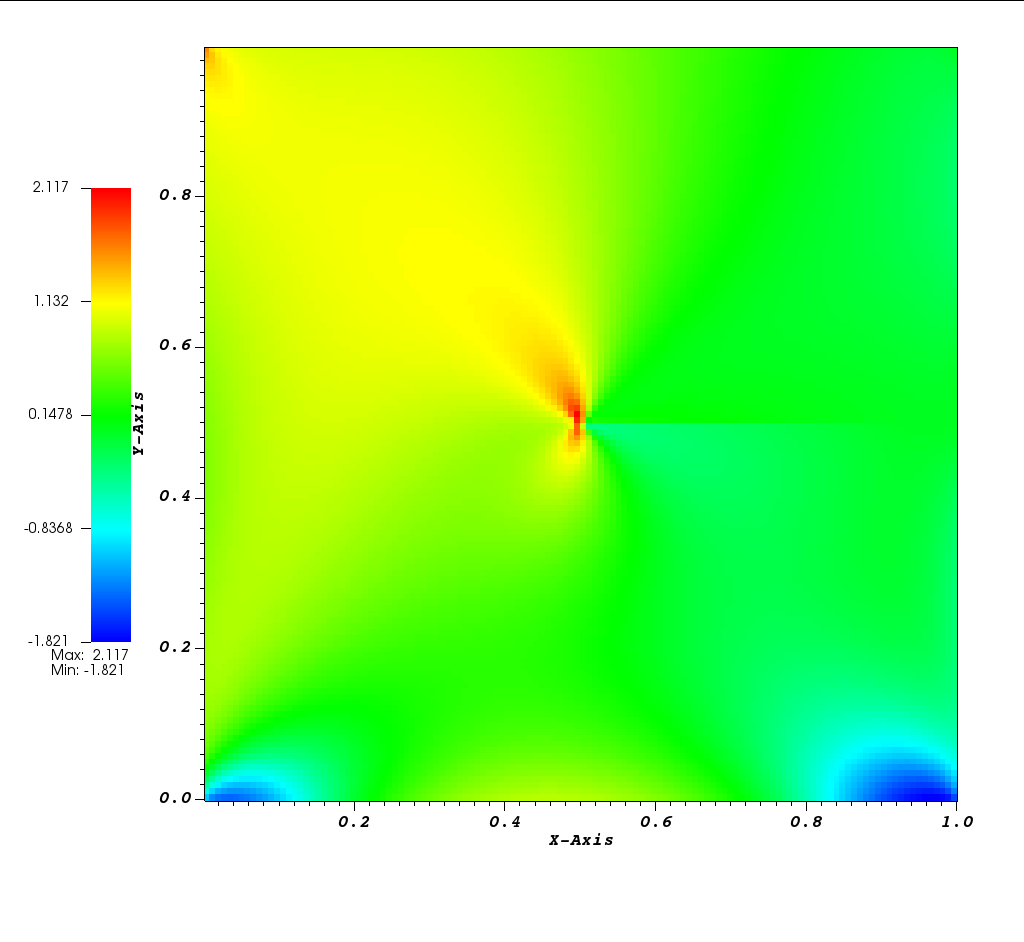}}
\caption{(Example 2 with \texttt{CASE} 2) Stress and Strain for the Nonlinear: (a) axial stress and (b) axial strain.}
\label{figs:S_E_Ex2_Case2_NL_A0p5_B0p01}
\end{figure}

\begin{figure}[H]
\centering
\includegraphics[width=1.0\textwidth]{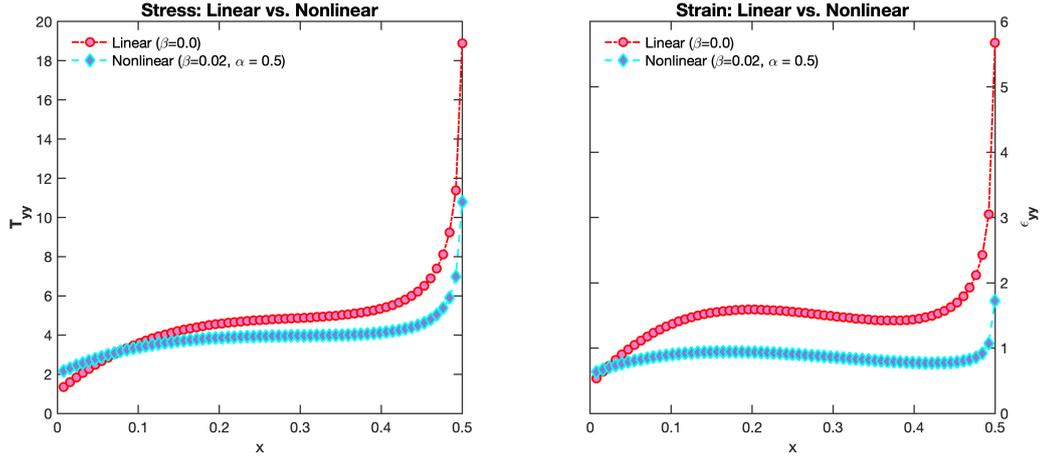}
\caption{(Example 2 with \texttt{CASE} 2) Linear vs. Nonlinear: stress {$\bfT_{yy}$} (left) and strain {$\bfeps_{yy}$} (right) on the reference line.}
\label{figs:S_E_Center_Ex2_Case2_L_NL}
\end{figure}

\noindent(\texttt{CASE} $2$) \textbf{Nonlinear\quad}Lastly, we study the nonlinear {strain-limiting} model with the boundary value problem of the parabolic temperature boundary condition. 
Figure \ref{figs:U_Ex2_Case2_NL_A0p5_B0p01}  depicts the displacement field and  Figure \ref{figs:S_E_Ex2_Case2_NL_A0p5_B0p01} delineates the stress-strain for the nonlinear model in the whole material body. {The differences in their values for displacement, stress and strain between the linear and the nonlinear models get larger for \texttt{CASE} $2$ than the corresponding differences for \texttt{CASE} $1$.} Figure \ref{figs:S_E_Ex2_Case2_NL_A0p5_B0p01} (a) clearly characterizes {the concentration of the stress} in the neighborhood of the crack-tip, {whereas Figure \ref{figs:S_E_Ex2_Case2_NL_A0p5_B0p01} (b) represents the bounded strain values which is more distributed over the domain compared to the one for the linear model (Figure~\ref{figs:S_E_Ex2_Case2_Linear} (b)). Also from its maximum and  minimum values, it can be inferred that the growth of strain is not in the same order as the one for stress but limited.}
{The contrast are clearly shown in Figure \ref{figs:S_E_Center_Ex2_Case2_L_NL} depicting the axial stress (left) and the axial strain (right) on the reference line approaching the tip with the two models.} 
{We confirm that} the growth rate at which stresses are increasing near the crack-tip is much larger than than rate at which the strain does, and there is a clear distinction of the strain behavior between the two models studied in this paper. The stress growth in the neighborhood suggests that the crack-tip is a singular energy sink, therefore one can make use of the local growth criterion to study the crack evolution within the framework of {strain-limiting model}  
proposed in this paper.  
{Likewise as seen in Example 1,  
\texttt{CASE} 2 has more differences between the linear and the nonlinear models than \texttt{CASE} 1 (See Figure~\ref{figs:S_E_Center_Ex2_Case1_L_NL} and Figure~\ref{figs:S_E_Center_Ex2_Case2_L_NL}).}

\section{Conclusion}
In the present paper, we have studied a finite element discretization of  some boundary value problems formulated within a relatively new class of strain-limiting nonlinear constitutive relations. The present study is an extension of constitutive theory developed in \cite{rajagopal2007elasticity} to the generalized thermoelasticity. We formulate a well-posed mathematical model to describe the response state of the elastic body under thermo-mechanical loading. {Since} the model consists of a nonlinear relationship between stress and strain, numerical discretization is not straightforward. A key {approach in its formulation} 
is the linearization using \textit{Newton's method}, therefore it require us to solve one elliptic boundary value {problem} at each Newton's step.  Numerical 
experiments provide a solid evidence in support for the efficacy of the proposed model. We observe the crack-tip stress concentration, similar to the classical model, but the growth of strain does not 
{behave similarly} to that of stress {but as expected with a slower order for growth near the tip}. {As a part of the} future work, 
{an immediate next step can include extending} the model studied in this paper for the quasi-static evolution of network of cracks {under thermo-mechanical loading}. {Furthermore, the formulation studied in this article can also be extended to study  
{a reliable constitutive relation for the heat such as the} wave propagation without energy dissipation for the Fourier's law.}

\section*{Acknowledgements}
Authors would like to thank the support of College of Science \& Engineering, Texas A\&M University-Corpus Christi. One of the author, HCY, appreciates the support received in part by the basic research fund (Project No. GP2020-006) of Korea Institute of Geoscience and Mineral Resources (KIGAM).

\bibliographystyle{model3-num-names}
\bibliography{NonLinear_thermo_paper}

\end{document}